\documentclass[12pt,fleqn]{report}
\usepackage{amsmath}
\usepackage{amsthm}
\usepackage{amsfonts}
\usepackage{epsfig}
\usepackage{psfig}
\usepackage{indentfirst}
\usepackage{multirow}
\usepackage{graphicx}
\usepackage{psfrag}

\theoremstyle{plain}
\newtheorem{Thm}{Theorem}[section]
\newtheorem{Cor}[Thm]{Corollary}
\newtheorem{Lm}[Thm]{Lemma}
\newtheorem{Prop}[Thm]{Proposition}
\newtheorem{Prob}[Thm]{Problem}
\newtheorem{Def}[Thm]{Definition}

\newtheorem{Rm}[Thm]{Remark}
\newtheorem{Example}[Thm]{Example}

\theoremstyle{definition}

\def\cal{\mathcal}
\def\bdef{\begin{Def}}
\def\endef{\end{Def}}
\def\bthm{\begin{Thm}}
\def\ethm{\end{Thm}}
\def\blm{\begin{Lm}}
\def\elm{\end{Lm}}
\def\bprop{\begin{Prop}}
\def\eprop{\end{Prop}}
\def\bcor{\begin{Cor}}
\def\ecor{\end{Cor}}
\def\brm{\begin{Rm}}
\def\erm{\end{Rm}}
\def\bprob{\begin{Prob}}
\def\eprob{\end{Prob}}
\def\bex{\begin{Example}}
\def\enex{\end{Example}}
\def\beal{\begin{aligned}}
\def\enal{\end{aligned}}

\def\beq{\begin{eqnarray}}
\def\eneq{\end{eqnarray}}
\def\Cal{\mathcal}
\def\pa{\partial}
\def\R{\mathbb R}
\def\C{\mathbb C}
\def\Z{\mathbb Z}
\def\om{\omega}
\def\~{\tilde}

\def\Bbb{\mathbb}
\def\R{\Bbb R}
\def\C{\Bbb C}
\def\Z{\Bbb Z}
\def\V{\Cal V}

\def\Cal{\mathcal}
\def\Si{\Sigma}

\def\asik{Khovanski}
\def\W{\Omega}
\def\land{\wedge}
\def\ell{\mathsf p}
\def\F{\cal F}
\def\o{\emptyset}
\def\A{\Cal A}

\def\J{\Cal J}
\def\1{\bf{1}}

\def\D{\Cal D}

\def\dt{\delta}

\def\al{\alpha}
\def\gm{\gamma}
\def\lb{\lambda}

\def\dt{\delta}
\def\eps{\varepsilon}

\def\inv{^{-1}}

\def\pa{\partial}

\begin{document}

\title{Around Hilbert-Arnold Problem}
\author{Vadim Kaloshin}
\date{}
\maketitle

\begin{abstract}
This lectures notes consists of four lectures. The first
lecture discusses questions around Hilbert-Arnold Problem
which is naturally arises from Quantitative Hilbert 16-th
problem. In the second lecture we outline author's solution
\cite{K1} of a weak form of  Local Hilbert-Arnold Problem.
This solution provides an independent proof of
Ilyashenko-Yakovenko Finiteness Theorem \cite{IY2}.
The third lecture discusses question of existence of
$a_P$-stratification of Thom \cite{T2} and presents a simple
geometric proof of existence of such a stratification for
polynomial functions, which was originally proven by
Hironaka \cite{Hi}. The forth lecture gives application of
Grigoriev-Yakovenko's construction to the problem of growth
of the number of periodic points and the problem of
bifurcation of spacial polycycles. The latter problem
naturally generalizes Local Hilbert-Arnold Problem.
\end{abstract}

\pagenumbering{roman}



\newpage
\addcontentsline{toc}{chapter}{List of Figures}

\newpage

\newpage
\pagenumbering{arabic}

\chapter{Around the Hilbert 16-th problem and an 
estimate for cyclicity of elementary polycycles.}
\label{introd}
\thanks{The author is partially supported by the Sloan Dissertation
Fellowship and the American Institute of Mathematics Five-Year 
Fellowship}

\medskip
 
\section{The Hilbert 16-th Problem and its offsprings}

Consider a polynomial vector field on the real $(x,y)$-plane
\begin{eqnarray} \label{H}
\left\{
\begin{aligned}
\dot y=P_n(x,y) \\
\dot x=Q_n(x,y)
\end{aligned}
\right.
\quad P_n,Q_n - {\textup{ polynomials}},\ \deg P_n,Q_n \leq  n.
\end{eqnarray}
A limit cycle of a polynomial vector field (\ref{H}) is an isolated 
periodic solution. Define 
\beq
H(n)= \boxed {\text{uniform bound for the number of limit 
cycle of}\ (\ref{H}).} \nonumber
\eneq

One way to formulate the Hilbert 16-th Problem is the following:

{\bf{The Hilbert 16-th Problem (HP).}} {\it{Estimate $H(n)$
for any $n \in \Z_+$.}}

To prove that $H(1)=0$ is an exercise, but find $H(2)$ is already
a difficult unsolved problem ( see \cite{DRR}, \cite{DMR} for
work in this direction). Below we discuss two of the most significant 
branches of research HP generated: Existential and Tangential 
Hilbert 16-th Problems. 

\subsection{The Tangential Hilbert 16-th Problem}
Consider a polynomials perturbation of a Hamiltonian polynomial
line field
\beq \label{TH}
\left\{
\begin{aligned}
\dot y=-\frac{\partial H}{\partial x}-\eps Q(x,y)\\
\dot y=\frac{\partial H}{\partial y}+\eps P(x,y).
\end{aligned} \right.
\eneq
For $\eps=0$ the line field (\ref{TH}) does not have any limit 
cycles at all (all cycles are nonisolated). An oval (topological 
circle) $\gm$ of the level curve $H(x,y)=h$ generates a limit 
cycle for small nonzero value of $\eps$ if the accumulated energy 
dissipation is zero in the first approximation, i.e. when 
\beq \label{abel}
\oint\ P(x,y)\ dx +\ Q(x,y)\ dy=0, \quad \gm \subseteq 
\{H(x,y)=h\}.
\eneq
The left-hand side expression is called {\it a complete Abelian
integral}. If polynomials $H,P,$ and $Q$ are fixed the integral
(\ref{abel}) defines a multivalued function $I(h)$. Multivaluity
appears when the corresponding level curve  $\{H(x,y)=h\}$ has 
several disjoint ovals.

{\bf Tangential Hilbert 16-th Problem (THP).} {\it \cite{A1}
 For any collection of polynomials $H,P,$ and $Q\in \R[x,y]$ 
of degree $\leq n$ give an upper bound $TH(n)$ on the number of real 
ovals over which the integral (\ref{abel}) vanishes, but not 
identically}.

In the latter case the perturbation (\ref{TH}) is a Hamiltonian 
system for $\eps\neq 0$ so it does not have limit cycles at all.
Even though Tangential Hilbert Problem is not solved yet, in contrast
to the Hilbert 16-th, there are several quite general results related
to it. Khovanski \cite{Kh1}--Varchenko \cite{V} proved 
{\it Finiteness Theorem: for any $n\in \Z_+$ the number of isolated 
zeroes of Abelian integrals is uniformly bounded over all Hamiltonian 
and forms of degree $\leq n$}.  

For various other results estimating $H(n)$ in various particular 
cases see \cite{Ga}, \cite{I1}, \cite{Mr}, \cite{NY1}, 
\cite{P}, and the lecture course \cite{NY2} in the present volume 
for more references.

If we consider (\ref{abel}) over the field of complex we have that 
an Abelian integral satisfies a fuchsian equation or Picard-Fuchs
equation (see e.g. \cite{AA}), i.e. an equation of the form
\beq \label{fuchs}
\dot z=\sum_j \frac{A_j}{t-\alpha_j} z_j, \quad
\textup{where the}\ A_j\ \textup{are constant matrices}  
\eneq
and $z=(z_1,\dots,z_p)\in \C^p$ is a complex vector for some $p$.
Investigation of various properties of fuchsian equations is the main 
topic of lectures of Bolibrukh \cite{Bo} in the present volume.

\subsection{From Existential Hil\-bert 16-th  Problem to
Hilbert-Ar\-nold Problem}

A qualitative form of  Hilbert 16-th Problem is the following:

{\bf{Existential Hilbert 16-th Problem (EHP).}}
{\it{Prove that $H(n)< \infty$ for any $n \in \Z_+$.}}

The problem about finiteness of number of limit cycles for an 
individual polynomial line field (\ref{H}) is called {\it{Dulac
problem}} after the pioneering work of Dulac \cite{Du} who claimed 
in 1923 to solve this problem, but an error was found by Ilyashenko 
\cite{I2} 60 years later. The Dulac problem was solved by two 
independent, rather different and incredibly complicated proofs 
given almost simultaneously by Ilyashenko and Ecalle 

{\bf Individual Finiteness Theorem (IFT).} {\it {\cite{I3}}, 
{\cite{E}}\ Any polynomial line field  (\ref{H}) has only a finite
number of limit cycles.}
  
However, neither proof allows any generalization to solve EHP.
Consider the equation (\ref{H}) for different polynomials
$(P_n(x,y),Q_n(x,y))\in \R^2$ as the family of vector fields on 
$\R^2$ depending on parameters of the polynomials. Using a central 
projection $\pi:\Bbb S^2 \to \R^2$ and homogeneity with respect 
to parameters of the equation (\ref{H}) (vector fields 
$(\lb P_n(x,y), \lb Q_n(x,y))$ and $(P_n(x,y),Q_n(x,y))$ for any 
$\lb \neq 0$ have the same trajectories) one can construct a 
{\it{finite parameter family of analytic line fields on the sphere 
$\Bbb S^2$ with a compact parameter base $B$ (see e.g. {\cite{IY2}} 
for details).}} After this reduction Existential Hilbert Problem
becomes a particular case of the following 

{\bf{Global Finiteness Conjecture (GFC).}}{\it{ \cite{R1} For 
any family of line fields on $\Bbb S^2$ with a compact parameter 
base $B$ the number of limit cycles is uniformly bounded over all 
parameter values.}}

We refer the reader to the volumes {\cite{IY2}}, \cite{S}, and a book
\cite{R2}, where various development of these and related problems 
are discussed. Families of analytic fields are difficult to analyze.  
In the middle of 80's Arnold {\cite{AA}} proposed to 
consider generic families of smooth vector fields as the first step
toward understanding families of analytic vector fields.
A smooth analog of Global Finiteness Conjecture is the following

{\bf{Hilbert-Arnold Problem (HAP).}} {\it{{\cite{I4}}
Prove that in a generic finite parameter of vector fields on 
the sphere $S^2$ with a compact base $B$, the number of limit cycles 
is uniformly bounded. }} 

Assume for a moment that an analytic (or a generic smooth) 
vector field on the sphere $\Bbb S^2$ has an infinite number of limit
cycles. By the Poincar\'e-Bendixon Theorem, any limit cycle should 
surround an equilibrium point and, since our vector field has at most
finitely many equilibria, there should be an infinite ``nested'' 
sequence around one of equilibria. Then this ``nested'' sequence of
limit cycles have to accumulate (in the sense of Hausdorff metric) to
a certain contour (polygon) consisting of equilibria (as vertices) 
and separatric curves (sides of that polygon) connecting them.
Such objects are called {\it polycycles}. 
It turns out that a possible solution to Hilbert-Arnold Problem 
reduces to investigation of bifurcation of polycycles. 
Let's give several definitions.

\bdef \label{polyc}
A polycycle $\gm$ of a vector field on the sphere $\Bbb S^2$
is a cyclically ordered collection of equilibrium points 
$p_1, \dots, p_k$ (with possible repetitions) and arcs
$\gm_1, \dots , \gm_k$ (distinct integral curves consisting possibly 
from equilibrium points) connecting them in the specific order: the 
j-th arc $\gm_j$ connects $p_j$ with $p_{j+1}$ for $j=1,\dots, k$.  

A polycycle $\gm$ is called monodromic if one can choose a segment
$\Sigma$ transversal to $\gm$ such that one side $U \subset \Sigma$ 
of $p=\Sigma\cap \gm$ a Poincare return map $\Delta_\gm:U \to \Sigma$
is defined with $\Delta_\gm(p)=p$. 
\endef

{\bf Nonaccumulation Theorem.} {\it  \cite{I3}, \cite{E}
For any analytic monodromic polycycle $\gm$ there is a tube 
neighborhood free from limit cycles or a Poincare return map 
$\Delta_\gm:\Sigma\supset U\to \Sigma$ can't have infinitely many 
fixed points accumulating to $p=\Sigma\cap \gm$}.

This Theorem along with above compactness arguments implies
IFT. Both proofs of Ilyashenko and Ecall\'e deal with analysis of
type of germs of maps arising as Poincare return maps of an analytic
monodromic polycycle. Lectures by van den Dries \cite{Dr1}
(see also \cite{Dr2}) in the present volume discusses the theory of 
o-minimality. This theory deals with classes of functions which 
satisfy certain axioms. A basic example of o-minimal class of
functions is polynomials and analytic functions. In particular,
if a map of a compact interval $\Delta:U\to \R$ belongs to 
an o-minimal class of functions, then the equation $\Delta(x)=x$ has 
a finitely many solutions. One of hopes is that deeper understanding 
of o-minimal structures would allow to include Poincare return maps
of monodromic polycycles into an o-minimal class and give an
independent proof of Nonaccumulation Theorem.
Finiteness Theorems for differentiable function fields are
discussed in lectures by Buium \cite{Bu} in the present volume.

\bdef \label{cycle}
Let $\{\dot x=v(x,\eps)\}_{\eps \in B^n}, \ x \in \Bbb S^2,$ be an
$n$-parameter family of vector fields on $\Bbb S^2$
having a polycycle $\gm$ for some parameter value $\eps_*\in B^n$.
The polycycle $\gm$ has cyclicity $\mu$ in the family 
$\{v(x,\eps)\}_{\eps \in B^n}$ if there exist neighborhoods 
$U$ and $V$ such that $\Bbb S^2 \supseteq U \supset \gm,
\ B \supseteq V \in \eps_*$ and for any $\eps\in V$ the field 
$v(\cdot,\eps)$ has no more than 
$\mu$ limit cycles inside $U$ and $\mu$ is the minimal number
with this property.
\endef

{\it{ Examples:}} \ \ 
1) In a generic $n$-parameter family, the maximal multiplicity of 
a degenerate limit cycle does not exceed $n+1$, e.g. in codimension 
$1$ a semistable limit cycle has cyclicity $2$. Thus, the cyclicity
of a trivial polycycle (a polycycle without singular points)
in a generic $n$-parameter family does not exceed $n+1$. 

2) (Andronow-Leontovich, 1930s; Hopf, 1940s).
A nontrivial polycycle of codimension $1$ has cyclicity at most $1$.

3) (Takens, Bogdanov, Leontovich, Mourtada, Grozovskii,
early 1970s-1993 (see {\cite{Gr}}, {\cite{KS}}, and references there)).
A nontrivial polycycle of codimension $2$ has cyclicity at most $2$.

\bdef The bifurcation number $B(k)$ is the maximal 
cyclicity of a nontrivial polycycle occurring in a generic 
$k$-parameter family.
\endef

The definition of $B(k)$ does not depend on a choice of the 
base of the family, it depends only on the number $k$ of parameters. 

{\bf{Local Hilbert-Arnold Problem (LHAP).}} {\it{ \cite{I4} 
Prove that for any finite $k$, the  bifurcation number $B(k)$ is 
finite and find an upper estimate for $B(k)$.}}

It turns out that a solution to Local Hilbert-Arnold Problem implies 
a solution to Hilbert-Arnold Problem.

Similarly to the generic smooth vector fields, in the case of analytic 
vector fields one can define so-called {\it a limit periodic set} 
{\cite{FP}}, {\cite{R1}}, which is either a polycycle or 
has an arc of equilibrium points\footnote{ generic vector fields
can not have an arc of equilibrium points}, and formulate 

{\bf Local Finiteness Conjecture (LFC).} {\it \cite{R1}
Prove that any limit periodic set occurring in an analytic family
of vector fields on $\Bbb S^2$ has finite cyclicity in this family}.

Smooth vector fields are more flexible then analytic vector fields
and easier to analyze. A strategy to attack Existential Hilbert
Problem, proposed by Arnold {\cite{AA}} (see also {\cite{IK}}), 
is first understand generic smooth vector fields and then 
try to apply developed methods to analytic vector fields.
Let us summarize the discussion in the form of the diagram:

\begin{figure}[htbp]\label{conj}
  \begin{center}
    \begin{psfrags}
     \includegraphics[width= 5in,angle=0]{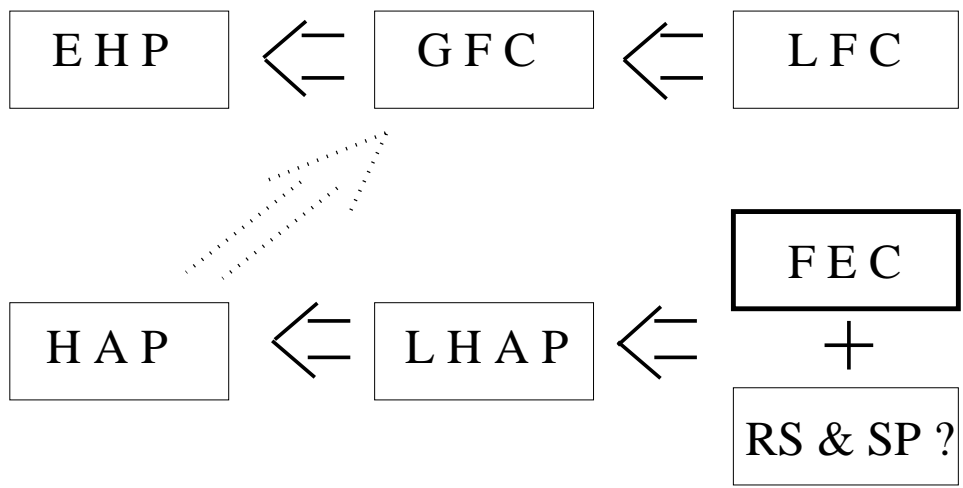}
    \end{psfrags}
    \caption{Existential Hilbert Problem and its offsprings}
  \end{center}
\end{figure}

\subsection{Cyclicity of Elementary Polycycles}

Now we shall formulate the Main Result of this course of lectures.

\bdef
An equilibrium point of a vector field on the two-sphere is called 
{\it{elementary}} if at least one eigenvalue of its linear part is 
nonzero. A polycycle is called an {\it{elementary}} polycycle if all 
its singularities are elementary.
\endef

The Local Hilbert-Arnold problem was solved under the additional 
assumption that a polycycle have elementary singularities only.

\bdef The elementary bifurcation number $E(k)$ is the maximal 
cyclicity of a nontrivial elementary polycycle occurring in 
a generic $k$-parameter family.
\endef

From examples 2) and 3) above it follows that
$$
E(1)=1, \quad E(2)=2.
$$
Information about behavior of the function $k \mapsto E(k)$ 
has been obtained recently. The First crucial step was done
by Ilyashenko and Yakovenko: 
 
{\bf{Finiteness of Elementary Cyclicity (FEC).}}\ {\it{\cite{IY2}
For any $n$ the elementary bifurcation number $E(n)$ is finite.}}

\bcor \label{gha} Under the assumption that families of vector fields 
have elementary singularities only the Global Hilbert-Arnold Problem 
is solved, i.e. any generic finite parameter family of vector fields 
on the sphere $\Bbb S^2$ with a compact base and only elementary
equilibria has a uniform upper bound for the number of limit cycles.
\ecor

{\bf{Main Theorem.}}\ {\it{ {\cite{Ka1}} For any $k\in \Z_+$}} 
\beq \label{eotkbound}
E(k) \leq 2^{25k^2}.
\eneq
This is the first explicit general estimate for cyclicity of polycycle.
The case of a polycycle consisting only one singular point
with no arcs at all is well known. An elementary equilibrium point
can generate limit cycles in its small neighborhood if it is a slow 
focus, that is the linearization matrix has a pair of 
two imaginary eigenvalues. This bifurcation was investigated by 
Takens {\cite{Ta}}. 

\bcor As in Corollary \ref{gha} under the assumption that all the
polycycles are elementary the Main Theorem gives a solution to 
the Local Hilbert-Arnold Problem.
\ecor

\subsection{Resolution of Singularities (RS) or Blow-up
of Singularities of Vector Fields and Singular Perturbations (SP)}

In this subsection we discuss Resolution of Singularities and 
Singular Perturbation  which might lead to generalization of 
the Main Result to a solution to Local Hilbert-Arnold Problem 
(see the box with {\bf RS \& SP?} in the diagram \ref{conj}). 

Let $\dot x=v(x)$ be a $C^\infty$ vector field on $\R^2$ such 
that $v(0)=0$. A vector field satisfies {\it a Lojasiewicz condition} 
if there exists $k\in \Z_+$ and $c>0$ such that
\beq
\|v(x)\|\geq c \|x\|^k 
\eneq
for all $x$ from some neighborhood of $0$.
It can be shown {\cite{D}} that any generic finite-parameter family 
of vector fields on the sphere $\Bbb S^2$ has only vector fields with 
equilibrium points satisfying a Lojasiewicz condition for some 
$k\in \Z_+$ and $c>0$. 

To define {\it a blow up} for a $C^\infty$ vector field $\dot x=v(x)$
on $\R^2$ with an equilibrium at $0$, i.e. $v(0)=0$ consider the map
\beq\label{polar}
\phi:\Bbb S^1 \times \R \to \R^2; \quad 
\phi(\theta,r)\mapsto (r \cos \theta, r \sin \theta).
\eneq
Then the pull-back $\hat v$ with $\phi(\hat v)=v$, is a 
$C^\infty$ vector field on $\Bbb S^1 \times \R$, i.e.
$d\phi_0 (\hat v(0))=X \circ \phi(0)$, where $\hat v$ is the 
blown-up vector field. 

{\bf Desingularization Theorem.} {\it {\cite{D}}
 If a $C^\infty$ vector field $\dot x=v(x)$ on $\R^2$ with
$v(0)=0$ satisfying a Lojasiewicz condition, then there is a finite 
sequence of blow-ups leading to a vector field with only elementary 
equilibria.}

Sometimes this theorem is called Bendixon-Seidenberg-Dumortier 
{\cite{Be}}, {\cite{Se}}, {\cite{D}}. Bendixon stated it without 
a proof it on the brink of the twentieth century. Seidenberg proved 
it in the analytic case and Dumortier did it for $C^\infty$
vector fields with a Lojasiewicz condition. A quantitative version
of the Desingularization Theorem, which estimates number
of necessary blow-ups, was obtained by Kleban {\cite{Kl}}.

This Theorem reduces consideration of {\it an individual} vector
field, occurring in a generic finite-parameter family, with equilibria 
without restriction to an individual vector field with only 
{\it elementary} equilibria.

However, in order to extend an estimate on cyclicity
of elementary polycycles (\ref{eotkbound}) to an estimate on
cyclicity of a generic nonelementary polycycle (LHAP) one needs 
{\it Desingularization Theorem  for families} of generic $C^\infty$
vector fields. Different approaches to attack this problem were
proposed by Denkowska and Roussarie {\cite{DeR}} and by Trifonov
{\cite{Tr}}.

An approach proposed by Trifonov leads the dynamical phenomenon
called {\it Singular Perturbation} (SP): in the simplest case one
needs to analyze families of vector fields on the plane, which for 
some values of parameters have a curve of equilibria. Certainly, 
a generic finite-parameter family of vector fields has no curve of 
equilibria, however, after even one step of blow-up such a curves can 
occur {\cite{Tr}}. Appearance of curves of equilibria after 
a desingularization in a family now seems to be the main obstacle
between an estimate on cyclicity of elementary polycycles 
(\ref{eotkbound}) and Local Hilbert Arnold Problem
(see \cite{Tr}, \cite{IY2}, and \cite{R2} for more). 
 
\section{Bifurcation of Spatial Polycycles and Multiplicity
of Generic Germs}
In this part we present by-product results the Main Theorem. 
The first result is an
extension of the Main Theorem on estimate of cyclicity of planar
elementary polycycle to an estimate on cyclicity of spatial
quasielementary polycycle (see section \ref{space}). The second result
gives an estimate on cyclicity of generic germs of smooth mappings
which is a partial answer to Arnold's question \cite{A2} (see section
\ref{germ}).

\subsection{Bifurcation of Spatial Polycycles}\label{space}

Definition  of a polycycle in a multidimensional case is 
the same as in the planar case verbal (see definition \ref{polyc}). 
When Ilyashenko and Yakovenko proved the Finiteness of Cyclicity 
for elementary polycycles, Arnold posed the question:
{\it{What can be said about bifurcations of spatial polycycles?}} 

Another sufficient reason to look at this problem is, because the 
planar argument (the Poincar\'e-Bendixon Theorem) imply that a
collection of an infinite number of limit cycles of uniformly bounded
length, located in a bounded domain, accumulate to a limit cycle.
Indeed, consider a vector field $\dot x=v(x)$ of a finite codimension 
in $\Bbb R^3$ (dimension $3$ can be replaced by any $N>2$ anywhere in
this section), i.e. a vector field which occurs in a generic 
finite-parameter family. Then $v(x)$ has only isolated singular points. 
Fix a positive number $L$ and assume that in a compact region of 
the phase space there are an infinite number of phase curves of length 
less $L$ corresponding to limit cycles of $v(x)$. Then a subset of 
these limit cycles must accumulate to a separatric polygon (polycycle). 

Bifurcation properties of spatial polycycles are much richer then 
those of planar polycycles. The first important $3$-dimensional
feature is existence of limit cycles that winds several times around 
a polycycle. It happens because a Poincar\'e return map is a
$2$-dimensional map and it might have not only fixed point, but also
periodic points of higher periods too. We call a periodic trajectory
that ``turns'' around a whole polycycle exactly $m$-times before 
closing up an $m$-{\it{cycle}}. Such a trajectory corresponds to
a periodic point of a corresponding Poincar\'e return map of minimal 
period $m$. On the plane because of topological reasons only $1$-cycles 
exist. Definition of cyclicity requires some additional care.

Consider an $n$-parameter family of flows 
$\{\dot x= v(x,\eps)\}_{\eps \in B^n}$ in $\Bbb R^3$. Let 
$\gm \subset \R^3$ be a polycycle of the field $\dot x=v(x,\eps^*)$
for some $\eps^* \in B^n$. Then $\gm$ can be represented as a union
of a finite number of equilibrium points $\{p_j\}_{j\in J}$
and connecting them phase curves $\{\gm_j\}_{j\in J}$. {\it A tube
neighborhood} $T_\gm$ of the polycycle $\gm$ is a union of 
 neighborhoods of equilibria $\{p_j\}_{j\in J}$ and tube 
neighborhoods $\{T_j\}_{j\in J}$ of phase curves 
$\{\gm_j\}_{j\in J}$.

\bdef Let $m \in \Bbb Z_+$. Then $m$-cyclicity of the polycycle $\gm$ 
in the family $\{\dot x= v(x,\eps)\}_{\eps \in B^n}$, denoted by
$\mu(m,\gm)$, is a minimal number $\mu(m,\gm)$ for which there is 
a tube neighborhood $T_\gm$ of the polycycle 
$\gm \subset T_\gm \subset \R^3$, a  neighborhood $V$ of the
parameter $\eps^* \in V \subset B^n$ and for each $j\in J$
Poincar\'e section $L_{\gm,j}$ (a hyperplane) transversal to 
the corresponding $\gm_j$ such that the following condition holds:

1. for any parameter $\eps \in V$ the corresponding vector field 
$\dot x=v(x,\eps)$ has at most  $\mu(m,\gm)$ limit cycles in $T_\gm$;

2. Each of those limit cycles $l_i(\eps)$ intersects each Poincar\'e
section $L_{\gm,j}$ in exactly $m$ different points;

3. In the sense of Hausdorff metric for each $j\in J$ distance between 
each part of $l_j(\eps)$, which lies between two consecutive 
intersections of $L_{\gm,j}$, and the polycycle $\gm$ tends to $0$ as 
$\eps$ tends to $0$. 
\endef

Now we discuss a classical example of a polycycle which has
infinite  $m$-cyclicity  for any $m \geq 1$.

\subsection{The Shilnikov polycycle}

Consider a flow $\phi_t$ in $\Bbb R^3$ with a hyperbolic 
equilibrium point $O$ that has one positive eigenvalue $\lambda$ 
and two complex conjugates $\mu\pm \omega$ with negative real part. 
Suppose that the sum of $\lambda +\mu$  is positive, and the unstable  
one-dimensional manifold returns to the stable one, which is 
two-dimensional. Thus, the equilibrium $O$ has a homoclinic orbit that
tends back to  $O$  along the unstable manifold as $t \to -\infty$, 
and spirals around $O$ on the stable manifold as $t \to +\infty$. 
In 1965 Shilnikov \cite{Sh} discovered that the Poincar\'e map along
this polycycle has a countable number of pairwise disjoint subdomains
so that a restriction to each of them gives a Smale horseshoe. Any of
such horseshoes is structurally stable, therefore, the polycycle
described above (the Shilnikov polycycle) has an infinite
$m$-cyclicity for all $m \in \Bbb Z_+$ (see e.g. {\cite{GH}},
{\cite{IL}}). Codimension of this polycycle is $1$.

\begin{figure}[htbp]
  \begin{center}
    \begin{psfrags}
      \psfrag{S}{$\Sigma$}
      \psfrag{O}{$O$}
      \includegraphics[width= 3in,angle=0]{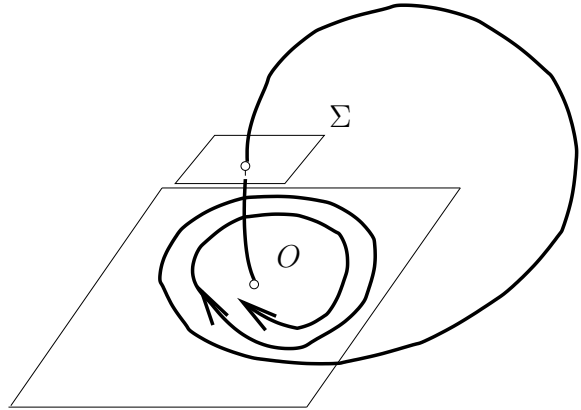}
    \end{psfrags}
    \caption{The Shilnikov Polycycle}
  \end{center}
\end{figure}

However, it seems reasonable to state the following 

{\bf{Conjecture}} {\it{(Arnold-Ilyashenko-Yakovenko)
If a spatial polycycle $\gm \in \Bbb R^3$ has finite codimension $k$
and all its equilibrium points are saddles with real eigenvalues or 
saddlenodes with at most one zero eigenvalue and the other eigenvalues 
are real, then the $m$-cyclicity of $\gm$, denoted $C(m,\gm)$, is 
finite for each $m \in \Bbb Z_+$.}}

Using the ideas and methods for the planar problem and a result 
of Grigoriev-Yakovenko {\cite{GY}} Arnold-Ilyashenko-Yakovenko's 
Conjecture has been solved in arbitrary dimension $N>2$ with 
additional nondegeneracy assumptions on polycycle's equilibria. 
 
\subsection{An estimate of the cyclicity of a quasielementary 
spatial polycycle}

In the planar case we considered polycycles with elementary equilibria 
only, now we define a class of points called {\it{quasielementary 
equilibria}}. The author has shown that polycycles with
quasielementary equilibria only have finite $m$-cyclicity for any 
$m\geq 1$. Moreover, there exists an explicit upper bound for
$m$-cyclicity.

Recall some standard definitions from the normal form theory.

\begin{Def}
The set of complex numbers $\lambda_1, \dots, \lambda_N \in \Bbb C$ 
is called :

a) nonresonant if there is no integral relation among the numbers
$\lambda_j$ of the form $\lambda_j=\sum_{i=1}^N k_i \lambda_i$,
where $k_i \in \Bbb Z_+$ for $i=1,\dots ,n$ and 
$\sum_{i=1}^N k_i \geq 2$.

b) {\it strongly simply resonant} if all the nontrivial resonance
relations $\lambda_j=\sum_{i=1}^N k_i \lambda_i$ follows from the
single one $\sum_{i=1}^N k^*_i \lambda_i=0$, where $k_i \in \Bbb Z_+,\ 
i=1, \dots ,n$ and $\sum_{i=1}^N k_i \geq 2$.
\end{Def}

\begin{Def}
We shall call an equilibrium point of a differential equation
quasi\-e\-le\-men\-tary, if the linearization matrix of the equation 
at this point has only {\it real} eigenvalues, at most one 
of them is zero, and they satisfy one of the following conditions:

1) they are nonresonant and we call such an equilibrium 
{\it a nonresonant saddle};

2) they form a strongly simply resonant set of numbers---
{\it a strongly simply resonant saddle};

3) one eigenvalue is zero with Lojasiewicz exponent 2 and the others 
form a nonresonant set.

A polycycle is called {\it quasielementary} if all its 
vertices are quasielementary. 
\end{Def}

Note that the class of quasielementary points in the case of 
the plane ($N=2$) coincides with the class of elementary points,
except of multiplicity two condition for saddlenodes.
In a sense, Theorem \ref{Space} below is a generalization of 
Theorem \ref{eotkbound}.

\begin{Def}  
The {\it quasielementary bifurcation number} $QE(N,n,m)$ is 
the maxi\-mal $m$-cyclicity of a quasielementary polycycle occurring 
in a generic $n$-para\-meter families of vector fields in $\R^N$.
\end{Def} 
 
\begin{Thm} \label{Space} \cite{Ka2}
For any positive integer $N$ (dimension of the phase space), $n$
(number of parameters), $m$ (number of turns around a polycycle), 
and $T=6Nnm$ we have 
$$
QE(N,n,m) \leq 2^{T^2}.
$$ 
\end{Thm}

In the next section we describe another by-product of the
Main Theorem.

\subsection{Geometric multiplicity of germs of generic maps}\label{germ}

Let $F:\Bbb R^n \to \Bbb R^n$ be a generic $C^{k}$ smooth map, 
$k \geq n+1$. Fix a point $a \in \Bbb R^n$ and denote $F(a)$ by $b$.

\bdef A geometric multiplicity of a map germ 
$F:(\Bbb R^n,a) \to (\Bbb R^n,b)$ at $a$, denoted by
$\mu^G_a=\mu^G_a(F)$, is the maximal number of isolated preimages
$F^{-1}\left(\tilde b\right)$ close to $a$:
\beq \label{geom}
\mu^G_a(F)=\limsup_{r \to  0} \sup_{\tilde b \in \Bbb R^n}
\#\{x \in B_r(a): F(x)=\tilde b\}.
\eneq
\endef
For example, the geometric multiplicity of the function $f:x \to x^2$ 
at $0$ is two, but the geometric  multiplicity of $f:x \to x^3$ at $0$
is one, even though $0$ is a degenerate point of the second order.

In the complex case the geometric multiplicity equals the usual
multiplicity (see e.g. \cite{AGV}). In the real case the first is 
no greater than the second.
 
\bdef Define geometric multiplicity of $n$-dimensional germs,
$\mu^G(n)$, as follows. Let $F:\Bbb R^n \to \Bbb R^n$ be a generic
map. The geometric multiplicity of $F$ equals the least upper bound of 
geometric multiplicities of $\mu^G_a(F)$ taken over all points 
$a \in \Bbb R^n$. Then the geometric multiplicity of $n$-dimensional
germs is the maximum of the geometric multiplicities of all generic 
maps $F$ from $\Bbb R^n$ to $\Bbb R^n$
\beq
\mu^G(n)=\sup_{F-\text{generic},\ a \in \Bbb R^n} \mu^G_a(F).
\eneq
It turns out that the geometric multiplicity of $n$-dimensional germs 
is finite for all positive integer $n$ and depends only
on dimension $n$.
\endef
\brm
For example, for $n=2$ the Whitney Theorem about maps of surfaces
states that a generic map of two dimensional manifolds $F:M^2 \to N^2$ 
can have only three different types of germs: 1-to-1, a fold, and a
pleat (see e.g. \cite{AGV}). This implies that $\mu^G(2)=3$.
\erm

A natural problem posed by Arnold \cite{A2} is
{\it{to give estimates for the geometric multiplicity  $\mu^G(n)$ of
$n$-dimensional germs.}}

In the case of complex analytic maps of $\C^n$ into $\C^n$
Gabrielov and Khovanskii {\cite{GK}} Thm.7 obtained an estimate
on $\mu^G(n)$ of the type $\mu^G(n)\leq n^n$. 
The upper bound for the geometric multiplicity for $n$-dimensional 
{\it smooth} germs of generic maps is given by

\bthm \label{working} {\cite{Ka1}} 
The geometric multiplicity of germs of a generic $C^k$ smooth
map  $F:\Bbb R^n \to \Bbb R^n,\ k >n$ admits the 
following upper estimate:
\beq \label{2nto2}
\mu^G_a(F) \leq 2^{n(n-1)/2+1}n^n,\  \forall\  a \in \Bbb R^n.
\eneq
\ethm

Using the same method one can prove 

\bthm \label{working1} {\cite{Ka1}} 
Let $F:\Bbb R^n \to \Bbb R^N$ be a generic $C^k$ smooth map
with $k >n,\ N \geq n$ and $P:\Bbb R^N \to \Bbb R^n$ be
a polynomial of degree $d$. Then the geometric multiplicity of 
germs of  a chain map $P \circ F:\Bbb R^n \to \Bbb R^n$
admits the following upper estimate:
\beq
\mu^G_a(P \circ F) \leq 2^{n(n-1)/2+1} (dn)^n,\  
\forall\  a \in \Bbb R^n.
\eneq
\ethm
An interesting feature of this theorem is that
the geometric multiplicity does not depend on dimension $N$
of the intermediate space.

The problem about an upper estimate of geometric multiplicity of germs
of generic smooth maps $F:(\R^n,0) \to (\R^n,0)$ or chain maps 
$P \circ F:(\R^n,0) \to (\R^n,0)$ is closely related to the problem
about an estimate cyclicity of elementary polycycles as the reader
will see below.

All the results (The Main Theorem, Theorem \ref{Space},
Theorem \ref{working}, and Theorem \ref{working1})
were announced in \cite{Ka3}.

\section{Three stages of the proof of the Main Theorem
and outline of the content of the lectures}
\label{program}

The Main Theorem is an quantitative extension of Ilyashenko-Yakovenko 
Finiteness Theorem. The paper of Ilyashenko-Yakovenko {\cite{IY2}} 
was a corner stone for the present paper. In {\cite{IY2}} 
the authors made an important step:
they found a pass from {\it{bifurcation theory \textup{to} 
singularity theory}} using the Khovanskii reduction method 
{\cite{Kh1}}. In \cite{Ka1} we follow this pass up at the beginning 
and using some new ideas get an estimate
for the cyclicity of elementary polycycles. Below we outline
the main steps of the proof of the Main Theorem and describe
the content of the coming lectures.

The proof of the Main Theorem consists of three steps.
Relation to the proof of the Finiteness Theorem {\cite{IY2}}
is discussed in section \ref{relation}

{\it{Step 1. Normal forms for local families of vector fields
and their integration.}}\ \ We use normal forms to establish an
explicit form for the Poincar\'e correspondence map near 
equilibrium points on the polycycle under consideration. In \cite{MR}
and later in \cite{IY2} it is noticed that these maps satisfy Pfaffian
(polynomial differential) equations with coefficient of polynomials 
depending smoothly on the parameters of the family. As the result a 
{\it{basic system}} of equations, determining the number of limit 
cycles, is obtained.

{\it{Step 2. The Khovanskii reduction method.}}\ \ We discuss a
variation of the Khovanskii method \cite{Kh2}. This method allows us to 
investigate systems of equations that involve functions satisfying 
Pfaffian equations. It turns out that the number of solutions to the
basic system can be estimated by the number of solutions to a
{\it{mixed functional-Pfaffian}} system. After an application of 
the Khovanskii method to the mixed functional-Pfaffian system we
obtain several {\it{chain maps}}: the maps of the form
\beq \label{chain}
x \mapsto (P_1,\dots ,P_n)\circ 
\left(x,f(x),f'(x),\dots,f^{(n)}(x)\right),
\eneq
where $x$ is a point nearby $0 \in \R^n$,  $f$ is a generic function,
$f^{(k)}(x)$ is collection of all derivatives of $f$ of order $k$, and 
$P=(P_1,\dots ,P_n)$ is a vector-polynomial given by its coordinate 
functions of known degree. 

It turns out that the problem of estimating the number of limit cycles 
reduces to estimating the number of small regular preimages of some 
{\it special} points by the chain map. Special points form an 
open cone-like semialgebraic set $K$ approaching to $0$ in the image, 
e.g. if $K \subset \R^2,$ then $K=\{(x_1,x_2): 0<x_2<x_1^m\}$ for
some $m\in \Z_+$.

Denote by $F$ the map $F:x\mapsto (x,f(x),f'(x),\dots,f^{(n)}(x))$
which is called the $n$-th jet of $f$. Denote by $L_F$ the 
linearization of $F$ at point $x=0$.

Lecture 2 highlights Steps 1 and 2 in a simple nontrivial case. 

{\it{Step 3. Bezout's theorem for the Chain maps.}}\ \  
We shall construct an algebraic set $\Sigma$ in the image of $F$ 
(in the space of $n$-jets) so that if $F$ is transversal to $\Sigma$,
then the number of preimages of any point $a$ from a set of special
points $K$ is {\it{the same}} for $F$ and its linearization $L_F$ at
zero: 
\beq \label{bez}
\#\{x: P\circ F(x)=a\}=\#\{x: P\circ L_F(x)=a\}
\leq \prod_{j=1}^k \deg P_j.
\eneq
Since $L_F$ is a linear map, one can apply Bezout's theorem to estimate 
the right-hand part of the equality. This observation completes the proof 
of the Main Theorem. 

In order to prove existence of such a set $\Sigma$ we need to apply 
stratification theory originated by works of Whitney \cite{Wh}, 
Thom \cite{Th1}, and Mather \cite{Ma}. More exactly, we need to prove
existence of so-called $a_P$-stratification introduced by Thom
in some special case \cite{Ka1}.  Lecture 3 presents necessary 
notions from stratification theory and states the required result
on existence of $a_P$-stratification.

\subsection{Multichain maps and bifurcation of spatial polycycles}

In order to get an estimate on cyclicity of spatial polycycles Theorem 
\ref{Space} we face the problem of estimating geometric multiplicity of 
{\it multichain maps} of the form
\beq \label{multichain}
P\circ (F,F):B^n \times B^n \to \R^{2n},
\eneq 
where $B^n\subset \R^n$ is a unit ball, $F:B^n \to \R^N$ is a generic 
map and $P:\R^{2N}\to \R^{2n}$ is a vector-polynomial of known degree. 
Appearance of this problem is described with many more details 
in lecture 4. It is no longer possible to treat the $2$-tuple map 
\beq
(F,F):B^n \times B^n \to \R^N \times \R^N
\eneq
as a generic map. 

{\it{ Step 4. Blow-up along the diagonal in the multijet space}}

Grigoriev and Yakovenko \cite{GY} constructed a so-called
{\it{space of divided differences or ${\Cal{DD}}_2$-space}} and the following
commutative diagram:

\begin{figure}[htbp]
  \begin{center}
    \begin{psfrags}
      \psfrag{Phase}{ $B^n \times B^n$}
      \psfrag{MultJ}{$\R^{2N}$}
      \psfrag{DivD}{${\Cal{DD}}_2(B^n,\R^N)$}
      \psfrag{pi}{$\pi_2$}
      \psfrag{N}{${\Cal D}_2(F)$}
      \psfrag{(jf,...,jf)}{$(F,F)$}
      \includegraphics[width= 3.5in,angle=0]{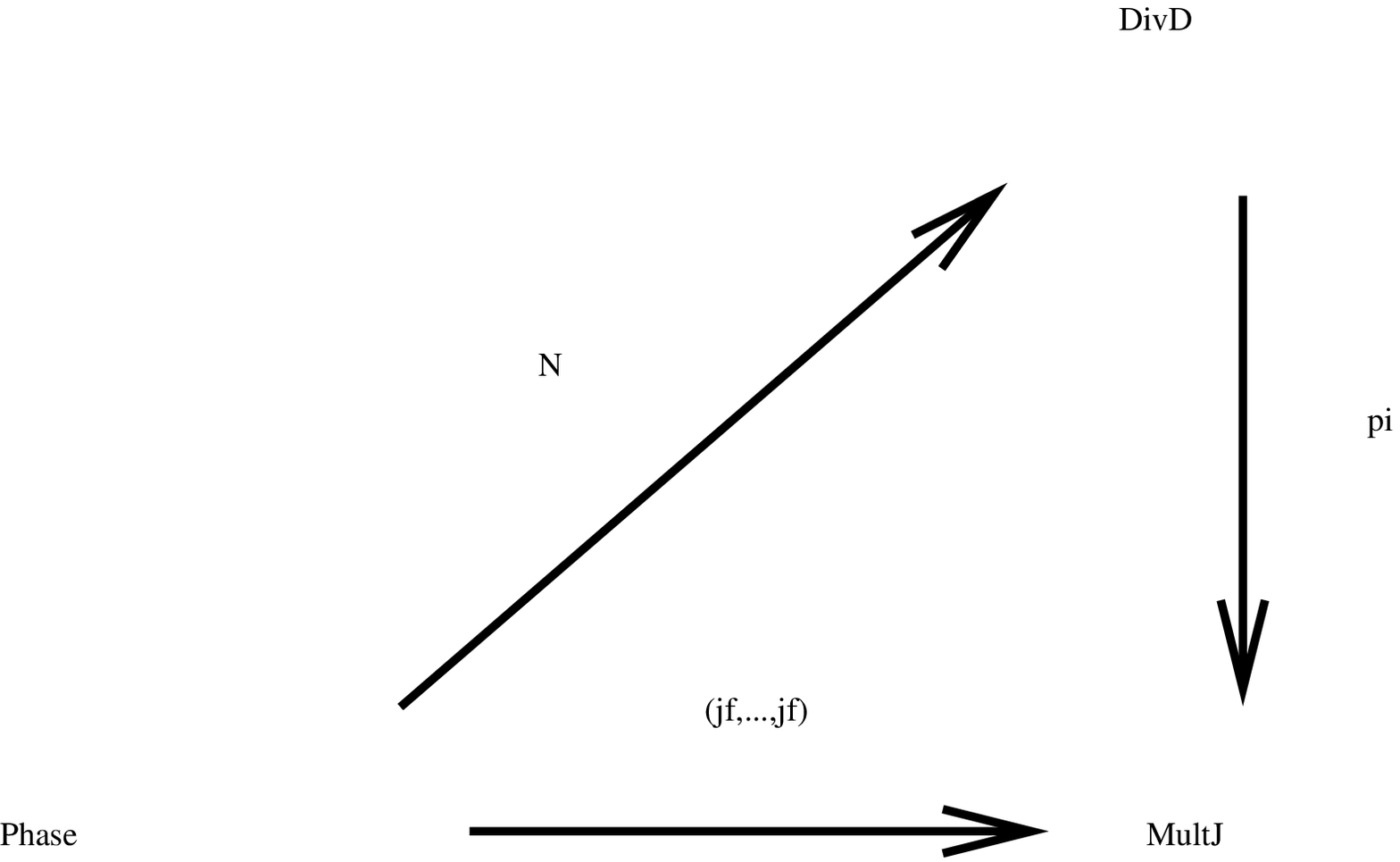}
    \end{psfrags}
    \caption{Polynomial blow-up of the multijet space}
  \end{center}
\end{figure}
where ${\Cal D}_2(F): B^n \times B^n \to {\Cal{DD}}_2(B^n,\R^N)$ is 
a smooth map, $\pi_2:{\Cal{DD}}_2(B^n,\R^N) \to \R^{2N}$ 
is an explicitly computable polynomial and
\beq
\pi_2 \circ {\Cal D}_2(F)=(F,F):B^n \times B^n \to \R^{2N}.
\eneq

It turns out that one can treat ${\Cal D}_2(F)$ as a generic map
for a generic $F$ and impose various transversality conditions. 
Therefore, we can represent the multichain map $P \circ (F,F)$, 
given by (\ref{multichain}), in the form
\beq \label{recomp}
P \circ (F,F)=(P \circ \pi_2) \circ {\Cal D}_2(F),
\eneq
where $P \circ \pi_2$ is a polynomial, since $\pi_2$ is a polynomial,
and ${\Cal D}_2(F)$ is a smooth map. Moreover, it turns out that 
${\Cal D}_2(F)$ is generic for a generic $F$. Now we can apply 
Bezout's Theorem to the chain map $(P \circ \pi_2) \circ {\Cal D}_2(F)$. 
In lecture 4 we describe this construction with details and in greater 
generality and exhibit application of this construction to an old problem 
of the rate of growth of the number of periodic points for generic
diffeomorphisms in smooth dynamic systems ( see e.g. \cite{AM} and 
\cite{Sm}).

%
%
%

\subsection{Relation of the proof of Main Theorem and 
Ilyashenko-Yakovenko Finiteness Theorem \cite{IY2}}\label{relation}
 
Steps 1 of both proofs \cite{IY2} and \cite{Ka1} are the same. 
We shall just present the table of
required normal forms from {\cite{IY2}}, which were obtained in
\cite{IY1}. Step 2 in this proof is slightly different from the one
in {\cite{IY2}} and this is the first novel point. After application
of the Khovanskii method  we obtain the same collection of chain maps 
of the form (\ref{chain}) as in \cite{IY2}. However, in \cite{IY2}
the authors investigate the number of regular preimages of points in 
the image by the chain maps {\it{without any restriction}} on those
points. In the present proof, using new additional arguments in 
the Khovanskii method, we reduce consideration to only preimages of 
{\it{special}} points, i.e. points from a tiny cone-like set in the
image. At this point our proof goes independently. 
The proof from \cite{Ka1} can be considered as an independent
simplified proof of Ilyashenko-Yakovenko's Finiteness Theorem 
by modulo of derivation of mixed functional-Pfaffian system. 

{\it Acknowledgments:}\ \ I would like to thank Dana Schlomiuk for 
giving an great opportunity to give a course of lectures in a workshop
held in Montreal during June 2000. Special thank goes to my teacher 
Yulij Ilyashenko who patiently taught me bifurcation theory and 
whose long lasting support and encouragement have been crucial source
of energy. Discussions with William Cowieson, Andrei Gabrielov, 
Askold Khovanskii, Pavao Mardesic, John Mather, Robert Moussu, 
Oleg Shelkovnikov, Sergei Yakovenko have been very helpful for me. 
I would like also to thank Christiane Rousseau and Robert Roussarie 
for inviting me to give series of lectures on the bifurcation of limit 
cycles in Montreal and Dijon and for useful discussions. These series 
of lectures have been extremely helpful to improve the presentation of 
the present here lectures. I would like to acknowledge financial support 
and very productive atmosphere of Institute for Physical Sciences and 
Technology University of Maryland and Courant Institute of Mathematical 
Science, NYU, where various parts of the work have been done. While in 
the Institute for Physical Sciences and Technology I had enjoyed 
fruitful discussions with James Yorke and Brian Hunt.

 \chapter{Normal Forms and The Khovanskii Method.}

\medskip

We explain the proof of the Main Theorem in the simplest
nontrivial case $n=2$. Consider a generic $2$-parameter family of
vector fields $\{\dot x=v(x,\eps)\}_{\eps \in B^2}$ and suppose that 
for $\eps=0$ the vector field $\dot x=v(x,0)$ has a polycycle $\gm$ 
which consists of two saddles $p_1,\ p_2$ and two separatrices 
connecting $\gm_1,\ \gm_2$ them. Consider a segment $\Sigma$
transversal to, say, $\gm_1$ and denote by 
$\Delta: \Sigma \supset U \to \Sigma$ the Poincare return map,
which is define on some open set $U$ in $\Sigma$. In order to 
estimate the number of limit cycles bifurcating from the polycycle
$\gm$ we need to {\it estimate the number of isolated fixed points  
$\# \{x \in U :\ \Delta(x)=x\}$.} 

Using the standard approach we decompose the Poincare map $\Delta$ 
into a composition of four maps: two local $\Delta_1$ and $\Delta_2$ 
in neighborhoods of equilibria $p_1$ and $p_2$ respectively and 
two semilocal $f_1$ and $f_2$ along connecting separatrices $\gm_1$ 
and $\gm_2$ respectively to be defined precisely below (see Fig.2.1). 
After that we replace the equation $\Delta(x)=x$ by the system
of equations corresponding to $\Delta \equiv \Delta_2 \circ f_2 \circ
\Delta_1 \circ f_1$. To understand properties of local maps
$\Delta_i,\ i=1,2$ we use normal forms theory.

\begin{figure}[htbp]
  \begin{center}
    \begin{psfrags}
      \psfrag{S2+}{$\Si_2^+$}
      \psfrag{S1+}{$\Si_1^+$}
      \psfrag{S1-}{$\Si_1^-$}
      \psfrag{S2-}{$\Si_2^-$}
      \psfrag{O2}{$p_2$}
      \psfrag{O1}{$p_1$}
      \psfrag{D1}{$\Delta_1$}
      \psfrag{D2}{$\Delta_2$}         
      \psfrag{f1}{$f_1$}
      \psfrag{f2}{$f_2$}         
     \includegraphics[width= 3.3in,angle=0]{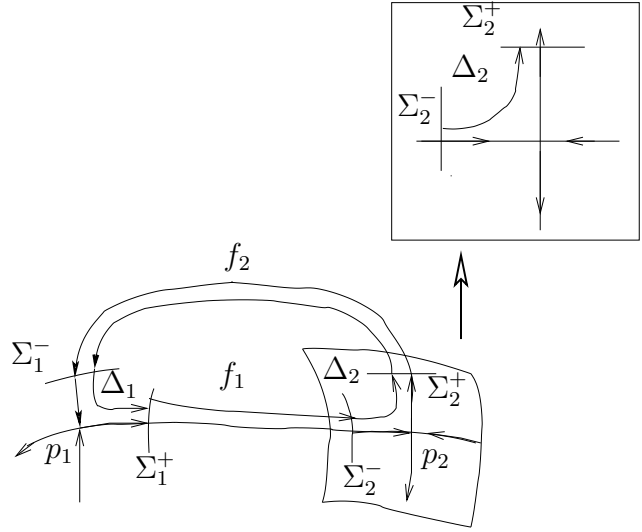}
    \end{psfrags}
    \caption{Construction of ``entrance'' and ``exit'' transversals}
  \end{center}
\end{figure}

\section{Normal forms and a Basic system determining the number
of limit cycles}
\subsection{Polynomial Normal forms of local families and 
Pfaffian Poincare return maps}

It turns out that in a small neighborhood of an elementary equilibrium 
point there exists a finitely differentiable normal coordinates
(in the Cartesian product of the phase space and the parameter space),
so-called normal forms of an equilibrium point. The list of
finitely differentiable normal forms was obtained in {\cite{IY1}}.
The main feature of the list: all normal forms are {\it{polynomial
and integrable.}} The smaller is the neighborhood of a normal form,
the higher is its smoothness. So smoothness can be chosen arbitrary
large. All normal forms are summarized in Table 1 below.

In a small neighborhood of an elementary equilibrium point one can
choose two small segments, say $\Si^-$ and $\Si^+$, transversal to 
the vector field for the critical value of parameter and explicitly 
calculate the Poincare (correspondence) map which maps a point from 
one segment say $\Si^-$ along the corresponding phase curve to a point 
from the other segment $\Si^+$ (see Fig.2.1). For an appropriate choice
of segments $\Si^-,\Si^+$ and coordinate functions $x,y$ in
$\Si^-,\Si^+$  respectively, and a smooth function $\lb(\eps)$ in 
the original parameter $\eps$ of the family the Poincare return map 
$\Delta_\eps:x\to y$ can be explicitly computed.
Moreover, there is a Pfaffian 1-form $\om$ of the form, i.e. 1-form
of the form
\beq \label{pfaff}
P(x,y,\lb(\eps))\ dx+\ Q(x,y,\lb(\eps))\ dy=0
\eneq
vanishing on the graph $y=\Delta_\eps(x)$, where 
$P(x,y,\lb(\eps))$ and $Q(x,y,\lb(\eps))$ are polynomials. 
This was first noticed by Moussu-Roche \cite{MR}.

\bex Consider a nonresonant saddle on the plane. There is a normal
form which after an appropriate rescaling is given by the equation

\beq \label{example}
\begin{cases}
\dot x=\lb_1(\eps) x\\
\dot y=-\lb_2(\eps) y,
\end{cases}
\eneq
where $\lb_1(\eps)$ and $\lb_2(\eps)$ are smooth functions
and two transversal ``exits''-''entrance'' sections are
$\Si^-=\{y=1\}, \quad  \Si^+=\{x=1\}.$
\enex

\begin{figure}[htbp]
  \begin{center}
   \begin{psfrags}
     \psfrag{x}{\small{$x$}}
     \psfrag{y=D(x,l)}{\small{$(y=\Delta(x,\lambda),1)$}}
     \psfrag{Dcmu}{\small{$D^c_\mu$}}
     \psfrag{Smu}{\small{$S_\mu$}}
     \psfrag{Dhmu}{\small{$D^h_\mu$}}
     \psfrag{tilq}{\small{$\tilde q$}}
     \psfrag{q}{\small{$q$}}
     \psfrag{p}{\small{$p$}}
    \includegraphics[width= 4in,angle=0]{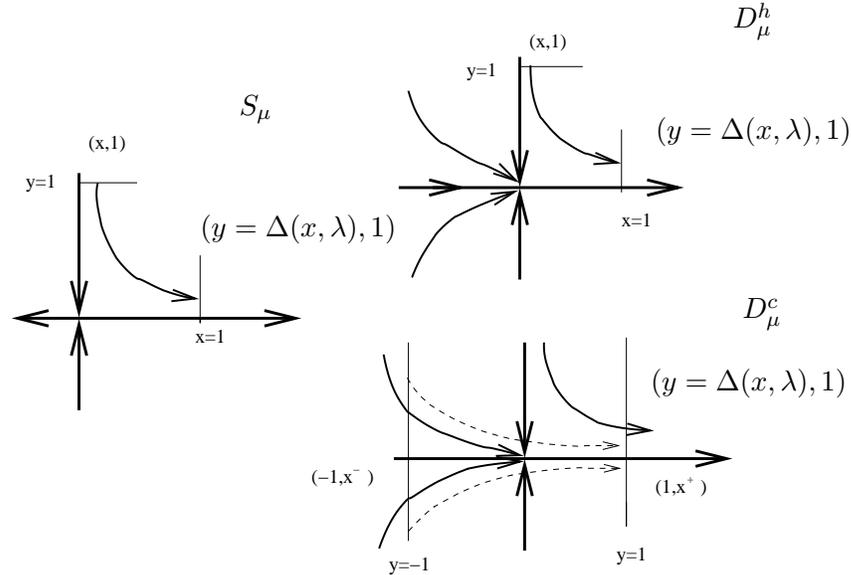}
    \end{psfrags}
   \caption{Poincare Correspondence maps}
  \end{center}
\end{figure}

Then for $\lb(\eps)=\lb_1(\eps)/\lb_2(\eps)$ the function 
$u(t)=x(t)y^{\lb(\eps)}(t)$ is the first integral. Therefore, if 
a trajectory starts form $(x(0),y(0)) \in \Sigma^-$ and ends 
at $(x(t^*),y(t^*))\in \Sigma^+$ (see Fig.2.2 $S_\mu$-case), then 
$y^{\lb(\eps)}(t^*)=x(0)$ or in the induced on $\Sigma^-$ and 
$\Sigma^+$ coordinates from $\R^2$ we have  $y^{\lb(\eps)}=x$. 
It is easy to see that the $1$-form $\om=x\ dy+\ \lb(\eps)y\ dx$ 
vanishes on the graph $y^{\lb(\eps)}=x$.
All necessary information about normal forms, Poincare correspondence
maps, and corresponding Pfaffian forms is in the Table 1 below.
For completeness let's give necessary definitions of from 
theory of normal forms.

\subsection{Definitions and a collection of Normal forms}

 A local family of planar vector fields is the germ of a map,
$$
v:(\R^2,0)\times(\R^k,0)\to(\R^2,0),\qquad
(x,y,\eps)\mapsto v(x,y,\eps).
$$
A $C^r$-smooth conjugacy between two local families $v$ and $w$ of the
above form is a map
\beq \nonumber
H:(\R^2,0)\times(\R^k,0)\to(\R^2,0),\qquad
(x,y,\eps)\mapsto H(x,y,\eps),
\eneq
such that
\beq \nonumber
\ \ \ 
H_* v(x,y,\eps)=w(H(x,y,\eps),\eps),
\eneq
where $H_*$ stands for the Jacobian matrix with respect to the
variables $x,y$.  (this definition does not yet allow for
reparametrization of a local family).  Two families are finitely
differentiably equivalent, if for any $r<\infty$ there exists a
$C^r$-conjugacy between them. The two families $v,w$ are orbitally
equivalent, if there exists the  germ of a  nonvanishing function
$\phi:(\Bbb R^2,0)\times(\Bbb R^k,0)\to\Bbb R^1$ such that $v$ is
equivalent to $\phi\cdot w$.

To allow for a reparametrization of local families, we say that a
family $v(\cdot, \eps)$ is induced from another family $w(\cdot,\lb),\
\lb\in(\Bbb R^m,0)$, if $v(\cdot,\eps)=w(\cdot,\lb(\eps))$, where 
$\lb(\eps)$ is the germ of a smooth map $(\Bbb R^k,0)\to(\Bbb R^m,0)$. 
The number of new parameters $m$ may be different from $k$.

Assume that the family $w(\cdot,\lb)$ is {\it global\/} (i.e. the
expression $w(x,y,\lb)$ makes sense for all $(x,y,\lb)\in\Bbb R^{m+2}$);
this happens in particular when $w$ is {\it polynomial\/} in all its
arguments. Restricting the parameters $\lb$ onto a small neighborhood of
a certain point $(0,0,{\bf c})\in \Bbb R^2\times\Bbb R^m$, we obtain a
{\it localization\/} of the global family $w$, which 
formally becomes a local family after the parallel translation 
$\lb\mapsto\lb-\bf c$.

\bdef
1. A local family $v=v(\cdot,\lb)$ is finitely smooth orbital versal
unfolding (in short, versal unfolding) of the germ $v(\cdot,0)$, if
any other local family unfolding this germ is finitely differentiable 
orbitally equivalent to a family induced from $v$.

2. A polynomial family $w(\cdot,\lb)$, $\lb\in\Bbb R^m$, is a {\it 
global finitely smooth orbital versal unfolding\/} (in short, global
versal unfolding) for a certain class of local families of vector
fields, if any local family from this class is finitely differentiable
orbitally equivalent to a local family induced from some localization
of $w$. 
\endef

To investigate a versal unfolding means to investigate at the same time
all smooth local finite-parametric families which unfold the same germ
$v(\cdot,0)$.  The main result describing versal unfoldings of germs of
elementary singularities on the plane, is given by the following

\bthm \label{nforms}{\cite{IY1}}
Suppose that a generic finite-parameter family of smooth vector fields
on the plane possesses an elementary singular point for a certain value
of the parameters. If this point has at least one hyperbolic sector,
than the family is finitely differentiable orbitally equivalent to a
family induced from some localization of one of the families given in
the second column of Table~1.
\ethm

\newpage
\begin{center}
Table 1. Unfolding of elementary equilibrium points on the plane.
\end{center}
$$
\ 
$$
{\small
\begin{tabular}{|c |c |c |c |}
\hline & & & \\
 Type &  Normal forms
 & Poincare & Pfaffian equations \\
 & & Correspondence maps & 
 \\ \hline & & & \\
 & $\dot x\ =\ x,$ & & \\
 $S_0$ & $\dot y\ =\ -\lb y.$ & $y=x^\lb,$ & 
 $x\ dy\ - \lb y\ dx\ =0$ \\ &  & & \\
 & $\lb=\lb_0 \in \R^1$& $x>0,\ y>0$ & \\ & & & \\
 \hline & & & \\
 & $\dot x\ =\ x\ \left( \frac{n}{m}+ P_\mu(u,\lb)\right),$ & & \\
 & $\dot y\ =\ - y.$ &  & \\
 $S_\mu$& & $0\ =\ m \log y\ +$ & $y\ P_\mu(y^n,\lb)\ dx\ -$ \\
 & $u=u(x,y)=x^m\ y^n,$ & & \\ & 
 & $\int_{x^m}^{y^n} \frac{du}{uP_\mu(u,\lb)}.$ 
 & $\left( \frac{n}{m}+ P_\mu(y^n,\lb)\right)\times$\\
 & $P_\mu(u,\lb)=\pm u^\mu(1+\lb_\mu u^\mu)$ & & \\ 
 & $+ W_{\mu-1}(u,\lb),$ & $x>0,\ y>0$ & 
 $x P_\mu(x^m,\lb)\ dy=0$\\ & $\lb=(\lb_1,\dots ,\lb_\mu)$ 
 & & \\ & & & \\ \hline & & & \\
 & $\dot x\ =\ Q_\mu(x,\lb),$ & & \\
 & $\dot y\ =\ - y.$ & $y\ =\ C(\lb) x,$ & \\
 $D^{c}_\mu$ &  & $C={\int_{-1}^1} \frac{du}{Q_\mu(u,\lb)},$ 
 & $x\ dy\ - y\ dx=0$ \\ & &  & \\ 
 & $Q_\mu(x,\lb)\ =\pm x^{\mu+1}(1+\lb_\mu x^\mu)$
 & $x,\ y\ \in \R^1$ & \\ \cline{1-1} \cline{3-4}
 & $+ W_{\mu-1}(x,\lb),$ &  & \\
 & & $0=\ \log y\ + $ & \\
 $D^{h}_\mu$& $\lb=(\lb_1,\dots ,\lb_\mu)$ & 
 ${\int_{x}^1}\frac{du}{Q_\mu(u,\lb)}$ &
 $Q_\mu(x,\lb)\ dy\ -$ \\
 & & $y>0,\ x\ \in \R^1$ & $y\ dx\ =0$ \\ & & & \\ \hline
\end{tabular}}
$$
\ 
$$

In the first column we use the following notation for elementary 
equilibria (the subscript indicates the degree of degeneracy):

$S_0$--- Nonresonant saddle;

$S_\mu$--- Resonant saddle whose quotient equation (the differential 
equation for $u=x^m\ y^n$ below) has the singular point of 
multiplicity $\mu+1$ at the origin, $\mu \geq 1$; if we want to specify 
explicitly the resonance between the eigenvalues, we use the extended 
notation $S_\mu^{(n:m)}$ assuming that the natural numbers $m,n$ are 
mutually prime;

$D_\mu$--- Degenerate saddlenode of multiplicity $\mu$;

$W_{\mu-1}(z,\lb)=\lb_0+\lb_1z+\cdots+\lb_{\mu-1}z^{\mu-1}$ is a
Weierstrass polynomial of degree $\mu-1$.

Different technical remarks concerning this table see in {\cite{IY2}}
$\S 1.1$. We just briefly describe each column.

The second  column has the corresponding normal forms. 
In the third column of the table the Poincare correspondence maps
$y=\Delta(x,\lb)$ for the polynomial normal forms are given. They 
are implicitly defined by the equations relating $x$ to $y$, these
equations depending explicitly on the parameters $\lb$ and thus
implicitly on the original parameters $\eps$. The
choice of segments transversal to the phase curves of the family
described in Fig. 2.1. The last column has  Pfaffian equations
vanishing on the graphs of corresponding Poincare maps.

\subsection{Singular-regular systems determining the number of 
li\-mit cycles} 

Recall that for simplicity we consider a $2$-parameter family of 
vector fields $\{\dot x=v(x,\eps)\}_{\eps \in B^2}$ and suppose 
that for $\eps=0$ the vector field $\dot x=v(x,0)$ has a polycycle 
$\gm$ which consists of two saddles $p_1,\ p_2$ and two 
separatrices connecting $\gm_1,\ \gm_2$ them (see Fig.2.2). For each 
saddle $\{p_j\}_{j=1,2}$ there is a neighborhood $\{U_j\}_{j=1,2}$ 
with a {\it{$C^r$-normal coordinate charts}}. 
Consider transversal segments ``entrance'' $\Si^-_j$ and ``exit''
$\Si^+_j$ which are parallel to coordinate axis of the normal chart
such that the phase curve $\gm_{j-1}$ enters the neighborhood $U_j$ 
through $\Si^-_j$ and the phase curve $\gm_{j}$ exists $U_j$ 
through $\Si^+_j$. The normal coordinates induce coordinates $x_j$ and 
$y_j$ on $\Si^-_j$ and $\Si^+_j$ respectively. For some parameter values
the corresponding vector field defines the following collection
of Poincare correspondence maps:
\beq
\begin{aligned}
\Delta_j(\cdot,\eps)& :x_j\to y_j=\Delta_j(x_j,\eps),\ j=1,2\\
f_j(\cdot,\eps)& :y_j \to x_{j+1}=f_j(y_j,\eps), \ j=1,2\ \ (\mod 2), 
\end{aligned}
\eneq
where $\Delta_j(\cdot,\eps)$ is a local Poincare map form the
``entrance'' segment $\Si^-_j$ to the ``exit'' segment $\Si^+_j$ and 
$f_j(\cdot,\eps)$ is a semilocal Poincare map along the phase curve
$\gm_j$ form the ``exit'' segment $\Si^+_j$ to the ``entrance''
segment $\Si^-_{j+1}$.

Now we decompose the monodromy map (the Poincare first return map)
along the polycycle $\gm$ into the chain of two the local singular maps  
$\{\Delta_j(\cdot,\eps)\}_{j=1,2}$ and the semilocal regular maps
$\{f_j(\cdot,\eps)\}_{j=1,2}$ of the total length $4$. Limit cycles 
correspond to the fixed points of the monodromy. But instead of
writing one equation for the fixed points of the monodromy we consider 
a system of $4$ equations, which we call {\it the preliminary basic
system}:
\beq
\begin{cases} \label{prebasic}
y_1 = \Delta_1 (x_1, \eps), \\
x_2-f_1(y_1,\eps)=0,\\
y_2 = \Delta_2 (x_2, \eps),\\
x_1-f_2(y_2,\eps)=0. 
\end{cases}
\eneq

Recall that $x_j$'s are $C^r$-normal coordinates on $\Si_j^-$ and 
$y_j$'s are  $C^r$-normal coordinates on $\Si_j^+$. Thus the system 
involves $C^r$-smooth regular functions $f_j$'s and the maps  
$\Delta_j$ from the list (modulo reparametrization 
$\eps \to \lb (\eps)$), that are essentially singular. The problem now 
is to estimate the number of small isolated solutions uniformly over all 
sufficiently small parameter values.

Suppose for $\eps=\eps^*$ the system (\ref{prebasic}) has the maximal
number of isolated solutions. Since each isolated solution of this
system corresponds to an isolated solution of the 1-dimensional
Poincare return map $\Delta(x_1,\eps)=x_1$, one can choose a
small $\dt_1$ so that the number  of regular
(nondegenerate) solutions of $\Delta(x_1,\eps)=x_1+\dt_1$ bounds the
number of isolated solutions to (\ref{prebasic}) from above
(see Fig. 2.3, c.f. Fig.8 \cite{IY2}.)

Recall that a point $x\in \R$ is  {\it{nondegenerate or regular}} for 
the map $\Delta$ if the derivative  $\Delta'(x)\neq 0$ in the 
$1$-dimensional case and $x\in \R^n$ is a regular point of a smooth 
map $F:\R^n \to \R^m$ if the rank of the linearization  $dF(x)$ at 
$x$ is maximal. 
Direct calculation shows that regular solution to
$\Delta(x_1,\eps)=x_1$ (resp. $\Delta(x_1,\eps)=x_1+\dt_1$)
corresponds to a regular solution to the system (\ref{prebasic}) 
(see lemma 3.3 \cite{IY2}).

\begin{figure}[htbp]
  \begin{center}
    \begin{psfrags}
      \psfrag{dt}{$\dt$}
      \psfrag{y=D(x,e)}{${y=\Delta(x,\eps)}$}
    \includegraphics[width= 3in,angle=0]{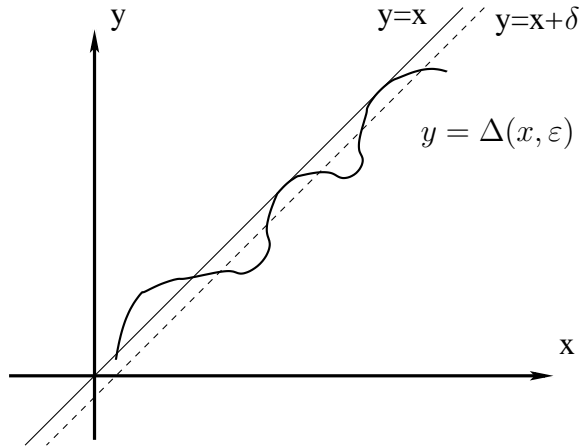}
    \end{psfrags}
    \caption{Isolated and regular solutions}
  \end{center}
\end{figure}

Moreover, if $\dt_2$ is nonzero and much smaller than $\dt_1$, then 
by the implicit function theorem the number of regular solutions of 
the system
\beq
\begin{cases} \label{basic}
y_1 = \Delta_1 (x_1, \eps), \\
x_2-f_1(y_1,\eps)=\dt_1,\\
y_2 = \Delta_2 (x_2, \eps),\\
x_1-f_2(y_2,\eps)=\dt_2 
\end{cases}
\eneq
is the same as the one for $\Delta(x_1,\eps)=x_1+\dt_1$.
Therefore, it suffices to estimate the number of small regular
solutions to the system (\ref{basic}) provided that
$1 \gg |\dt_1| \gg |\dt_2| \geq 0$.

\section{The Khovanski reduction method.}\label{asik}
\subsection{A Mixed Singular-Regular Functional System}

The system (\ref{basic}) is not easy to analyze, because it has 
the singular functions $\Delta_j$. The first key idea of the second 
step of the Main Result \cite{MR}, \cite{IY2} is to replace 
these singular equations in (\ref{basic}) by the singular 
functional-Pfaffian  equations which have polynomial differentials of 
the form (\ref{pfaff}). As a result we can obtain the {\it{mixed}} 
functional singular-regular system of the following form 
\beq
\begin{aligned} \label{funpfaf}
&
\begin{cases}
\F_1(x_1,y_1,\eps) = 0  \\
F_1(x,y,\eps)=\dt_1,\\
\F_1(x_1,y_1,\eps) = 0,\\
F_2(x,y,\eps)=\dt_2,
\end{cases}
\\
d\F_j(x_j,y_j,\eps) & =P_j(x_j,y_j,\eps)\ dx_j + 
Q_j(x_j,y_j,\eps)dy_j,\\
F_j(x,y,\eps)= &\ x_{j+1}-f_j(y_j,\eps),\ \ \ 
X=(x_1, y_1, x_2, y_2) \in (\R^{4}, 0),\ \ \  \eps \in (\R^2, 0),
\end{aligned}
\eneq
where $\F_j$ are such a functions that their differentials are
polynomial $1$-forms one of the type from column 4 of Table 1 and 
$1\gg |\dt_1|\gg |\dt_2|\geq 0$. In order to simplify considerations
we replaced functions eigenfunction $\lb_j(\eps)$ 
(see (\ref{example})) by $\eps$.
What we are interested in is the upper estimate for the number 
of small regular solutions to (\ref{funpfaf}), uniform over all 
parameters and all sufficiently small values of $\dt$'s with 
$1\gg |\dt_1| \gg |\dt_2| \geq 0$.

\subsection{Reduction of the Mixed functional system (\ref{funpfaf}) 
to Chain maps $P \circ F$}

Let $F=(F_1,F_2):\Bbb R^4 \to \Bbb R^2$ be a smooth map formed by
functions $F_1$ and $F_2$. Denote by $J^m(\Bbb R^4,\Bbb R^2)$ 
the space of $k$-jets of maps of $\Bbb R^4$ to $\Bbb R^2$. Fix 
coordinates in the source $X=(x_1,\dots , x_4)$ and the target 
$(\dt_1,\dt_2)$. Then the space $J^m(\Bbb R^4,\Bbb R^2)$ consists of
coordinates in the source, the target, and all partial derivatives
of $F$ of order at most $k$ 
\beq \label{jets}
\begin{aligned}
 \left\{(x_1,\dots, x_4);\ (F_1(X),F_2(X)); \right. \\ 
\left. \left(\frac{\partial^{\alpha} F_i}{\partial^{\alpha_1} x_1 
\dots \partial^{\alpha_n} x_4},\ \forall i=1,2,\ \alpha_j \geq 0,\ 
{\textup{such that}}\ \sum_{j=1}^4 \alpha_j \leq m \right)\right\}.
\end{aligned}
\eneq
We shall call these coordinates on the $m$-jets space 
$J^m(\Bbb R^4, \Bbb R^2)$ {\it{the natural coordinates}}.
With this coordinate system the space of $m$-jets has a natural
linear structure. We also denote by $j^mF$ the $m$-th jet of the map
$F$. Denote also by $B_r(0)\subset \R^4$ the $r$-ball centered at the
origin. We call a polynomial map $P:\R^N \to \R^n$ nontrivial if
the image $P(\R^N)$ has a nonempty interior in $\R^n$.
 
Our goal now is to realize Step 2 of our program outlined in section
\ref{program}, i.e. estimate the small number of solutions to 
(\ref{funpfaf}) via geometric multiplicity of the chain maps 
(\ref{chain}) or prove the following 

\bthm \label{Khovans} (cf. {\cite{Ka1}} Thm.10) Suppose that 
degrees of polynomial $1$-forms from (\ref{funpfaf})
are bounded by some $d\in \Z_+$. Then for a sufficiently small
$r>0$ there exists a set of $3$ (= the number of singular equations
$+1$) explicitly computable nontrivial polynomials 
$P^k=(P_1^k,\dots, P_2^k):J^2(\R^4,\R^2)\to \R^2$, $k=0,1,2$ defined 
on the space of $2$-jets $J^2(\Bbb R^4,\Bbb R^2)$ such that for a
generic $C^3$ smooth \footnote{required smoothness $3$ = the number 
of singular equations $+1$} map $F:\R^4 \to \R^2$ the number of 
regular solutions to the system (\ref{funpfaf}) inside the ball
$B_r(0)$ is bounded by the number of small regular solutions  
\beq \label{Khovappl}
\begin{aligned}
\# \{X \in B_r(0): (F_1,F_2)(X)=(\dt_1,\dt_2),
(P_1^0,P_2^0) \circ j^2F(X) =(\dt_3, \dt_4) \}\\
+\frac{1}{2} \sum_{k=1,2} \# \{X \in B_r(0): (F_1,F_2)(X)=(\dt_1,\dt_2),
(P_1^k,P_2^k) \circ j^2F(X) =(\dt_3, \dt_4)\},
\end{aligned}
\eneq
where $1\gg|\dt_1|\gg\dots \gg |\dt_4|\geq 0$ decrease to zero 
sufficiently fast. The degrees of the polynomials satisfy  
inequalities $\deg P^k_i \leq 2^i(d+1)$ for all $k$ and $i$. 
\ethm
\brm
We can not find a direct reference in the book of Khovanskii
{\cite{Kh2}}, but this Theorem is in the spirit of the results about 
perturbations discussed in section $5.2$ of this book. 
In fact this Theorem is due to Khovanskii.
\erm

\subsection{An Application of Khovanskii Method to the System 
(\ref{funpfaf}) or a Proof of Theorem \ref{Khovans}}

The method is based on the following version of Rolle's lemma

\blm \label{Rolle} 
Consider $C^2$ functions $f:S^1 \to  \Bbb R^1$ on the circle 
and $g:[0,1] \to  \Bbb R^1$ on the segment,  with a finite number 
of critical points. Then  for any $a \in \Bbb R$ and any sufficiently 
small $\dt>0$ 
\beq \label{Rol}
\beal
 \#\{ x:\ f(x)=a\} \leq & \#\{ x:\ f'(x)=\dt \} \\
\#\{ x:\ g(x)=a\} \leq & \#\{ x:\ g'(x)=\dt \}+1.
\enal
\eneq
\elm
{\it{Proof:}}\ \ One proves first the formula for $\dt=0$ using the
fact that between any two consecutive preimages there is a point where
derivative is zero. Then one uses nondegeneracy of critical points. 
Q.E.D.

Now using this lemma we shall prove Theorem \ref{Khovans}.

{\it Proof of Theorem \ref{Khovans}:}\ \  Denote by 
$\rho_r(x)=r-\sum_j x_j^2$ the function which measures distance to 
the boundary of the $r$-ball
$B_r(0)$ and vanishes on the boundary $\partial B_r(0)$.
Recall that $r$ is sufficiently small.

Denote by $G_1:\R^4 \to \R^3$ the map defined by coordinate functions
$(\F_2,F_1,F_2)$. Then the system  (\ref{funpfaf}) under investigation
becomes the map $(\F_1,G_1):\R^4 \to \R^4$, given by its coordinate
functions. In terms of this map we need to estimate the number of small 
preimages of points of the form
$$
\#\{(\F_1,G_1)\inv(0,0,\dt_1,\dt_2)\cap B_r(0)\},
$$ 
where $1\gg |\dt_1| \gg |\dt_2| \geq 0$.

Let's estimate the number of small preimages of a point 
\beq \label{mult}
\#\{(\F_1,G_1)\inv(a_1,a_2,\dt_1,\dt_2)\cap B_r(0)\},
\eneq 
where $1\gg |\dt_1| \gg |\dt_2|\geq 0$ and $a_1,\ a_2$ are arbitrary.
Since there is no restriction on $a_1,\ a_2$ the number of 
solutions may only increase.

{\it Step 1.}\ Eliminate one singular equation, say, $\F_1=0$ and
replace it by two chain-type equations $\{P^i_1 \circ j^1F\}_{i=0,1}$ 
so that for a sufficiently small $ |\dt_2|\gg |\dt_3| \geq 0$ 
the number of small regular preimages  
\beq
\# \{(G_1,P^i_1\circ j^1F(X))\inv
(a_2,\dt_1,\dt_2,\dt_3)\cap B_r(0)\}_{i=0,1}
\eneq
is at least the number of small regular preimages of 
(\ref{mult}) for any $a_1$, i.e.
\beq \label{step1}
\begin{aligned}
\#\{X\in B_r(0):\ (\F_1,G_1)(X)=(a_1,a_2,\dt_1,\dt_2)\}\leq \\
\#\{X\in B_r(0):\ (G_1,P^0_1 \circ j^1F)(X)=(a_2,\dt_1,\dt_2,\dt_3)\}+\\
\frac 12 \#\{X\in B_r(0):\ (G_1,\rho_r)(X)=(a_2,\dt_1,\dt_2,\dt_3)\}
\end{aligned}
\eneq

Consider a regular value $(a_2,\dt_1,\dt_2)\in \R^3$ of the map $G_1$.
By the Rank Theorem {\cite{GG}} the level set 
$L_{(a_2,\dt_1,\dt_2)}=G_1\inv(a_2,\dt_1,\dt_2) \cap B_r(0)\subset B_r(0)$
is a smooth $1$-dimensional manifold in the $r$-ball. 
It consists of a finite number of connected parts either compact --- 
topological circles, denoted by $\{S_i\}_{i\in I(a_2,\dt_1,\dt_2)}$, 
or noncompact --- curves $\{L_j\}_{j\in J(a_2,\dt_1,\dt_2)}$ 
reaching the boundary $\partial B_r(0)$. 
It is easy to see that 
\beq \label{step11}
\begin{aligned}
\#\{X \in B_r(0):\ (\F_1,G_1)(X)=(a_1,a_2,\dt_1,\dt_2)\} = \\
\sum_{i\in I(a_2,\dt_1,\dt_2)} \#\{X \in S_i:\ \F_1(X)= a_1\}+\\
\sum_{j\in J(a_2,\dt_1,\dt_2)} \#\{X\in L_j:\ \F_1(X)= a_1\}.
\end{aligned}
\eneq

\begin{figure}[htbp]
  \begin{center}
    \begin{psfrags}
      \psfrag{F1}{$\F_1$}
      \psfrag{dB}{$\partial B(0)$}
     \includegraphics[width= 4in,angle=0]{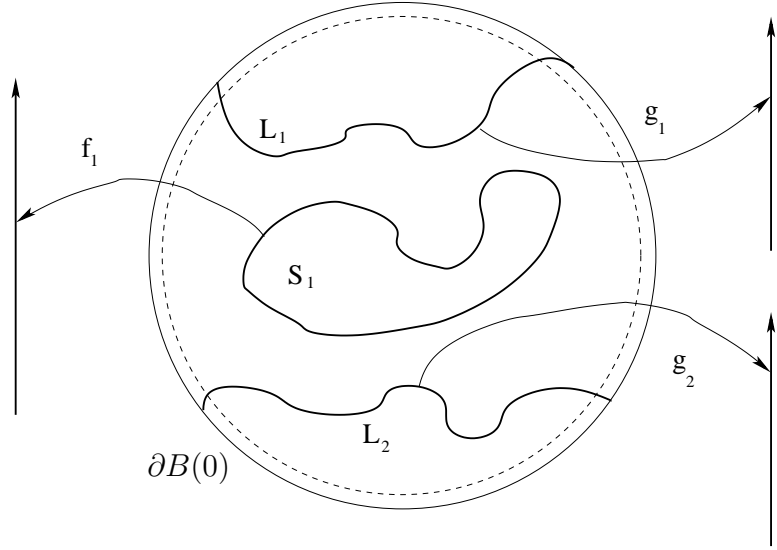}
    \end{psfrags}
    \caption{Application of Rolle's lemma}
  \end{center}
\end{figure}

Let us estimate the first sum on the right-hand side. Fix a circle,
say, $S_1$. Restrict the function $\F_1$ to $S_1$ and 
denote the result by $f_1=\F_1|_{S_1}:S^1 \to  \Bbb R$ 
(see Fig.4).  We get a function $f_i$ on the circle. 
Notice that the condition $f_1'(X)=0$ is equivalent to the
condition the Jacobian of the map $(\F_1,G_1)$, denoted by 
$J_{\F_1,G_1}(X)$, is zero.
\beq
\boxed{f_1'(X)=0} \qquad \Leftrightarrow \qquad
\boxed{J_{\F_1,G_1}(X)=0}
\eneq
Recall now that the differentials 
$d\F_j(X)=P_j\ dx_j + Q_j\ dy_j,\ j=1,2$ 
are polynomial, therefore, we have  
\beq \label{jacob1}
\begin{aligned}
J_{\F_1,G_1}(X)=*(d\F_1(X)\wedge d\F_2(X)\wedge dF_1(X) 
\wedge dF_2(X))=\\
\det(\nabla\F_1(X),\nabla\F_2(X),\nabla F_1(X),\nabla F_2(X))=
P^0_1 \circ j^1 F(X).
\end{aligned}
\eneq
where $*$ is a natural isomorphism between the space of functions 
on $\R^4$ and $4$-forms and $\nabla F(X)$ is the gradient vector of 
a function $F:\R^4 \to\R$. Since $\deg P_j,Q_j \leq d$, 
degree of $P^0_1$ is bounded by $2(d+1)$. 

Now we can apply Rolle's lemma \ref{Rolle} with $f=f_1$ and get
that for any $a_1$ and a sufficiently small $\dt_3\neq 0$
\beq\label{step12}
\begin{aligned}
\sum_i \#\{X \in S_i:\ \F_1(X)= a_1\} \leq
\sum_i \#\{X \in S_i:\ J_{\F_1,G_1}(X) =\dt_3\}.
\end{aligned}
\eneq

The second sum can be estimated in almost the same way. Instead of 
using lemma \ref{Rolle} with $f=f_i$ we need to use lemma \ref{Rolle} 
with $g=g_j=\F_1|_{L_j}:[0,1] \to  \Bbb R$ (see Fig. 2.4). Denote 
the number of component reaching the boundary $|J|$ by $k$. Then
\beq\label{step13}
\begin{aligned}
\sum_{j=1}^k \#\{X \in L_j:\ \F_1(X)= a_1 \} \leq
\sum_{j=1}^k \#\{X \in L_j:\ J_{\F_1,G_1}(X) =\dt_3\}+k
\end{aligned}
\eneq

In order to find the number of component reaching the boundary notice 
that each such component intersects the sphere $\rho_r^{-1}(\dt_3)$ for 
$\dt_3>0$ in at least two points. So $P^1_1 \circ j^1F(X)=\rho_r(X)$ 
and the second term in inequality (\ref{step1}) corresponds to 
the number of noncompact component (the boundary term).
This completes the proof of Step 1 or proves (\ref{step1}). 
 
For $i=0,1$ denote by $G^i_2:\R^4 \to \R^3$ the maps defined by its 
coordinate functions $(F_1,F_2,P^i_1\circ j^1F)$. Let's fix $i=0$ or
$1$. 
 
{\it Step 2.}\ Eliminate the second singular equation $\F_2=0$ and
replace it by two chain-type equations $\{P^i_2 \circ j^2F\}_{i=0,1,2}$ 
so that for a sufficiently small $0\leq |\dt_4| \ll |\dt_3|$ the
number of small regular preimages  
$\{(\F_2,G_2^i)\inv(a_2,\dt_1,\dt_2,\dt_3)\}_{i=0,1}$
is at least the number of small regular preimages 
$(G_2^i, P^i_2\circ j^2F)\inv (\dt_1,\dt_2,\dt_3,\dt_4)$ for any
$a_2$, i.e.
\beq \label{step2}
\begin{aligned}
\#\{X\in B_r(0):\ (\F_2,G_2^i)(X)=(a_2,\dt_1,\dt_2,\dt_3)\}
\leq \\
\#\{X\in B_r(0):\ (G_2^i,P^i_2\circ
j^2F)(X)=(\dt_1,\dt_2,\dt_3,\dt_4)\}+ 
\\ \frac 12 \#\{X\in B_r(0):\
(G_1,\rho_r)(X)=(\dt_1,\dt_2,\dt_3,\dt_4)\}
\end{aligned}
\eneq

Proof of this inequality is very similar to the proof of step 1.
We reproduce a shortened version of it in order to show why {\it the
condition $|\dt_4| \ll |\dt_3|$ it is necessary} for the Khovanskii 
method to work. 
  
Let's choose a regular value $(\dt_1,\dt_2,\dt_3)$ for the map $G_2^i$
and consider the level set 
$L_{(\dt_1,\dt_2,\dt_3)}=\left(G_2^i \right)\inv (\dt_1,\dt_2,\dt_3)$
which is by the rank theorem is a smooth $1$-dimensional manifold
consisting of a finite number of connected components either
compact--- topological circles, denoted by 
$\{S_i\}_{i\in I(\dt_1,\dt_2,\dt_3)}$, or noncompact --- curves 
$\{L_j\}_{j\in J(\dt_1,\dt_2,\dt_3)}$ reaching the boundary 
$\partial B_r(0)$. 

Then we restrict $\F_2$ to $L_{(\dt_1,\dt_2,\dt_3)}$ and get a finite
collection of functions $\{f_i=\F_1|_{S_i}:S^1 \to  \Bbb R\}_{i\in I}$
on circles and $\{g_j=\F_1|_{L_j}:[0,1] \to  \Bbb R\}_{j\in J}$
on the interval $[0,1]$. In order to use Rolle's lemma \ref{Rolle} 
we need to compute the condition $f_i'(X)=0$ (resp. $g_j'(X)=0$).
This is equivalent to the Jacobian $J_{\F_2,G_2^i}(X)$ of the map 
$(\F_2,G_2^i)$ being equal to $0$

\beq \label{jacob2}
\begin{aligned}
J_{\F_2,G_2^i}(X)=\*(d\F_2(X)\wedge dF_1(X) \wedge dF_2(X) 
\wedge d(P^i_1 \circ j^1 F)(X))= \\
\det(\nabla\F_2(X),\nabla F_1(X),\nabla F_2(X),
\nabla (P^i_1 \circ j^1 F)(X))=P^i_2 \circ j^2 F(X).
\end{aligned}
\eneq
Since $d\F_2(X)=P_2\ dx_2 + Q_2\ dy_2$ and
$\deg P_2, Q_2 \leq s$, degree of $P^i_2$ is bounded
by $4(s+1)$. An easy calculations show that each time
we take a Jacobian of a chain-map $P \circ j^k F$ its
degree at most doubles.

Now we would like to apply  Rolle's lemma \ref{Rolle} with $f=f_i$ 
(resp. $g=g_j$) and substitute a singular equation $\F_2$ by the
equation $J_{\F_2,G_2^i}(X)=\dt_4$. This equation have to be
equivalent to the fact that the derivative $f'(X)$  (resp. $g'(X)$) 
is small or covectors $\nabla \F_2(X), \nabla F_1(X),\nabla F_2(X),$ 
and $\nabla (P^i_1 \circ j^1 F)(X)$ have to be almost linear dependent. 
However, the determinant (\ref{jacob2}) can be almost zero {\it not 
because gradient vectors are almost linear dependent, but because one
of gradient vectors is small}. In order to avoid this problem 
let's make the following remark: for a fix regular value 
$(\dt_1,\dt_2,\dt_3)$ the level set $L_{(\dt_1,\dt_2,\dt_3)}$
is a smooth compact $1$-dimensional manifold possibly with a boundary
and lengths of the gradient vectors $\nabla F_1(X),\nabla F_2(X),$ and 
$\nabla (P^i_1 \circ j^1 F)(X)$ have to be bounded away from zero.
Knowing how far these lengths from zero we can choose $\dt_4$ of much 
smaller size to guarantee almost linear dependence of the gradient 
vectors. This proves inequality (\ref{step2}).

This argument allow to apply Rolle's lemma \ref{Rolle} is
the described fashion inductively in any dimension and eliminate
arbitrary number of singular Pfaffian equations. This completes 
the proof of Theorem \ref{Khovans}. See \cite{IK} and \cite{Ka1} 
for more general treatment.

\subsection{Geometric Multiplicity of Chain Maps.}
 
Let $P:\R^N \to \R^n$ be a nontrivial vector-polynomial, i.e. 
the image $P(\R^N) \subset \R^n$ has nonempty interior, 
$B^n \subset \R^n$ be a unit ball, and
$F:B^n \to \R^N$ be a generic sufficiently smooth map with 
$N \geq n$. We call the composition of a generic smooth map
and a nontrivial polynomial   
\beq \label{simchain}
P \circ F: B^n \to  \R^n
\eneq
{\it a chain map}. More generally, let $P:J^n(B^n,\R^N)\to \R^n$
be a nontrivial vector-polynomial, defined on the space of $m$-jets
for some $m \in \Z_+$. Then {\it a chain map} is
\beq \label{jetchain}  
P \circ j^m F: B^n \to  \R^n
\eneq
As the result of application of Theorem 
\ref{Khovans} to the system (\ref{funpfaf}) we need to estimate 
the number of small regular preimages of a special point of
{\it a chain map} or geometric multiplicity, defined in (\ref{geom}), 
of it. Actually application of Theorem \ref{Khovans} gives not a
generic smooth map, but a jet of a generic smooth map. To simplify 
discussion we consider the case of a smooth map (\ref{simchain}).
The general jet case can be treated using the same method.

The next two lectures are devoted to a proof of Bezout's Theorem
for chain maps. Recall that 
$B_r(0)\subset \R^n$ denotes the $r$-ball centered at the origin. 

\bthm (cf. \cite{Ka1}, Thm. 3)\label{bezout} Let $P=(P_1,\dots,P_n)$ 
be a nontrivial polynomial defined on the space of $m$-jets 
$P:J^m(B^n,\R^N)\to \Bbb R^n$ and let $F:B^n \to \Bbb R^N$ be 
a $C^k$ smooth mapping, $k>m$, and $N>n$. Suppose $F$ satisfies a
transversality condition 
depending only on $P$. Then for a sufficiently small $r$ to 
find a geometric multiplicity of the chain map (\ref{simchain})
at $0$ one can replace $j^mF$ by $L_{F,0,m}$ its linear part at $0$. 
Namely,
\beq \label{linear}
 \begin{aligned} 
 \# \{X \in B_r(0):\ 
\ P_1 \circ j^mF (X)=\dt_1, \dots , P_n \circ j^mF(X) =\dt_n \}=\\
 \# \{X\in B_r(0):\ \ P_1\circ L_{F,0,m}(X) =\dt_1,\dots , 
 P_n \circ L_{F,0,m}(X)=\dt_n\},
 \end{aligned}
\eneq
where $1\gg |\dt_1|\gg \dots \gg |\dt_n| \geq 0$. 
\ethm
\brm
By Bezout's theorem the number of isolated solutions to the equation 
in the right-hand side of (\ref{linear}) can be bounded by the product 
$\prod_{i=1}^n \deg P_i$. 
\erm

The classical transversality theorem {\cite{AGV}} says that 
a generic mapping $F$ satisfies any ahead given transversality 
condition and generic mappings form an open dense set in the space 
of smooth mapping of $B^n$ to $\R^N$. Moreover, a mapping $F$
with ``probability one'' satisfies any ahead given transversality
condition. For definitions of ``probability one'' or ``prevalence''
see \cite{HSY} and \cite{Ka7}.

{\it Acknowledgments:}\ \ I would like to thank Askold Khovanski
whose deep insight helped me to make significant simplification 
of application of Khovanski's method.

\chapter{Stratifications and Bezout's Theorem  for Chain Maps}

\medskip

In this lecture, first in section \ref{heur} we describe a geometric 
picture behind Bezout's Theorem  for chain map (Theorem \ref{bezout}) 
formulated in the end of the last lecture. It turns out that the proof 
of this Theorem reduces to a question about existence of a certain, 
so-called, $a_P$-stratification for the outer part $P$ of the chain map
(\ref{jetchain}). Then in section \ref{state} we define necessary notions 
from stratification theory, including $a_P$-stratification and discuss 
the question of existence of $a_P$-stratification. In general, it is 
not always exists as examples from \ref{examples} of Thom and Grinberg 
show.  At the end of this section we state  Hironaka's Theorem on
existence of $a_P$-stratification for polynomial functions, i.e for 
polynomial maps with $1$-dimensional image, and its extension
a Theorem on existence of $a_P$-stratification for maps with 
a mutlidimensional image proven in \cite{Ka1}.
Such a Theorem is required for the proof of the Main Result.
Finally, in section we present a geometric proof of Hironaka's Theorem
is based on the author's proof of existence of Whitney's 
stratification \cite{Ka4}. A proof of existence of Whitney's 
stratification is also outlined.

\section{An Heuristic Description} \label{heur}

Consider a chain map $P \circ F: \Bbb R^2 \to \Bbb R^2$, where 
$F:\Bbb R^2 \to \Bbb R^N$ is a generic $C^k$ smooth map, $N,\ k>2$
and $P=(P_1,P_2):\Bbb R^N \to \Bbb R^2$ is a polynomial of degree $d$. 
Fix a small positive $r$. We would like to estimate the maximal
number of small preimages 

\beq \label{numb}
\# \{x \in B_r(0):\ P_1 \circ F(x)=\eps,\ P_2 \circ F(x)=0\}
\eneq
for a small enough $\eps$.

To show the idea put $N=3$, $P_1(x,y,z)=x^2+y^2$, and
$P_2(x,y,z)=xy$. Assume also that $F(0)=0$. 
Denote the level set by $V_\eps=\{P_1=\eps,\ P_2=0\}$. The level set
$V_\eps$ for $\eps>0$ consists of $4$ parallel lines (see Fig. 3.1).

Notice that in our notation the number of intersections of
$F(B_r(0))$ with $V_\eps$ equals the number of preimages of the point 
$(\eps,0)$ under $P \circ F$ see (\ref{numb}).

It is easy to see from Fig. 3.1 that if $F$ is transverse to 
$V_0$ it is transverse to $V_\eps$ for any small $\eps>0$. Moreover,
the number of intersections $F(B_r(0))$ with $V_\eps$ equals 4
(see the points $P_1,\dots,P_4$).

Another way to calculate the same number is
as follows. Let us replace $F$ by its linear part $L_F$ at zero. Then 
\beq \nonumber 
\beal
\#\{x \in B_r(0):\ P_1 \circ F(x)=\eps,\ P_2 \circ F(x)=0\}=\\
\#\{x \in B_r(0):\ P_1 \circ L_F(x)=\eps,\ P_2 \circ L_F(x)=0\}
\enal
\eneq
and solving this polynomial system also yields $4$.

\begin{figure}[htbp]
  \begin{center}
    \begin{psfrags}
      \psfrag{Ve}{$V_\eps$}
      \psfrag{V0}{$V_0$}
      \psfrag{P1}{$P_1$}     
      \psfrag{P2}{$P_2$}
      \psfrag{P3}{$P_3$}
      \psfrag{P4}{$P_4$} 
      \psfrag{F(B)}{$F(B_r(0))$}            
     \includegraphics[width= 3in,angle=0]{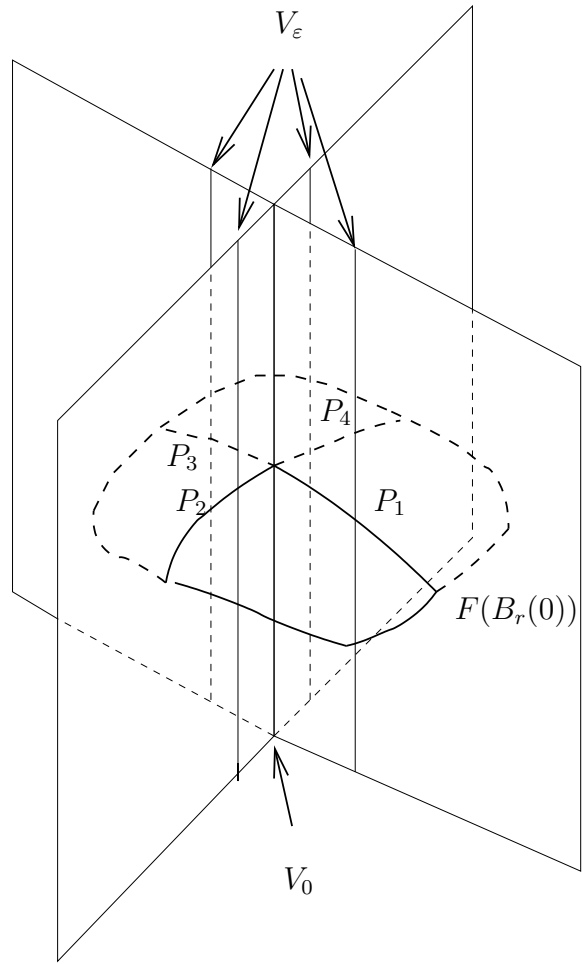}
    \end{psfrags}
    \caption{The Idealistic Example}
  \end{center}
\end{figure}

The idea behind this picture is the following: Consider an arbitrary
$N$ and a polynomial $P=(P_1,P_2):\Bbb R^N \to \Bbb R^2$ of degree at
most $d$ and $N>2$. Define the algebraic variety 
$V_\eps=(P_1,P_2)^{-1}(\eps,0)$ as the level set.

Assume for simplicity that for any small $\eps \neq 0$ the level 
set $V_\eps$ is a manifold of codimension 2. We shall get rid of 
this assumption later (see Theorem \ref{existence} b). 
It turns out that there exists a special partition 
$\Cal V_0=\{V_i\}_{i\in \Cal A}$ of $V_0=\bigsqcup_{i\in \Cal A} V_i$
into semialgebraic parts which are attached to their neighbors 
``regularly'' see definition \ref{strdef}) such that 
it depends on $P$ only and satisfies the following condition. 
We say $F$ is transverse  to a stratified set $(V_0,\Cal V_0)$
if $F$ is transverse to each stratum $V_i \in \Cal V_0$, then

\begin{multline}\label{apstrat}
\boxed{F\ \text{is transverse to}\ (V_0,\Cal V_0)}
\implies \boxed{F \ \text{is transverse to}\  V_\eps}
\qquad\ 
\end{multline}

\blm \label{lin} Let $B_r(0)$ be the $r$-ball centered at the point 
$0 \in \Bbb R^2$ and let $L_{F,0}$ denote the linearization of $F$ at
the point $a$. Under condition (\ref{apstrat}), the number of 
intersections of the image $F(B_r(0))$ with $V_\eps$ coincides with 
the number of intersections of the image $L_{F,0}(B_r(0))$ with
$V_\eps$, provided $r$ is small enough. That is  
\beq \label{replace}
\begin{aligned}
\#\{x \in B_r(0):\ (P_1,P_2) \circ F(x)=(\eps,0)\}=\\
\#\{x \in B_r(0):\ (P_1,P_2) \circ L_{F,0}(x)=(\eps,0)\}.
\end{aligned}
\eneq
\elm
\brm
The argument below is independent of codimension of $V_\eps$.
We only need condition (\ref{apstrat}) and the fact that codimension 
of $V_\eps$ coincides with dimension of the preimage of a chain map 
$P \circ F$.
\erm 
{\it{Proof:}}\ \ Consider the $1$-parameter family of maps 
$F_t=tF+(1-t)L_F$ deforming the linear part of $F$ into $F$. 
Clearly, $F_1\equiv F$ and $F_0 \equiv L_F$. Fix a small $r>0$.
Since, $F$ is transverse to $V_0$ at $0$ all $F_t$ are 
transverse to $V_0$ at $0$. Condition (\ref{apstrat})
implies that for all small $\eps$ and all $t \in [0,1]$ we have
$F_t$ is transverse to $V_\eps$. 

Therefore, {\it{the number of intersections of $F_t(B_r(0))$ with 
$V_\eps$ is independent of $t$.}} Indeed,
assume that $\#\{F_{t_1}(B_r(0)) \cap V_\eps\} \neq  
\#\{F_{t_2}(B_r(0)) \cap V_\eps\}$ for some $t_1<t_2$.
Then as $t_1$ increases to $t_2$ there is a point $t^*$ where the
number of intersections either drops or jumps.
At this point $t^*$ the condition of transversality of
$F_{t^*}$ and $V_\eps$ must fail.
This completes the proof of the lemma. Q.E.D.

\section{Basic definitions of stratified sets, maps, and etc}
\label{state}

\subsection{Stratified sets}

A stratification of a set, e.g. an analytic variety, is, 
roughly, a partition of it into manifolds so that these manifolds
fit together ``regularly''. Stratification theory was originated 
by Thom \cite{Th1} and Whitney \cite{Wh} for algebraic and analytic 
sets. It was one of the key ingredients in Mather's proof of
the topological stability theorem {\cite{Ma}}.  
For the history and further applications of stratification theory 
see {\cite{GM}} and {\cite{PW}}.

We consider here only the category of real (semi)algebraic 
sets for simplicity. Theorems on existence of stratifications proven
here in the category of semialgebraic sets can be proven for 
the categories of complex or real (semi)analytic sets using similar
methods. Call a subset $V \subset \R^N$ 
a {\it{semivariety}} if locally at each point $x\in \R^N$ it is 
a finite union of subsets defined by equations and inequalities
\beq
f_1= \dots =f_k=0  \quad  g_1>0, \dots, g_l> 0
\eneq
where $f_i$'s and $g_j$'s are real algebraic depending on. 
Semivarieties are closed under Boolean operations. 

\bdef  \label{regular} (Whitney)
Let $V_i,V_j$ be disjoint manifolds in $\R^N$,
$\dim V_j>\dim V_i$, and let $x \in V_i \cap \overline{V_j}$. 
A triple $(V_j,V_i,x)$ is called $a$(resp. $b$)- regular if 

$a$) when a sequence $\{y_n\} \subset V_j$ tends to $x$
and $T_{y_n} V_j$ tends in the Grassmanian bundle to 
a subspace $\tau_x$ of $\R^N$, then 
$T_x V_i\subset \tau_x$;

$b$) when sequences $\{y_n\} \subset V_j$ and 
$\{x_n\} \subset V_i$ each tends to $x$, the unit vector 
$(x_n-y_n)/|x_n-y_n|$ tends to a vector $v$, and 
$T_{y_n}V_j$ tends to $\tau_x$, then 
$v\in \tau_x$\footnote{This way of defining $b$-regularity is 
due to Mather \cite{Ma}. Whitney's definition \cite{Wh}
is equivalent to this one provided of $a$-regularity}.

$V_j$ is called $a$(resp. $b$)- regular over 
$V_i$ if each triple $(V_j,V_i,x)$ is $a$(resp. $b$)- regular.
\endef
\brm
Since the Grassmanian manifold of $\dim V_j$-panes in
$m$-dimensional space is compact, existence of limits in 
the definition above can be reached by choosing a subsequence
$\{x_{n_k}\}_{k \in \Z_+}$ or $\{y_{n_k}\}_{k \in \Z_+}$ if necessary. 
\erm

\bdef \label{strdef} (Whitney) Let $V$ be a semivariety in 
$\R^N$.  A disjoint decomposition 
\beq\label{decomp}
V=\bigsqcup_{i \in I} V_i, \quad V_i\cup V_j =\emptyset \ \ \ 
\textup{for} \ \ \ i\neq j  
\eneq
into smooth semivarieties $\Cal V=\{V_i\}_{i \in I}$,
called strata, is called an $a$(resp. $b$)-regular 
stratification if

1. each point has a neighborhood intersecting only finitely 
many strata;

2. the frontier $\overline{V_j}\setminus V_j$ of each stratum $V_j$ 
is a union of other strata $\bigsqcup_{i \in J(i)} V_i$;

3. any triple $(V_j,V_i,x)$ such that $x \in V_i \subset
\overline{V_j}$ is $a$(resp. $b$)-regular.
\endef

The classical example of a stratified algebraic set in $\R^3$ is 
so-called {\it Whitney umbrella}. It is defined as follows 
\bex \label{umbr}
Consider the $2$-dimensional algebraic variety in $\R^3$, defined by
\beq
V=\{(x,y,z)\in \R^3:\ y^2=zx^2\}.
\eneq 
\enex

\begin{figure}[htbp]\label{umbrella}
  \begin{center}
    \begin{psfrags}
     \includegraphics[width= 3in,angle=0]{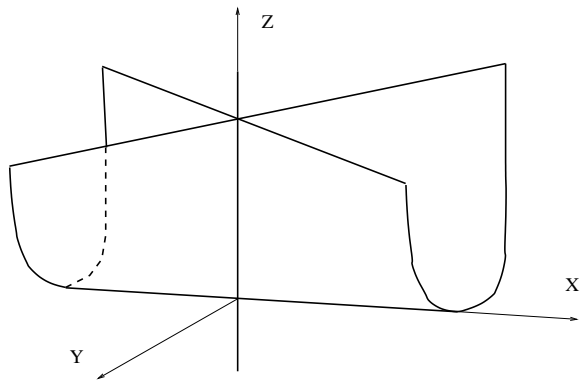}
    \end{psfrags}
    \caption{Whitney's umbrella}
  \end{center}
\end{figure}

The first natural partition of $V$ into smooth parts (strata) is the 
vertical line $V_1=\{x=y=0\}$ and the complement $V_2=V\setminus V_1$.
However, $V_2$ does not fit regularly to $V_1$ at the origin.
To see that consider the sequence of the form $(x_n,0,0)\in V_2$ with  
$x_n\to 0$ as $n\to \infty$. It is easy to see that after we refine 
$V_1$ into $V_{0'}=\{0\}$ and $V_{1'}=V_1\setminus V_0$ and put 
$V_2=V_{2'}$ the partition $V=\bigsqcup_{i'=1,2,3} V_{i'}$ becomes 
the stratified manifold $(V,\{V_{0'},V_{1'}, V_{2'}\})$.

\bthm \label{main} \cite{Wh}, \cite{Th2}, {\cite{Lo1}} For any 
semivariety  $V$ in $\R^N$ there is an $a$(resp. $b$)-regular
stratification.
\ethm
\brm
This Theorem is not true for smooth sets. To see that one can
construct a $2$-surface in $\R^2$ which looks like a corkscrew.
\erm
Existence of stratifications in the complex analytic case was proved 
by Whitney {\cite{Wh}}. Later Thom published a sketch of 
a proof {\cite{Th2}}. Then Lojasiewicz {\cite{Lo1}} extended 
these results to the semianalytic case. Later other proofs were 
found. In \cite{Hi} Hironaka found a nice proof using his
resolution of singularities. J.Bochnak, M. Coste,
and M.-F. Roy \cite{BCR} and Z. Denkowska, K. Wachta \cite{DeW}
follow the classical route of Whitney \cite{Wh} via the wing lemma. 
\cite{BCR} uses a Nash wing lemma and \cite{DeW} apply the parameterized 
Puiseux Theorem of W. Pawlucki \cite{Pa}. T. Wall \cite{Wa} and 
S. Lojasiewicz, J. Stasica, and K. Wachta \cite{LSW} found proofs 
which use on Milnor's curve selection lemma \cite{Mi1}. 
The latter proof also uses the subanalyticity of the tangent map 
( for which an elementary proof was given by Z. Denkowska and 
K. Wachta \cite{DeW}). In \cite{Ka4} the author gives a geometric 
proof based on a simple observation that 
regularity of stratifications is related to uniqueness of the limit 
of the tangent planes to a bigger stratum as they approach to a
smaller stratum. This proof is outlined in section \ref{whit}.
For a nice exposition of the theory of semianalytic 
and subanalytic sets see \cite{Lo2}.

\subsection{Stratified maps and $a_P$-stratification}

First we define a smooth map of a stratified set $(V,\Cal V)$:

\bdef \label{stratdef}
Let $(V,\Cal V)$ be stratified in $\R^N$, $V \subseteq \R^N$, then 
a map $P:V \to \R^n$ is called $C^2$-smooth if it can be extended to 
a $C^2$ smooth map of an open set $V \subset U \subset \R^N$, denoted 
by ${\bf P}:U \to \R^n$ whose restriction to $V$ coincides with $P$.

A stratification $V=\cup_i V_i$ stratifies a smooth 
map $P:V \to \Bbb R^n$ if the restriction of $P$ to a stratum 
$V_i$ has a constant rank or rank $dP|_{V_i}(x)$ is 
independent of $x \in V_i$.

A map $G:\R^n \to \R^N$ is called transverse to a stratified set  
$(V,\Cal V)$ if $G$ is transverse to each strata $V_i \in \Cal V$.
\endef

\bex With the notations of the example \ref{umbr} of the Whitney 
umbrella consider the Whitney umbrella $V$ and the projection 
$P=\pi|_{V}:V \to \R^2$ along the $z$-coordinate restricted to it. 
Then the stratification $(V,\{V_{0'},V_{1'}, V_{2'}\})$ of $V$
is also an $a_P$-stratification of $P:V \to \R^2$. 
\enex

Let $V_i$ be a stratum of a stratification $(V,\Cal V),\ 
\Cal V=\{V_i\}_{i \in \A}$ and $a\in N$. Denote by 
$L_{a,i}=(P\inv(a) \cap V_i)$ the level sets of $P$ in $V_i$. 
By the Rank Theorem {\cite{GG}}, if a stratification $(V,\Cal V)$ 
stratifies a smooth map $P$, then for each strata $V_i$ the number  
$d_i(P)=\dim V_i- {\textup{rank}}\ dP|_{V_i}$
is well defined and equals to dimension of any nonempty level 
set $L_{a,i}$. 

Roughly speaking, $a_P$-stratification is a stratification of a map 
$P:V \to N$ such that it is also an $a$-stratification of its level
sets, i.e. for any sequence of points $\{b_k\} \subset P(V_j)$
converging to a point $a \in P(V_i)$ the corresponding level sets 
$L_{b_k,j}=(P\inv(b_k) \cap V_j)\subset V_j$ approach the limiting 
level set $L_{a,i}\subset V_i$ ``regularly''. A precise definition 
is the following

\bdef \label{a-p} Let $P:\R^N \to \R^n$ be a $C^2$ smooth map and 
let $V_j$ and $V_i$ be submanifolds of $\R^N$ such that 
$V_i \subset \overline{V_{j}}\setminus V_{j}$ the restrictions 
$P|_{V_j}$ to $V_j$ and $P|_{V_i}$ to $V_i$ have constant ranks. $V_j$ 
is called $a_P$-regular over $V_i$ with respect to the map $P$ at a
point $x\in V_i \cap \overline{V_j}$ if for any sequence 
$\{x_n\} \subset V_j$ converging to $x$ the sequence of tangent planes 
to the level sets $T_n=ker\ dP|_{V_j}(x_n)$ converges in 
the corresponding Grassmanian manifold of 
$\dim ker\ dP|_{V_j}$-dimensional planes to a plane $\tau$ and 
\beq \label{apcond}
\lim ker\ dP|_{V_j}(x_n)=\tau \supseteq ker\ dP|_{V_i}(x)
\eneq
\endef

\bdef \label{apstratif} (Thom)
A $C^2$ smooth map $P:V \to \R^n$ of a stratifiable set $V$ to $\R^n$
is called $a_P$-stratifiable if there exists a stratification $(V,\Cal V)$ 
such that the following conditions hold:

a) $(V,\Cal V)$ stratifies the map $P$ (see definition \ref{stratdef});

b) for all pairs $V_j$ and $V_i$ from $\Cal V$ such that 
$V_i \subseteq \overline{V_j} \setminus V_j$ the stratum $V_j$ 
is $a_P$-regular over the stratum $V_i$ with respect to $P$ at 
point $x$ for all $x \in V_i \cap \overline{V_{j}}$.
\endef

{\bf Remarks} \ 
1. The original definition of $a_P$-stratification requires 
an appropriate stratification of the image too {\cite{Ma}}, 
but for simplicity we do not require existence of stratification of
the image and it turns out to be sufficient for our purposes.

2. With the notations above for an $a_P$-stratification to exist we
must have $d_i(P) \leq d_j(P)$ for each 
$V_i \subseteq \overline{V_j} \setminus V_j$, i.e. nonempty level 
sets $L_{b,j}$ inside the bigger stratum $V_j$ have dimension
$d_j(P)$ greater or equal to dimension $d_i(P)$ of nonempty 
level sets $L_{a,i}$ in the smaller stratum $V_i$.
Otherwise $\dim ker\ dP|_{V_j}(x_n)<\dim ker\ dP|_{V_i}(x)$
and  condition (\ref{apcond}) can't be satisfied.

\subsection{For $a_P$-stratifications condition (\ref{apstrat})
holds.}

An heuristic description given above shows that the key to a proof 
of Bezout's Theorem is condition (\ref{apstrat}) (lemma \ref{lin}). 
Now we prove that existence of an $a_P$-stratification of 
a polynomial $P$ is sufficient for condition (\ref{apstrat}) to hold.
 
Let $P=(P_1,P_2):\Bbb R^N \to \Bbb R^2$ be a nontrivial polynomial, 
i.e. the image $P(\R^N)$ has nonempty interior. Denote by
$V=P_2^{-1}(0)$ and  $V_0=(P_1,P_2)^{-1}(0)$ the level sets. Assume 
that there exists a stratification $(V,\V)$ which stratifies 
the map $P|_V$ and such that the zero level set $V_0$ is also 
stratified by a stratification $(V_0,\Cal V_0)$ with  
$V_0=\bigsqcup_{i \in \A_0} V_i$

\blm \label{app} With the above notation if each stratum 
$V_j \in \V \setminus \Cal V_0$ is $a_P$-regular over each 
stratum $V_i \in \Cal V_0$ with respect to the polynomial $P$, 
then any $C^2$ smooth map $F: \Bbb R^2 \to \Bbb R^N$ which is 
transverse to $(V_0,\Cal V_0)$ is also transverse to 
each level set $\V_{b,j}$ for any small $b$ and this is equivalent 
to condition (\ref{apstrat}). 
\elm

{\it{Proof:}}\ \ Pick a point $x$ in $V_i \subset V_0$ and a point 
$y \in V_i$. Notice that $ker\ dP|_{V_i}(x)$ is the tangent 
plane to the level set $\{P^{-1}(P(x)) \cap V_i\}$ 
at the point $x$ and $ker\ dP|_{V_j}(y)$ is the tangent plane
to the level set $\{P^{-1}(P(y)) \cap V_j\}$.

By condition (\ref{apcond}) if a map $F:X \to \Bbb R^N$ is transverse 
to $ker\ dP|_{V_i}(x)$ at a point $x$, then $F$ is transverse to
$ker\ dP|_{V_j}(y)$ for any $y \in V_j$ nearby $x$.
This completes the proof of the lemma. Q.E.D.

\section{Existence of $a_P$-stratifications 
for polynomial maps}\label{polynom}

\subsection{Examples of nonexistence due to Thom and Grinberg}
\label{examples}

Existence of $a_P$-stratifications is a nontrivial question.
There are some obvious obstacles. For example, let $V \subset \R^N$ 
be an algebraic variety and let $P:\Bbb R^N \to \Bbb R^n$ be 
a polynomial map. Assume that $(V,\Cal V)$ stratifies $P$.
Take two strata $V_i$ and $V_j$ so that $V_j$ lies ``over''
$V_i$, i.e. $V_i\subseteq \overline{V_j} \setminus V_j$, then 
condition (\ref{apcond}) can't be satisfied if dimension of the level 
sets $d_i(P)$ in the upper stratum $V_i$ is strictly less than that
of $d_j(P)$ in the lower stratum $V_j$, i.e., 
$\dim ker\ dP|_{V_i}(y)< \dim ker\ dP|_{V_j}$. In this case a
plane $ker\ dP|_{V_j}(x)$ of the lower stratum $V_j$ should belong 
to a plane $\tau$ of smaller dimension by condition (\ref{apcond}), 
which is impossible. Thom constructed the first example when 
this happens {\cite{GWPL}}.

\bex (Thom) Consider the vector-polynomial $P$ of 
the form 
\beq
P:{x\choose y} \to {x\choose xy}.
\eneq 
The line $\{x=0\}$ is the line of critical points of $P$. Outside of 
the line $\{x=0\}$ the map $P$ is a diffeomorphism, therefore, the
preimage $P^{-1}(a)$ of any point $a \neq 0$ is $0$-dimensional. 
On the other hand, the preimage of $0$ is the $1$-dimensional line 
$\{x=0\}$. Thus, $a_P$-regularity  fails to exist.
\enex

\bdef \label{compat} Let us call an algebraic set $V$ rank compatible 
with respect to a polynomial $P$ if there exists a stratification 
$(V,\Cal V)$ which stratifies $P$ and for any pair $V_i$ and
$V_j$ from $\Cal V$ such that $V_i \subseteq \overline{V_j} 
\setminus V_j$ dimension of the level sets $d_i(P)$ in the lower 
stratum $V_i$ does not exceed dimension of the level sets
$d_i(P)$ in the upper stratum $V_j$.
\endef 

It turns out that even if an algebraic set $V$ is rank compatible 
with respect to a polynomial $P$, then $a_P$-stratification still 
does not always exist. Let us present an example with this property
due to M. Grinberg. It seems that the existence of a counterexample 
was known before, but we did not find an appropriate reference.

Let $V=\{(x,y,z,t)\in \R^4:\ x^2=t^2y+z\}$ be the three-dimensional 
algebraic variety and $P:V\to \R^2$ be the natural projection to 
the last two coordinates, i.e. $P:(x,y,z,t)\to (z,t)$.

\blm \label{Grinb} With the above notations the set $V$ is 
rank compatible with the polynomial map $P$ and does not 
have an $a_P$-stratification.
\elm

{\it{Proof:}} \ \ Consider a rank stratification of $V$.
Such a stratification consists of three stratum:
$V_1=\{x=t=z=0\},\ V_2=\{t=0,\ x^2=z, \ x \neq 0\},$ and 
$V_3=\{t \neq 0\}.$  On each stratum  rank$P|_{V_i}=i-1$.
The level sets $P^{-1}(t,z)$ are parabolas for $t \neq 0$ 
and lines for $t=0$. 

Show that for each point ${\bf{a}}=(0,a,0,0) \in V_1$
there exists a $1$-parameter family of level sets such that 
at the point ${\bf a}$ the property $a_P$-regularity of $V_3$ 
over $V_1$ fails. 

Consider the preimage  of the curve $\{z=-at^2\} \subset \R^2$.
This is an algebraic  variety of the form $W_a=\{x^2=t^2(y-a)\}$.
One can see that $W_a$  is the Whitney umbrella (see Fig. 3.2).
The level $x^2=t_0^2(y-a)$ is the parabola. As $t_0 \to 0$ 
this parabola tends to semiline $\{x=t=z=0,\ y \geq a\}$. At 
the point ${\bf{a}} \in V_1$ the property $a_P$-regularity of 
$V_3$ over $V_1$ clearly fails. This completes 
the proof of the lemma. Q.E.D.

\subsection{Existence of $a_P$-stratifications (Hironaka's 
Theorem and its extension)} \label{apstratific}

As we have seen above sometimes $a_P$-stratifications exist, 
sometimes they are not. Let us state a positive result 
on existence of it.

\bthm (Hironaka) \cite{Hi} \label{Hi}
If $V \subset \Bbb R^N$ is a semialgebraic variety and
$P:\Bbb R^N \to \Bbb R$ is a polynomial function, then 
there exists an $a_P$-stratification $(V,\Cal V)$ of $V$
with respect to $P$ by semialgebraic strata.
\ethm

In the next section we give a geometric proof of this result 
based on the proof of existence of Whitney's stratifications 
due to the author \cite{Ka4}. Below we describe an extension on
Hironaka's Theorem to maps with a multidimensional image
proven in \cite{Ka1}. This extension is sufficient to prove 
Bezout's Theorem for chain maps (Theorem \ref{bezout}).

{\bf The Tarsky-Seidenberg Principle} {\it (see e.g. {\cite{BCR}}, 
{\cite{Ja}}) For any semialgebraic $V$ in $\R^N$ and a polynomial 
map $P:\R^N \to \R^n$ the image $P(V)$ is semialgebraic.}

Let $\Bbb R^N$ and $\Bbb R^n$ be Eucledian spaces with the 
fixed coordinate systems $x=(x_1,\dots, x_N) \in \Bbb R^N$ and 
$a=(a_1,\dots, a_k) \in \Bbb R^n$ with $N \geq n$ and a 
non-trivial vector-polynomial $P:\Bbb R^N  \to \Bbb R^n$.
Recall that $P$ is a nontrivial if the image $P(\R^N)$ has nonempty
interior. In what follows we call {\it{vector-polynomial}} 
by {\it{polynomial}} for brevity.  

\bdef
Let $m \in \Z_+$ and $\dt > 0$. We call the $(m,\dt)$-cone 
$K^n_{m,\dt}$ the following set of points
\beq \label{cone}
\begin{aligned}
K^n_{m, \dt} = \{a=(a_1, \dots, a_N) \in 
\R^N:\  0 < |a_1| < \dt,\\
0 < |a_{j+1}| < |a_1 \dots a_j|^{m}\ 
\text{for}\ j=1, \dots, N-1\}.
\end{aligned}
\eneq
Let $m'\in \Z^N_+$. If $m' \geq m$ and $\dt'\leq\dt$, then we say 
that the $(m', \dt')$-cone $K^n_{m',\dt'}$ is a refinement of the
$(m,\dt)$-cone $K^n_{m,\dt}$.
\endef

Define the following sets
\beq \label{limset}
V_{m, \dt, P} = 
\textup{closure} \{ P^{-1}(K^n_{m, \dt})\},\  
V_{0, m, P} = \cap_{\dt > 0} V_{m, \dt, P}
\eneq
Then one has 
\bthm \cite{Ka1} \label{existence} For any nontrivial polynomial 
$P:\R^N \to \R^n$ there exist an integer $m \in \Bbb Z_+$ 
and a positive $\dt$ such that the following conditions hold

a) the set  $V_0=V_{0,m,P}$ (see (\ref{limset})) is a
semialgebraic set of codimension  at least $n$. 

b) the set  $V_{{m},\dt,P}$ consists of regular 
points of $P$, i.e. if $b \in V_{m,\dt,P}$, then 
the level set $P^{-1}(P(b))$ is a manifold of codimension $n$.

c) there exists a stratification of $V_0$ by semialgebraic 
strata $(V_0,\Cal V_0)$ satisfying the property: 
$V_{m,\dt,P}$ is $a_P$-regular over any strata 
$V_i \in \Cal V_0$ with respect to $P$.
\ethm

\brm In order to have compatibility condition for the limiting 
set $V_{m, \dt, P}$ with regular level sets $P\inv(a)$ in 
the definition of the $(m,\dt)$-cone $K^n_{m, \dt}$ (\ref{cone}) 
it is necessary that range of values (``smallness'') of $a_{j+1}$ 
depends on all $a_i$'s with $i=1,\dots,j$. Indeed, consider 
the following  
\erm
\bex
Let $x=(x_1,x_2,x_{3})$ denote a point in $\R^3$ and 
$P=(P_1,P_2, P_3):\R^3\to \R^3$ be a polynomial map, given by 
\beq
P_1(x)=x_1,\ P_2(x)=x_1x_2,\ 
P_3(x)=x_1 x_2 x_3.
\eneq 
If definition of the $(m,\dt)$-cone is
\beq \label{cone'}
K^3_{m, \dt} = \{a=(a_1, a_2, a_3) \in \R^3:\  0 < |a_1| < \dt,
\ 0 < |a_{2(\textup{resp.} 3)}| < |a_1|^{m}\},
\eneq
then the limiting set $V_{m, \dt, P}$, defined by (\ref{limset}),
is $1$-dimensional for any positive $m$. However, all level sets 
$P\inv(a)$ with $a_1a_2\neq 0$ are $0$-dimensional. In this case
compatibility condition \ref{compat} fails.
\enex

\section{A Proof of Hironaka's Theorem on existence of 
$a_P$-stratifications for polynomial functions}

In this section we present a geometric proof of Hironaka's
Theorem based on a proof of Whitney's Theorem \ref{main} on 
existence of $a$-stratifications due to the author \cite{Ka4}.
First, we briefly outline the latter proof and then prove
Hironaka's Theorem following the same path.

\subsection{An Outline of a Proof of Whitney's Theorem \ref{main} 
on existence of $a$-stratifications}\label{whit}

The outline given below works to prove $b$-stratifications
too after  a slight modification \cite{Ka4}. 

A semivariety $V$ has well-defined dimension, say $d \leq N$. 
Denote by $V_{reg}$ the set of points, where  $V$ is locally 
a real algebraic submanifold of $\R^N$ of dimension $d$. $V_{reg}$ 
is a semivariety, moreover, $V_{sing}=V \setminus V_{reg}$ is 
a semivariety of positive codimension in $V$, i.e. 
$\dim V_{sing}< \dim V$. In the algebraic case they are not 
difficult (see e.g. {\cite{Mi1}}).

{\it{Step 1.}} There is a filtration of $V$ by semivarieties
\beq
V^0 \subset V^1 \subset \dots \subset V^d=V,
\eneq
where for each $k=1, \dots, d$ the set $V^k \setminus V^{k-1}$ 
is a manifold of dimension $k$. This is not difficult see e.g.
\cite{Mi1}. Indeed, consider $V_{sing} \subset V$, then 
$V\setminus V_{sing}$ is a manifold of dimension $d$ and 
$\dim V_{sing}<d$. Inductive application of these arguments
completes the proof.

A {\it{refinement}} of a decomposition $V=\bigsqcup_{i\in I}V_i$
is a decomposition $V=\bigsqcup_{i'\in I'}V_{i'}$ such that 
any stratum $V_j$ of the first decomposition is a union of some 
strata of the second one, i.e. there is a set
$I'(j) \subset I'$ such that $V_j=\bigsqcup_{i' \in I'(j)}V_{i'}$.  
 
{\it{Step 2}}. Let $V\subset \R^N$ be a manifold and $W\subset V$ be 
a semivariety. Denote by $Int_V(W)$ the set of interior points of 
$W$ in $V$ w.r.t. the induced from $\R^N$ topology. Let $V_i$ and 
$V_j$ be a pair of distinct strata. For each point 
$x\in V_i\cap \overline{V_j}$ denote by $V_j^{con,x}$ a local 
connected component of $V_j$ at $x$, i.e. a connected component of 
intersection of $V_j$ with a small ball centered at $x$ and call it 
{\it essential} if the closure of $V_j^{con,x}$ has $x$ is in 
the interior, $x\in Int_{V_i}(V_i\cap \overline{V^{con,x}_j})$. 
Denote by $V_j^{ess,x}$ the union of all local essential components 
of $V_j$. A semialgebraic set $V_j$ can have only a finitely 
many local connected components (see e.g. \cite{Mi1}). 
 
\bthm \label{sing} For any two disjoint strata $V_j$ and 
$V_i$ the set of points 
\beq \nonumber
Sing_{a}(V_j,V_i)=
\{ x \in V_i\cap\overline{V_j}:\ (V_j^{ess,x},V_i,x) \ 
\textup{is not}\ a -\textup{regular}\},
\eneq
is a semivariety in $V_i$ and 
$\dim  Sing_{a}(V_j,V_i)<\dim V_i$.
\ethm
Let us show that this Theorem is sufficient to prove the $a$-regular
case of Theorem \ref{main}.  Consider a decomposition 
$V=\bigsqcup_{i\in I}V_i$ and
split the strata into two groups: the first group consists of strata 
of dimension at least $k$ and the second group is of the rest. 
Suppose that each stratum from the first group is $a$-regular over 
each stratum from the second group. Then by definition of 
$a$-regularity any refinement of a stratum from the second group 
preserves this $a$-regularity. 

Now apply this refinement inductively. Consider strata in 
$V^d\setminus V^{d-1}$ of dimension $d$. Using Theorem \ref{sing} and
the result of Lojasiewicz \cite{Lo1} that a frontier of a semivariety 
has dimension less than a semivariety itself, refine $V^{d-1}$ so 
that each $d$-dimensional stratum is $a$-regular over each stratum 
in $V^{d-1}$. The above remark shows that any further refinement of the 
strata in $V^{d-1}$ preserves the $a$-regularity of strata from 
$V^d\setminus V^{d-1}$ over it. This reduces the problem of 
existence of stratification for $d$-dimensional semivarieties to 
the same problem for $(d-1)$-dimensional semivarieties. Induction on 
dimension completes the proof of Theorem \ref{main}.

Our proof is based on the observation that if $V_i\subset
\overline{V_j}$ are a pair of strata $a$-regularity of
$V_j$ over $V_i$ at $x$ in $V_i$ is closely related to whether the
limit of tangent planes $T_yV_j$ is unique or not as $y$ from $V_j$ 
tends to $x$. The rest of the paper is devoted to the proof of Theorem
\ref{sing} which consists of two steps. In lemma \ref{one}
we relate $a$-regularity with (non)uniqueness of limits of tangent
planes $T_yV_j$, then based on it and Rolle's lemma in lemma \ref{two}
one can prove Theorem \ref{sing}.

Let $V_i$ and $V_j$ be a pair of distinct strata in $\R^N$. Define
\beq 
\begin{aligned}\label{uniquelim}
Un(V_j,V_i)=\{ x \in V_i\cap \overline{V_j}&:
\textup{for any}\ V_j^{con,x},\ \textup{there exists}\ 
\tau_x\subset T_x\R^N\\ 
\textup{such that for any}\ & \{y_n\} \subset V_j^{con,x}\ 
\textup{tending to}\  x, \ \  T_{y_n}V_j \to \tau_x \},
\end{aligned}
\eneq

The proof consists of two lemmas.

\blm \label{one}
With the above notations we have 
\beq
Sing_a(V_j,V_i) \subset V_i\setminus Un(V_j,V_i).
\eneq
\elm

\blm \label{two} With notations above there is a set of strata 
$\{V_j^p\}_{p\in \Z}$ (resp. $\{V_i^p\}_{p\in \Z}$) in $V_j$ 
(resp. in $V_i$) each of positive codimension in $V_j$ (resp. in
$V_i$) such that
\beq \label{include}
\beal
Sing_a(V_j,V_i) & \subset \ \ 
\bigcup_{p\in \Z} Sing_a(V_j^p,V_i)
\bigcup_{p\in \Z} V_i^p. \\
V_i\setminus Un(V_j,V_i) & \subset \ \ 
\bigcup_{p\in \Z} V_i\setminus Un(V_j^p,V_i)
\enal
\eneq
\elm
{\bf Remarks.} \ 
1. Inductive application of this lemma to the right-hand side of
the first line of (\ref{include}) reduces dimensions of $V_j^p$'s 
up to $\dim V_i$.

2. Dimension of the frontier of a semivariety ($Sing_a(V_j^p,V_i)
\subset V_i\cap \overline{V^p_j}$) has dimension strictly smaller 
that a semivariety ($V_j^p$) itself.

3. By lemma \ref{one} the set $Sing_a(V_j,V_i)$ is a semivariety. 
Since a countable union of semivarieties of positive codimension in 
$V_i$ contains $Sing_a(V_j,V_i)$ we have $Sing_a(V_j,V_i)$ has 
a positive codimension in $V_i$ which proves Theorem \ref{sing}.

Since this proves Theorem \ref{sing}, as a consequence this proves 
Theorem \ref{main} too. We are not going to prove these lemmas,
however, we would like to exhibit geometric idea behind the proof
of lemma \ref{two}. The section below is devoted to the idea of
construction of proper subvarieties in the bigger stratum $V_j$
approaching to non-$a$-regular points $Sing_a(V_j,V_i)\subset V_i$.

\subsection{Separation of Planes and dimension reduction
in lemma \ref{two}}

Let $\tau_0$ and $\tau_1$ be two distinct orientable $k$-dimensional 
planes in $\R^N$. An orientable $(m-k)$-dimensional plane $l$ in $\R^N$ 
{\it separates} $\tau_0$ and $\tau_1$ if $l$ is transverse 
to $\tau_0$ and $\tau_1$ and the orientations induced by $\tau_0+l$ 
and $\tau_1+l$ in $\R^N$ are different. Notice that there always 
exists an open set of orientable $(m-k)$-planes separating 
any two distinct orientable $k$-plane.

{\bf Rolle's Lemma.} {\it If a continuous family of orientable 
$k$-planes $\{\tau_t \}_{t \in [0,1]}$ connects $\tau_0$ 
and $\tau_1$ and an orientable $(m-k)$-plane $l$ separates 
$\tau_0$ and $\tau_1$. Then for some $t^*\in (0,1)$ 
transversality of $\tau_{t^*}$ and $l$ fails.}

In what follows we use the transversality theorem \cite{GG} which 
says :\ {\it if $V \subset \R^N$ is a manifold, then almost every
plane of dimension $k$ is transverse to $V$}.  

{\it{An Outline of the Proof of lemma \ref{two}}} \ \   
Let $x \in Sing_{a}(V_j,V_i)$, then by lemma \ref{one}
there are sequences $\{y'_n\},\ \{y_n\}\subset V_j^{con,x}$ 
with different limiting tangent planes $\tau=\lim T_{y_n}V_j$ 
and $\tau'=\lim T_{y_n'}V_j$. Choose an orientation of 
$T_{y_0}V_j$. By connecting $y_0$ locally with 
all other points $\{y'_n\}$ one can induce an orientation on all 
other $T_{y'_n}V_j$ so that the orientations of $\tau_0$ and $\tau_1$
coincide with the orientations of the limits.

Denote $\dim V_j$ by $k$. There is an orientable $(N-k)$-plane $l_j$ 
separating $\tau_0$ and $\tau_1$ and transverse to $V_j$ (by the
transversality Theorem). Consider the orthogonal projection
$\pi_{l_j}$ along $l_j$ onto its orthogonal complement $l_j^\perp$.
Denote by $p_{l_j,j}$ its restriction to $V_j$, 
$p_{l_j,j}=\pi_{l_j}|_{V_j}:V_j \to l_j^\perp$.
Denote by $Crit(l_j,V_j)$ the set of critical points of 
$p_{l_j,j}$ in $V_j$ where the rank of $p_{l_j,j}$ is not maximal.
Then $Crit(l_j,V_j)$ is a semivariety in $V_j$ and 
$\dim Crit(l_j,V_j)<\dim V_j$. Connect two points $y_n^0$ and $y_n^1$ 
by a curve in $V_j$, then $T_{y_n^0}V_j$ deformates continuously  
to $T_{y_n^0}V_j$. Then  by Rolle's Lemma there is a critical point 
of $p_{l_j,j}$ in $V_j^{con,x}$ arbitrarily close to $x$. Thus
$x \in \overline{Crit(l_j,V_j)}$.

By the transversality Theorem there is a countable dense set of 
orientable $(N-k)$-planes $\{l^p_j\}_{p \in \Z_+}$ transverse to 
$V_j$ and separating any two distinct orientable $k$-planes $\tau_0$ 
and $\tau_1$. Therefore, we have that 
\beq \label{unique}
Sing_a(V_j,V_i)\subset V_i \setminus Un(V_j,V_i) 
\subset \bigcup_{p\in \Z_+} 
\left\{\overline{Crit(l^p_j,V_j)}\setminus Crit(l^p_j,V_j)\right\}.
\eneq
These sets $\{Crit(l^p_j,V_j)=V_j^p\}_{p \in \Z_+}$
are proper subsets in $V_j$ we are looking for. 
Using some additional simple argument given in \cite{Ka4}
one can complete the proof of lemma \ref{two}. Q.E.D.

\subsection{A Proof of $a_P$-stratifications for polynomials functions}
\label{hirproof}

The proof below also consists of two steps.

{\it Step 1. Construct a rank stratification of P.}
Consider an $a$-regular stratification $(V,\Cal V^0)$ of $V$
by semialgebraic strata. Such a stratification always exists
by Whitney's Theorem \ref{main} proved above.
Now we refine a stratification $\Cal V^0$ to a stratification
$\Cal V^1$ so that $\Cal V^1$ is a rank stratification of $P$
or restriction of $P$ to any stratum $V_i \in \Cal V^1$ 
is a map of constant rank. Notice that it is sufficient to refine 
each strata $V_i \subset \Cal V^0$ so that $P$ restricted to 
each strata $V_i^j \subset V_i$ has a constant rank.
 
There are two cases: if $P(V_i)$ is a point, then rank of $P|_{V_i}$
is identically zero and $V_i$ stays unchanged and if $P(V_i)$  contains 
an open set, then denote by $\Sigma_{i,P} \subset V_i$ the set of 
critical points of $P|_{V_i}$. By Sard's lemma for algebraic sets 
\cite{Mu} the set $\Sigma_{i,P}$ is a semialgebraic set of positive 
codimension in $V_i$. Refine now each $\Sigma_{i,P}$ to be an $a$-regular 
stratification of $\Sigma_{i,P}$. This is  possible by Whitney's 
Theorem \ref{main}. Denote such a stratification by $(V,\Cal V^1)$. 
By our construction $\Cal V^1$ is an $a$-regular rank stratification of 
$P|_{V}$, i.e. $P$ has constant rank on each stratum and strata ``fit''
$a$-regularly.

{\it Step 2}. It is sufficient to prove the following 
 
\bthm \label{strataher} Let $V_i,V_j \subset \Cal V^1$
be two strata in $\R^N$ and $P:\R^N \to \R$ has a constant 
on each strata $V_i$ and $V_j$. Then the set of singular points
\beq \nonumber
Sing_{a,P}(V_j,V_i)=\{ x \in V_i \cap \overline{V_j}:\ 
V_j \ \textup{is not}\ a_P \textup{-regular over}\ V_i \ 
\textup{at}\ x\ \textup{w.r.t.}\ P\} 
\eneq
is semialgebraic and has positive codimension in $V_i$.
\ethm

Inductive refinement arguments from section \ref{whit} along
with Theorem \ref{strataher} complete the proof of Hironaka's 
Theorem \ref{Hi}. The rest of the section is devoted to 
a proof of Theorem \ref{strataher}.

{\it Proof of Theorem \ref{strataher}:} Similarly to the proof of 
existence of Whitney`s Theorem \ref{main} above we define the set 
with a unique limit of tangent planes to level sets of $P$
\beq
\begin{aligned}
Un_P(V_j,V_i)=\{x \in V_i: \textup{for any}\ V_j^{con,x}\ 
\textup{there is}\ \tau_{x,P}\ \textup{such that}  \\ 
 \lim_{y_n\to x} ker \  dP|_{V_j^{con,x}}(y_n)=\tau_{x,P}\  
\textup{ and is unique}\}.
\end{aligned}
\eneq

\blm \label{aone}
With the above notations $Un_P(V_j,V_i)$ and 
$Sing_{a,P}(V_j,V_i)$ are semivarieties and
\beq
Sing_{a,P}(V_j,V_i) \subset V_i\setminus Un_P(V_j,V_i).
\eneq
\elm

\blm \label{atwo} With notations above there is a set of 
strata $\{V_j^p\}_{p\in \Z}$ (resp. $\{V_i^p\}_{p\in \Z}$) in 
$V_j$ (resp. in $V_i$) each of positive codimension in $V_j$ 
(resp. in $V_i$) such that
\beq
Sing_{a,P}(V_j,V_i) \subset 
\bigcup_{p\in \Z} Sing_{a,P}(V_j^p,V_i)
\bigcup_{p\in \Z} V_i^p. 
\eneq
\elm
\brm 
Similarly to the remarks after lemma \ref{two} this lemma allows 
to reduce dimension of $V_j$'s and prove that 
$Sing_{a,P}(V_j,V_i)$ has positive codimension in $V_i$.
This would prove Theorem \ref{strataher} and as a consequence
it would prove Theorem \ref{Hi}. So what is left to prove is 
lemmas \ref{aone} and \ref{atwo}.
\erm

{\it Proof of lemma \ref{aone}:} \ The proof goes by contradiction. 
Let $x\in Un_P(V_j,V_i)\cap Sing_{a,P}(V_j,V_i)$. Then for any
local connected component $V^{x,con}$ and any sequence
$\{y_n\}\subset V^{x,con}$ there a limiting plane 
$\lim ker \ dP|_{V_j}(x)=\tau_{x,P}$. Moreover, we have 
$ker \ dP|_{V_i}(x) \not\subset \tau_{x,P}$. 
Thus, there is a unit vector $v \in  ker \  dP|_{V_i}$ and
$v \notin \tau_{x,P}$. Contradiction we are going to get is
to find a sequence of points $\{y_n\}\subset V^{x,con}$
such that $\lim ker\ dP|_{V_j}(y_n) =\tau' \supset v$.
The rest of the proof is devoted to construction of such
a sequence.

By Theorem on implicit function one can straighten $V_i$ along 
with nonempty level sets $P\inv(a) \cap V_i$. Then the ray  
$l_v(x)=\{y\in\R^N:\ (y-x)/|y-x|=v\}\subset P\inv(P(x))
\subset V_i$ 
belongs to the level set $P\inv(P(x))$. By an extension
of Wall \cite{Wa} of Milnor's curve selection lemma there is a
$2$-dimensional ``wing'' $V_{j,v}\subset V_j$ such that 
$l_v(x) \subset \overline{V_{j,v}}$. 

By lemma \ref{two} the set of points with nonunique limit 
$Un(V_{j,v},l_v(x))$ is $0$-di\-men\-si\-o\-nal. Therefore, by 
lemma \ref{one} there is a neighborhood $U_x$ of $x$ such that 
any $y \in U_x\cap l_v(x)$, may be distinct from $x$, $V_{j,v}$ 
is $a$-regular over $l_v(x)$ at $y$ and the limit $\tau_y$ of tangent 
planes $T_{y_n}V_{j,v}$ as $y\to x$ is unique. The last two 
properties imply that $\tau_y$ depends continuously on $y$ as long as 
the limit $tau_y$ is unique. Therefore, there is a neighborhood 
$U_y\subset U_x \setminus x$ of $y$ such that $V_{j,v} \cap U_y$ is 
a $C^1$-manifold with a boundary.

Consider a $C^1$-smooth one-sided chart in $U_y\cap V_{j,v}$ and 
$\pi_L$ is the map from $U_y\cap V_{j,v}$ into the $2$-dimensional
plane $\R^2$. The image $\pi_L(V_i)$ is a line 
in $\R^2$ and $\pi_L(V_{j,v})$ is a one-sided neighborhood of this line. 
Using Rolle's type of argument it is easy to show that for
a sequence $\{y_n\} \subset \pi_L(V_{j,v})$ from
a semineighborhood of $\pi_L(x')\in \R^2$ such that 
$y_n \to \pi_L(y)$ and 
$T_{y_n'}\{\pi_L\circ P\inv (P(\pi_L\inv(y_n'))\} \to \pi_L(v)$.
This implies that $\lim ker\ dP|_{V_{j,v}}(y_n')$ tends to $v$, 
however, $V_{j,v}\subset V^{con,x}_j$. So 
$v=\lim_{n \to \infty} ker\ dP|_{V_{j,v}}(\pi_L\inv(y_n))\subset$
\newline
$\lim_{n \to \infty} ker\ dP|_{V_{j}^{con,x}}(\pi_L\inv(y_n))=\tau$. 
This is a contradiction with $v\notin \tau$. Q.E.D.

{\it Proof of lemma \ref{atwo}:} The proof is almost the same as 
the proof of lemma \ref{two} outlined above \cite{Ka4}. Q.E.D.

This completes the proof of Theorem  \ref{strataher} of Hironaka.
Q.E.D.

{\it Acknowledgments:}\ \ I would like to thank my thesis advisor 
John Mather and David Nadler for stimulating discussions and numerous 
remarks on mathematics of this lecture.

 \chapter{Bifurcation of Spatial Polycycles and 
Blow-up along the diagonal of the space of Multijets}

\medskip

\def\asik{Khovanski}
\def\W{\Omega}
\def\land{\wedge}
\def\ell{\mathsf p} 
\def\F{\cal F}
\def\J{\bf J}
\def\j{\bf j}
\def\o{\emptyset}
\def\X{\bf X}

In this lecture we discuss an essential ingredient of the proof of 
Theorem \ref{Space} \cite{Ka2} about an estimate on cyclicity of 
spatial quasielementary polycycles. First, in section \ref{spatial}
we motivate appearance of multichain maps (\ref{jetchain}) to get 
an estimate on cyclicity of spatial polycycles. Similarly to 
the planar case the question of estimating cyclicity of 
a quasielementary polycycle reduces to estimating geometric 
multiplicity of a multichain map of the form (\ref{jetchain}). 
To get an estimate on geometric multiplicity of a multichain 
map (\ref{jetchain}) one needs to prove a Bezout's type Theorem for 
multichain maps. However, a straightforward way to prove Bezout's type 
Theorem for multichain maps faces a typical in singularity theory 
problem, namely, the problem that transversality fails on the diagonal 
in the space of multijets (see e.g. \cite{GG} or \cite{Ma}).  
We shall overcome this problem using a construction of 
Grigoriev-Yakovenko \cite{GY} of blow-up along the diagonal in 
the space of multijets and a special Multijet Transversality 
Theorem \cite{GY}. This construction is described in section
\ref{blown-up} and its relation to Newton Interpolation Polynomials. 
In section \ref{growth} we describe the problem of rate of growth 
of the number of periodic points from Smooth Dynamical Systems 
( see e.g. \cite{AM} and \cite{Sm}) and outline the main result 
of the author along with Brian Hunt \cite{KH} and \cite{Ka9} 
in this direction. Finally, in 
section \ref{Newton-usage} we outline how  Newton Interpolation 
Polynomials can be applied to perturb trajectories and control
the number of periodic points of diffeomorphisms.

\section{Multichain maps and spatial polycycles}\label{spatial}

Consider the simplest example of a polycycle $\gm$ in $\R^3$
consisting of a saddle equilibria $p$ and a connecting
separatrix $\gm_p$ (see Fig.4.1).

 \begin{figure}[htbp]
  \begin{center}
    \begin{psfrags}
      \psfrag{S1}{$\Sigma_+$}
      \psfrag{O}{$p$}
      \psfrag{S2}{$\Sigma_-$}
      \psfrag{F}{$F$}
      \psfrag{D}{$\Delta_p$}
     \includegraphics[width= 4in,angle=0]{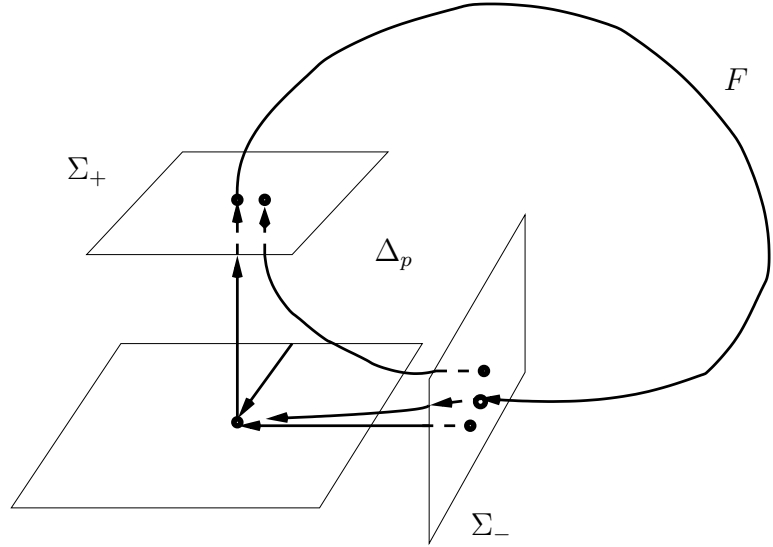}
    \end{psfrags}
    \caption{A saddle loop polycycle}
  \end{center}
\end{figure}

Let $\Sigma^-$ and $\Sigma^+$ be ``entrance'' and ``exit'' 
transverse sections to $\gm_p$ chosen so that in $C^r$-normal
coordinates the Poincare return map $\Delta_p$ along the 
polycycle $\gm$ has a ``nice'' form. We decompose the Poincare 
return map $\Delta$ along the polycycle $\gm$ into the composition 
of a local Poincare map $\Delta_p$ in a neighborhood of $p$ and 
a semilocal map $f$ along $\gm_p$. Consider $2$-cycles bifurcating 
from $\gm$. Denote by $x_1, x_2$ and $y_1, y_2$ the first and 
the second intersection of each of $2$-cycles with $\Sigma^+$ and 
$\Sigma^-$ respectively. Then the equation determining the number 
of $2$-cycles has the form
\beq \label{2-cycle}
\begin{cases}
y_1=f(x_1,\eps) \\
x_2=\Delta(y_1,\eps)\\
y_2=f(x_2,\eps) \\
x_1=\Delta(y_2,\eps).
\end{cases}
\eneq

Important that the first and the third equation have
the same functional parts. Notice also that each of the equations
from (\ref{2-cycle}) is an equality in $\R^2$ so it consists of two 
$1$-dimensional equalities itself. So the total number of equations 
in (\ref{2-cycle}) is 8. Following the strategy of the planar 
case from section \ref{asik} lecture 2 we apply the Khovanskii
method to the system (\ref{2-cycle}). Compare this system
with the system (\ref{prebasic}) or (\ref{funpfaf}). It is not 
difficult to see that the result of application of the Khovanskii
method give the map of the form 
\beq
P \circ (j^7 f,j^7 f):\R^8 \to \R^8.
\eneq
To simplify our considerations we redenote the 
seventh jet $j^7 f$ of $f$ by a map $F$ and consider 
the multichain map 
\beq
P \circ (F,F):\R^8 \to \R^8,
\eneq
where $F$ is a {\it{generic map}} in a sense that it satisfies 
{\it{any ahead given transversality condition}}. It is clear 
that even if $F$ is generic we can not assume that 
the $2$-tuple $(F,F)$ is a generic map. Simply because the first 
and the second components are the same. 
Let's explain why genericity fails for a $2$-tuple mapping by
an example.

\subsection{Genericity (resp. Transversality) fails for $2$-tuple
mappings!}

{\bf The Classical Transversality Theorem} {\it (e.g. \cite{AGV}, 
\cite{GG}) Let $N$ and $m$ be positive integer and $M$ be a smooth 
compact manifold in $\R^m$. Then for an open dense set of smooth 
mappings $F:\R^N \to \R^m$ we have that  $F$ is transverse to $M$.
In particular, it means that $F\inv(M)$ is a smooth manifold.}

{\bf Remarks} \ 
1. It is an exercise from calculus to construct
a set on the unit interval $[0,1]$ which is open dense and
has an arbitrary small positive measure. To justify that 
transversality property is, indeed, generic there is so-called 
prevalent extension of the Classical Transversality Theorem 
which says that for a.e. mapping $F:B^n \to \R^N$ we have that 
$F$ is transverse to $M$. More exactly, for a generic 
finite-parameter family of mappings 
$\{F_\eps:B^n \to \R^N\}_{\eps \in B^k}$ for a.e.
parameter value $F_\eps$ is transverse to $M$.
See \cite{HSY} and \cite{Ka7} for more.
\newline   
2. The fact that transversality of $F$ to $M$ implies
that $F\inv(M)$ is a smooth manifold follows from
the Theorem on implicit function (see {\cite{AGV}, \cite{GG}}). 

Since a genericity condition on $F$ we need is that $F$ 
has to satisfy a transversality condition, to show impossibility
of application of the classical transversality Theorem we give 
a trivial example when a transversality fails for an open set 
of $2$-tuples $(F,F)$.

\bex In the Classical Transversality Theorem put $n=m=1$.
Consider the function $f:x \to x^2$ for $x\in I=[-1,1]$ and 
the corresponding $2$-tuple $f\times f:I \times I\to 
\R \times \R$, given by $f\times f:(x_1,x_2) \to (x_1^2, x_2^2)
=(y_1,y_2)$. Let $M=\{y_1=y_2\} \subset \R \times \R$ be the
diagonal.  Then for each $\tilde f$ which is $C^1$-close to $f$
the preimage $(\tilde f,\tilde f)\inv(M)$ is a topological
cross (not a manifold). This, in particular, 
implies that $(\tilde f,\tilde f)$ is not transverse to $M$,
otherwise the preimage of a manifold should be a manifold.
\enex

To see that notice that a function  $\tilde f$ close to $f$
has to have a local minima $\tilde x$ close to $0$ and $\tilde x$ 
is a nondegenerate local minimum, i.e.
$\tilde f:x \to \~\eps+ \~a (x-\tilde x)^2 + 
HOT\left((x-\tilde x)^2\right)$ with $\~a \neq 0$. Then 
$\tilde f(x_1)-\tilde f(x_2)=0$ has two intersecting curve of 
solutions $x_1=x_2$ and $x_1-\tilde x \approx -(x_2 - \tilde x)$ 
which form a cross.
This completes the proof of the claim in the example.
 
To explain what happens in this example and we derive 
a general frame due to Grigoriev-Yakovenko \cite{GY}. 
 
\subsection{Blow-up along the diagonal for $2$-tuples
in the $1$-di\-men\-si\-o\-nal case}

For a smooth function $\tilde f:\R \to \R$ consider the map 
\beq \label{ddmap}
\begin{aligned}
\ &(x_1,x_2) \stackrel{\cal{D}_2 \tilde f}{\longrightarrow} 
 \left(x_1,x_2,\tilde f(x_1), 
\frac{\tilde f(x_2)-\tilde f(x_1)}{x_2-x_1} \right) 
= (x_1,x_1,u_1,u_2) \subset \R^4,\\
\ & (x_1,x_2,u_1,u_2) \stackrel{\pi_2}{\longrightarrow} 
(x_1,x_2,u_1,u_1+u_2 (x_2-x_1)).
\end{aligned}
\eneq 
Direct calculation shows that  
$\pi_2 \circ \cal{D}_2 \tilde f \equiv (\tilde f,\tilde f)$.
Therefore,  
\beq
(\tilde f,\tilde f)\inv(M)=
(\cal{D}_2 \tilde f)\inv \circ \pi_2\inv (M).
\eneq
This is incorporated into diagram 1.3 of lecture 1 with
$n=N=1$. By definition ${\cal D}_2(f): I \times I \to 
{\cal{DD}}_2(I,\R)$ is a smooth map, 
$\pi_2:{\cal{DD}}_2(I,\R)\to\R$ is an explicitly computable 
polynomial map, and  
$\pi_2 \circ {\cal D}_2(f)=(f,f):I \times I \to \R^2$.
Notice that outside of the diagonal $\{x_1=x_2\}$ the map
$\pi_2$ is one-to-one. However, the preimage of the set
$\pi\inv_2\{x_1=x_2,\ f(x_1)=f(x_2)\}$ is of dimension $3$
while the set $\{x_1=x_2,\ f(x_1)=f(x_2)\}$ itself
is of dimension $2$. This, in particular, means that 
$\pi_2$ is a blow-up along the diagonal.

Consider $\pi_2\inv (M)=\{u_2(x_2-x_1)=0\}\subset \R^4$.
This is the union of two intersecting hyperplanes. If the map 
$\cal{D}_2 \tilde f$ is transverse to $\{u_2(x_2-x_1)=0\},$ 
then the preimage $\left(\cal{D}_2 \tilde f\right)\inv 
\left( \{u_2(x_2-x_1)=0\}\right)$
has to be the union of two intersecting curves.
It turns out that assumption that $\cal{D}_2 \tilde f$
is {\it{generic for a generic}} $\tilde f$ is satisfied
or for a generic $\tilde F$ the map $\cal{D}_2 \tilde f$
is transverse to both hyperplanes of $\{u_2(x_2-x_1)=0\}$.
Let's justify that.

{\it A ``Proof'' of the Classical Transversality Theorem}.
(see e.g. \cite{AGV}, \cite{GG})
Transversality is an open property, i.e. if $F$ is transverse 
to $M$, then for all $\tilde F$ sufficiently close to $F$
we have $\tilde F$ is transverse to $M$ too.
So it is sufficient to show that by an arbitrary small 
perturbation of any mapping $F:\R^N \to \R^m$ one can reach 
transversality to a compact manifold $M$. Let's prove it now.

Consider a smooth mapping $F:B^n \to \R^N$. Include this mapping 
into the $m$-parameter family
${\bf F}:\R^n \times \R^N \to \R^n \times \R^N$, given by 
${\bf F}(x,\eps)=(x,F(x)+\eps )$. The determinant of the
linearization of the mapping (the Jacobian) $J_{\bf F}(x,\eps)$
is constant and identically equals $1$. Therefore,
${\bf F}$ is a diffeomorphism and $M_F={\bf F}\inv(\R^n \times M)$
is a manifold in the preimage $\R^n \times \R^N$.

{\bf Fact.}\ {\it If $\eps$ is a regular point of the projection
$\pi_{M,F}=\pi|_{M_F}:\R^n \times \R^N \to \R^N$ along the 
$x$-coordinate, restricted to $M_F$, then 
$F_\eps=F(x)+\eps $ is transverse to $M$. }

This follows from the implicit function theorem.

{\bf Sard's Lemma.} (e.g. \cite{Mi2}) {\it A.e. $\eps \in \R^N$ 
value is regular for the projection map $\pi_{M,F}$.} 

Thus, one can choose a regular value $\eps$ arbitrary close to
$0$. For such an $\eps$ the mapping $F_\eps$ is transverse to
$M$. This completes the proof of the Classical Transversality
Theorem. Q.E.D.

Now we are ready to state the main result of this section

\subsection{Multijet Transversality Theorem due to 
Grigoriev-Ya\-ko\-ven\-ko}

\bthm \cite{GY} Let $M \subset \R^N \times \R^N$ be an algebraic 
manifold (or variety) and $B^n\subset \R^n$ be a unit ball. Then 
for an open dense set of smooth mappings $F:B^n \to \R^N$ the set 
$\left(F \times F \right)\inv (M)$ is stratified.

Moreover, let $k \in \Z_+$ and 
$M \subset \R^N \times  \dots \times \R^N$ ($k$ times)
be an algebraic manifold (or variety). Then for an open 
dense set of smooth mappings $F:B^n \to \R^N$ the set 
$\left(F \times \dots \times F \right)\inv (M)$, with $k$ 
repetitions,  is stratified.

Moreover, let $n,\ k \in \Z_+$ and 
$M \subset J^m(B^n,\R^N) \times  \dots \times J^m(B^n,\R^N)$ 
($k$ times) be an algebraic manifold (or variety). Then for 
an open dense set of smooth mappings $F:B^n \to \R^N$ the set 
$\left(F \times \dots \times F \right)\inv (M)$, with $k$ 
repetitions,  is stratified.
\ethm 

{\it{A Proof of the Theorem for the model example $n=m=1$ and
$k=2$.}}\ \  Consider the map 
${\cal D}_2(F):(x_1,x_2,\eps_1,\eps_2)\mapsto (x_1,x_2,u_1,u_2),$
defined by the formula (\ref{ddmap}). Direct calculations show
that the determinant of the linearization (the Jacobian)
$J_{{\cal D}_2(F)}(x_1,x_2,\eps_1,\eps_2)\equiv 1$
and is formed by an upper triangular matrix with units on
the diagonal. Since, ${\cal D}_2(F)$ is a diffeomorphism
one can apply arguments of the proof of the Classical
Transversality Theorem given above. Q.E.D.

A Proof of the Theorem in the general case follows along 
the same lines. The main difficulty is to construct diagram 1.3
in the general case. This is the subject of the next subsection.

\bcor For an open dense set of smooth functions $F:I \to \R$ 
the preimage $(F,F)\inv(M)$ is a $1$-dimensional stratified \
manifold, i.e. locally finite union of points and curves.
\ecor

\section{Newton Interpolation Polynomials and Blow-up Along
the diagonal in the space of Multijets}\label{blown-up}

This section is devoted to description of Grigoriev-Yakovenko
construction of Blow-up along the diagonal in the space of 
Multijets in the general case. 
Let $F:B^n\to \R^N$ be a smooth map of a unit ball $B^n\subset \R^n$, 
$j^mF:B^n \to J^m(B^n,\R^N)$ be an $m$-th jet of $F$, and
${\Cal J}^{m,k}(B^n,\R^N)=J^m(B^n,\R^N)\times\dots \times J^m(B^n,\R^N)$ 
($k$ repetitions) be the space of $k$-tuple $m$-jets. 
Denote $k$-tuple of $m$-jet of a map $F:\R^N\to \R^m$ by
${\Cal J}^{m,k} F(x_1,\dots,x_k)=(j^mF(x_1), \dots, j^mF(x_k))$.
The goal of this section is to define entries of an extension
of diagram 1.3:
The, so-called, space of divided differences ${\cal{DD}}_k^m(B^n,\R^N)$,
the Newton map 
$\pi^m_k:{\cal{DD}}_k^m(B^n,\R^N)\to {\Cal J}^{m,k}(B^n,\R^N)$,
${\Cal D}_k^m(F):\underbrace{B^n \times \dots 
\times B^n}_{k\ \textup{repetitions}} \to {\cal{DD}}_k^m(B^n,\R^N)$.
We use the exposition from {\cite{GY}}.

\begin{figure}[htbp]
  \begin{center}
    \begin{psfrags}
      \psfrag{Phase}{ $\underbrace{B^n \times \dots 
      \times B^n}_{k\ \textup{repetitions}}$}
      \psfrag{MultJ}{${\Cal J}^{m,k}(B^n,\R^N)$}
      \psfrag{DivD}{${\cal{DD}}_k^m(B^n,\R^N)$}
      \psfrag{pi}{$\pi_k^m$}
      \psfrag{N}{${\Cal D}_k^m(F)$}
      \psfrag{(jf,...,jf)}{${\Cal J}^{k,m}F$}
      \includegraphics[width= 4in,angle=0]{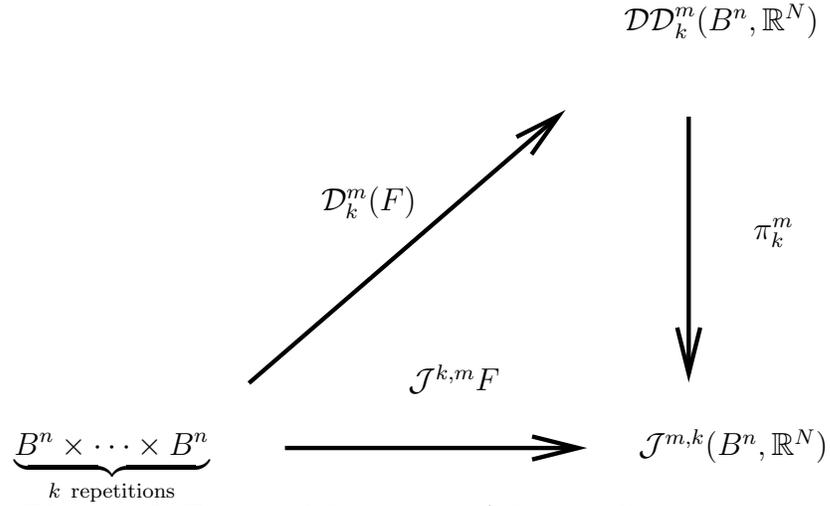}
    \end{psfrags}
    \caption{Polynomial blow-up of the multijet space}
  \end{center}
\end{figure}

\subsection{Divided Differences}\label{divided}

In order to extend the above construction we need to 
define so called {\it divided differences}.
Let $g:\R \to \R$ be a sufficiently smooth function of one real
variables. 
\bdef \label{divdif}
The {\it first order divided difference} of $g$ is  defined as
\beq \begin{aligned}
\Delta g(x_1, x_2)=\frac {g(x_2)-g(x_1)}{x_2-x_1}
\end{aligned}
\eneq
for $x_2 \neq x_1$ and extended by its limit value as
$g'(x)$ for $x=x_2=x_1$. Clearly, if $g$ is a $C^r$-smooth function, 
then  $\Delta g$ is at least a $C^{r-1}$-smooth
function of its arguments.

Iterating this construction we define {\it divided differences 
of the $s$-th order for any $s\in \Z_+$} as
\beq \nonumber 
\Delta^s g(x_1,\dots , x_{s+1})=
\frac {\Delta^{s-1} g(x_1,\dots, x_{s-1}, x_{s+1})-
\Delta^{s-1} g(x_1,\dots, x_{s-1}, x_s)}{x_{s+1}-x_s}
\eneq
for $x_{s+1} \neq x_s$ and extended by its limit value as
$\frac{\partial \Delta^{s-1} g}{\partial x_s}(x)$ for 
$x=x_{s+1}=x_s$. Clearly, if $g$ is a $C^r$-smooth function, 
then  $\Delta g$ is at least a $C^{r-s}$-smooth
function of its arguments.
\endef  
Notice that $\Delta^s$ is linear as a function of $g$, and
one can show that it is a symmetric function of $x_1,\ldots,x_{s+1}$; 
in fact, by induction it follows that
\beq
\Delta^s g(x_1, \dots, x_s)=
\sum_{i=1}^{s+1} \frac{g(x_i)}{\prod_{j \neq i} (x_i - x_j)}
\eneq
Another identity that is proved by induction will be more important
for us, namely
\beq
\Delta^s\ x^l(x_1,\dots, x_{s+1})=p_{l,s}(x_1,\dots, x_{s+1}),
\eneq 
where $p_{l,s}(x_1, \dots, x_{s+1})$ is $0$ for $s> l$ and for 
$s\leq l$ is the sum of all degree $l-s$ monomials in $x_0,\dots,x_s$ 
with unit coefficients,
\beq \label{homogpolyn}
p_{l,s}(x_1,\dots, x_{s+1})=\sum_{r_0+\dots + r_s=l-s} 
\quad \prod_{j=1}^{s+1} x_j^{r_j}.
\eneq 

The divided differences form coefficients for the Newton 
interpolation formula.  For all $C^\infty$ functions 
$g : \R \to \R$ we have
\beq 
\beal\label{int}
g(x) = & \Delta^0 g(x_1) +\Delta^1 g(x_1,x_2) (x-x_1)+ \dots \\
& + \Delta^{k-1} g(x_1,\dots, x_{k-2}) (x-x_1) \dots (x-x_{k-3}) \\
& + \Delta^{k} g(x_1,\dots, x_{k-1},x) (x-x_1) \dots (x-x_{k-2})
\enal
\eneq
identically for all values of $x, x_1, \dots, x_{k}$. 
All terms of this representation are polynomial in $x$
except for the last one which we view as a remainder term.
The sum of the polynomial terms is the degree $(k-1)$ {\em Newton
interpolation polynomial\/} for $g$ at $\{x_s\}_{s=1}^k$.  To
obtain a degree $2k-1$ interpolation polynomial for $g$ and its
derivative at $\{x_s\}_{s=1}^k$, we simply use (\ref{int}) with
$k$ replaced by $2k$ and the $2k$-tuple of points
$\{x_{s(\textup{mod}\ k)}\}_{s=1}^{2k}$. Similarly one can 
construct an interpolation polynomial for $g$ and its
derivatives up to any finite order. 
 
All terms of this representation, except for the last one, 
are polynomial in $x$ and their sum is the $k$-th order 
{\it Newton Interpolation Polynomial} denoted by 
$\mathcal P_{k-1}(x,{\X}_k)$, where ${\X}_k=(x_1,\dots ,x_k)$.

Now we can define entries of diagram 4.2 in the case $m=0$
Let ${\cal{DD}}_k(I,\R)=\underbrace{I \times \dots 
\times I}_{k\ \textup{times}} \times \R^k=
(x_1,\dots,x_k;u_0,u_1,\dots,u_{k-1})$. It is called 
{\it the space of divided differences}. Then
\beq
\begin{aligned}
{\cal D}_k(f): \underbrace{I \times \dots 
\times I}_{k\ \textup{times}} & \to {\cal{DD}}_k(I,\R),\\
{\cal D}_k(f):(x_1,\dots,x_k) & \mapsto 
(x_1,\dots,x_k;u_0,\dots,u_{k-1}),
\quad u_\al=\Delta^\al f(x_1,\dots,x_{\al+1}) \\
\pi_k:  {\cal{DD}}_k(I,\R) & \to \R^{2k},
\end{aligned}
\eneq
where ${\cal D}_k(f)$ is smooth, provided that $F$ is smooth, 
$\pi_n:{\cal{DD}}_k(I,\R) \to \R^{2k}$ is a Newton Interpolation 
Polynomial of the form (\ref{int}), and 
\beq
\pi_k \circ {\cal D}_k(f)=(f,\dots,f):
I \times \dots \times I \to \R^{2k},
\eneq
where $f$ and $\R$ are repeated $k$ times.

\subsection{ Language of divided differences and
the Newton interpolation formula}  

In this section we introduce construction of divided differences 
space ${\cal {DD}}_k(B^n,\R^N)$ and the corresponding map $\cal D_k(F)$ 
and the polynomial $\pi_k$ presented on diagram 1.3. 

Let $F:{\Bbb R}^n \to {\Bbb R}$ be a smooth function in $n$ 
real variables $x_1,\dots, x_n$. 

\bdef \label{divdif1}
The {\it first order divided difference} of $F$ in the variable
$x_k$ is the function of $n+1$ variables $x_1,\dots, x_{k-1},
x'_k,x''_k,\dots, x_n$  defined as
\beq \begin{aligned}
\Delta_{x_k} F(x_1,\dots, x_{k-1},x'_k,x''_k,\dots, x_n)=
\\ \frac {F(x_1,\dots, x_{k-1},x'_k,\dots, x_n)-
F(x_1,\dots, x_{k-1},x''_k,\dots, x_n)}{x'_k-x''_k}
\end{aligned}
\eneq
for $x'_k \neq x''_k$ and extended by its limit value as
$\frac {\partial F}{\partial x_k}(x_1,\dots, x_{k-1},x'_k,\dots, x_n)$
for $x'_k=x''_k=x_k$. Clearly, if $F$ is $C^r$ function, then (e.g., 
by the Hadamard lemma), $\Delta_{x_k}F$ is at least $C^{r-1}$-smooth
function of its arguments.
\endef 
It turns out that iterating this construction is possible {\cite{GY}}
which leads to 

\bdef
Let $\alpha=(\alpha_1, \dots, \alpha_n) \in {\Bbb Z}^n_+$ be 
a multiindex, let $F$ be as above. Then
$\Delta_x^\alpha F = \Delta_{x_1}^{\alpha_1}\dots 
\Delta_{x_n}^{\alpha_n} F$ is called the mixed divided difference 
of order $|\alpha|=\alpha_1+ \dots + \alpha_n$. This divided
difference is a smooth function of $n + |\alpha|$ arguments subdivided
into $n$ groups of $\alpha_1+1, \dots ,\alpha_n + 1$ variables, 
symmetric with respect to permutations of variables within 
the same groups.
\endef 
As direct calculations show the operators $\Delta_{x_j}$ and 
$\Delta_{x_i}$ commute for $i\neq j$, and, therefore, we can use 
the multiindex notation for divided differences.

\subsection{The Newton interpolation formula (in multivariables)}
\label{Newtmulti}
                    
Let\ \ $X^1=(x^1_1,\dots,x^1_N) \subset {\Bbb R}, \dots, 
X^n=(x^n_1,\dots,x^n_N) 
\subset {\Bbb R}$ be a subsets consisting of the same number of 
points, each $X^j$ belonging to the corresponding $j$-th coordinate 
axis of points in ${\Bbb R}^N.$ Then, given a multiindex 
$\alpha \in \Bbb Z^N_+$ and a smooth function $F(x)=F(x^1,\dots,x^N)$ 
in $N$ variables we can form the divided difference 
$\Delta^\alpha_{x}F(X^1,\dots,X^N)$.

In terms of the divided differences one can write the  
Newton interpolation polynomial as follows:

\beq \beal \label{Newt}
\mathcal P(t^1,\dots,t^N)= & 
\sum_{0\leq  \alpha_i \leq n} 
\Delta^\alpha_{x}F(X^1,\dots,X^n) \prod_{i_1=1}^{\alpha_1}
\dots \prod_{i_n=1}^{\alpha_n}(t^1 - x^1_{i_1}) \dots 
(t^n - x^n_{i_n}).
\enal
\eneq

The polynomial $\mathcal P(t^1,\dots,t^n)$ has degree $\leq kn$
in variables $t=(t^1,\dots,t^n)$. The Newton interpolation formula 
implies that the difference $F(t)-\mathcal P(t^1,\dots,t^n)$
vanishes at all points of the Cartesian product
grid ${\bf X}=X^1 \times \cdots \times X^n \subset {\Bbb R}^n$.
Moreover, if for each $X^j=(x^j_1,\dots,x^j_n)$ we denote by 
$\text{diag}^k(X^j)$ the set $(x^j_1,\dots,x^j_n)$ repeated 
$(m+1)$ times
\beq
(\underbrace{x^j_1,\dots,x^j_n},\dots ,\underbrace{x^j_1,\dots,x^j_n})
\quad (m\ \textup{times}),
\eneq
then to obtain interpolation of the $m$-th jet of $F$ we replace 
each $X^j=(x^j_1,\dots,x^j_n)$ by $\text{diag}^m(X^j)$.
The degree of interpolating polynomial will be $\leq nk(m+1)$.

In the case of a multivariate function $F: B^n \to {\Bbb R}^N$
interpolating polynomial  $\mathcal P(t^1,\dots,t^n)$ becomes 
$N$-dimensional vector and is interpolating by coordinate 
functions of $F=(F^1,\dots ,F^N)$. 

\bdef \label{multdiv}
Let ${\cal {DD}}^m_k(B^n,\R^N)$ be the collection of all divided 
differences with $m$ repetitions, 
$\{\Delta^\alpha_xF (\text{diag}^m(X^1),\dots, 
\text{diag}^m(X^n))\}_\alpha,\ \alpha_i \leq (m+1)k,\ i=1,\dots,n$.
This is a linear space naturally equipped with the coordinates    
$\{ x_i, u_\alpha:\ 0 \leq i \leq N,\ \alpha_i \leq (m+1)\}$,
where $x_i$ (resp. $u_\alpha$) are vectors from ${\Bbb R}^n$ 
(resp., ${\Bbb R}^N$).  Dimension of this space is equal to 
$kn + N((m+1)k)^n$.
 
The map ${\cal D}^m_k F$ is defined by 
\beq
\beal
{\cal D}^m_k F: (x_1,\dots,x_k) \ \to \
(x_1,\dots,x_k, \{u_\alpha\}_{\alpha}), \\   
\textup{where}\ \ u_\alpha=\Delta^\alpha_xF,\ \ 
\forall i\ \ \alpha_i\leq (m+1)n.
\enal
\eneq
\endef

The multivariate interpolation 
formula together with its derivatives in $t_j$ evaluated at the 
points of the grid, can be interpreted as a polynomial map  
restoring multijets from divided differences.

{\bf Newton Interpolation on ${\Bbb R}^m$ (abstract version)} {\it  
The multivariate Newton interpolation formula (\ref{Newt}) 
defines a polynomial interpolation map 
$\pi^m_k: {\cal {DD}}^m_k(B^n,\R^N)$\ $\to 
{\mathcal J}^{m,k}(B^n,\R^N)$ such that 
$\mathcal J^{m,k} f=\pi^m_k \circ \mathcal D^m_k F$. 
Degrees of the components of $\pi^m_k$ do not exceed $(k+1)nN.$}

In the next section we present an application of Newton
Interpolation Polynomials and diagram 4.2 to an old problem
in dynamical systems: the problem of rate of growth of the number
of periodic points for generic diffeomorphisms (see e.g.
\cite{AM} and \cite{Sm}).

\section{Rate of growth of the number of periodic points for generic
diffeomorphisms and Newton Interpolation Polynomials}
\label{growth} 

\subsection{Statement of the problem}

Let Diff$^r(M)$ be the space of $C^r$ diffeomorphisms  of 
a finite-dimensional smooth compact manifold $M$ with the uniform 
$C^r$-topology, $\dim M \geq 2,$ and let $f \in {\textup{Diff}}^r(M)$.
Consider the number of {\it{isolated}} periodic points of period $n$

\beq \label{grow}
P_n(f)=\# \{\textup{isolated}\ \ x \in M:\ \  x=f^n(x)\}.
\eneq 
The main question of this paper is:
\beq \nonumber
\boxed{\textup{How quickly can}\  P_n(f)\  {\textup{grow with}
\ n \ \textup{for a ``generic'' diffeomorphism}\ f?}}
\eneq    
We put the word ``generic'' in brackets because as the reader
will see the answer depends on notion of genericity.

We call a diffeomorphism $f \in \textup{Diff}^r(M)$ an 
{\it Artin-Mazur diffeomorphism} (or simply {\it A-M diffeomorphism)} 
if the number of isolated periodic orbits of $f$ grows  at most
exponentially fast, i.e.
for some number $C>0$ 
\beq
P_n(f) \leq \exp(Cn) \ \  {\textup{for all}}\ \  n \in \Bbb Z_+.
\eneq
Artin \& Mazur  {\cite {AM}} proved the following
\bthm For any $0\leq r\leq \infty$, A-M diffeomorphisms form a dense 
set of diffeomorphisms in $\textup{Diff}^r(M)$ with the uniform 
$C^r$-topology. 
\ethm 
In {\cite{Ka5}} an elementary proof of the following extension 
of the Artin-Mazur result is given 
\bthm For any $0\leq r<\infty$ A-M diffeomorphisms with all periodic 
points hyperbolic are dense in $\textup{Diff}^r(M)$ with the uniform 
$C^r$-topology.
\ethm

According to the standard terminology a set in Diff$^r(M)$ is 
called generic ( or residual) if it contains a countable 
intersection of open dense sets and a property is called (Baire) 
generic if diffeomorphisms with that property form a residual set. 
It turns out the A-M property is not generic, as it is shown
in {\cite{Ka6}}. Moreover:

\bthm {\cite{Ka6}} \label{supergrowth} For any $2\leq r <\infty$
there is an open set $\Cal N \subset$ Diff$^r(M)$ such that for 
any given sequence $a=\{a_n\}_{n \in \Bbb Z_+}$ there is a Baire
generic set $\cal R_a$ in $\Cal N$ depending on the sequence $a_n$ 
with the property if $f \in \cal R_a$, then for infinitely many 
$n_i \in \Bbb Z_+$ we have $P_{n_i}(f)>a_{n_i}$.
\ethm

Since any two residual sets have nonempty intersection Theorem
\ref{supergrowth} implies that A-M diffeomorphisms are not 
generic. The proof of this Theorem is based on a result of 
Gonchenko-Shilnikov-Turaev \cite{GST1}. Two slightly
different detailed proofs of their result are given
in {\cite{Ka6}} and {\cite{GST2}}. The proof in {\cite{Ka6}}
relies on a strategy outlined in \cite{GST1}.

However, it seems unnatural that if you pick a diffeomorphism
at random then it may have an arbitrarily fast growth of 
number of periodic points. Moreover, Baire generic sets in Euclidean 
spaces can have zero Lebesgue measure. Phenomena which 
are Baire generic, but have a small probability are well-known
in dynamical systems, KAM theory, number theory, and etc. 
(see {\cite{O}}, {\cite{HSY}}, and {\cite{Ka7} for various examples). 

This partially motivates the problem posed by Arnold {\cite{A2}}: 
{\it{Prove that ``with probability one'' $f$ is an A-M diffeomorphism.}}
Arnold suggested  the following interpretation of ``with probability 
one'': {\it{for a (Baire) generic finite parameter family of 
diffeomorphisms $\{f_\eps\}$, for Lebesgue almost every $\eps$ we have 
that $f_\eps$ is A-M.}} (cf. \cite{Ka7}). 
As Theorem \ref{supergrowth} shows, a result on the genericity of 
the set of A-M diffeomorphisms based on (Baire) topology is likely 
to be extremely subtle, if possible at all\footnote{ For example, using 
techniques from \cite{GST2} and \cite{Ka6} one can prove that for 
a Baire generic finite-parameter family $\{f_\eps\}$ and a Baire generic 
parameter value $\eps$ the corresponding diffeomorphism $f_\eps$ is 
not A-M. Unfortunately, how to estimate from below the measure of 
non-A-M diffeomorphisms in a Baire generic finite-parameter family is
so far an unreachable question.}.  We use instead
a notion of ``probability one'' based on prevalence \cite{HSY,Ka7},
which is independent of Baire genericity.  We also are able to state
the result in the form Arnold suggested for generic families using this
measure-theoretic notion of genericity.
The main result in this direction is a partial solution to 
Arnold's problem. It says that {\it{ For a prevalent 
diffeomorphism $f\in$ Diff$^r(M),\ r>1,$ and any $\dt>0$ 
there exists $C=C(\dt)>0$ such that for all $n \in \Z_+$
\beq \label{growthbound}
P_n(f) \leq \exp(Cn^{1+\dt})
\eneq
}}
This Theorem is announced in \cite{KH}. A major part of the proof 
is worked out in \cite{Ka9}. We omit the precise statement
which requires an additional discussion.

\section{Dynamical Usage of Newton Interpolation Po\-lynomials}
\label{Newton-usage}

\subsection{ Perturbation of recurrent trajectories by
Newton Interpolation Polynomials}

Let us start with several remarks which were the starting 
point of this paper. In order to keep notations and formulas 
simple we consider the 1-dimensional maps, but the reader should 
always have in mind that our consideration is designed for 
multidimensional diffeomorphisms. 

Consider a map $f:I \hookrightarrow I$ of the
interval $I=[-1,1]$. Recall that a trajectory 
$\{x_k\}_{k \in \Z}$ of $f$ is called {\it recurrent} if 
it returns arbitrarily close to its initial position ---
that is, for all $\dt>0$ we have $|x_0-x_n|<\dt$ for some $n>0$. 
A very basic question
of Closing lemma type is how much one should perturb $f$ to 
create a periodic point $x_0$. Let us give a ``baby'' answer

{\bf Baby Closing lemma.}\ {\it Let $\{x_k=f^k(x_0)\}_{k=0}^{n}$ 
be a trajectory of length $n+1$ of a map $f:I \hookrightarrow I$. 
Let $u=(x_n-x_0)/\prod_{k=0}^{n-2}(x_{n-1}-x_k)$.
Then $x_0$ is a periodic point of period $n$ of the map 
\beq \label{closing}
\~f(x)=f(x)+u\prod_{k=0}^{n-2}(x-x_k) 
\eneq}

Of course $\~f$ is close to $f$ only if $u$ is 
sufficiently small, meaning that  $|x_0-x_n|$ is small
compared to $\prod_{k=0}^{n-2}(x_{n-1}-x_k)$. However,
this product is likely to contain small factors for a 
recurrent trajectories. In general, it is difficult to control 
the effect of perturbations for recurrent trajectories. The simple 
reason why is because {\it one can not perturb $f$ at two 
nearby points independently}.

It is important for the proof in \cite{Ka9} to control
on derivative of $f$ along periodic orbits. If for some $x\in I$ 
$\gm>0$ and some positive integer $n$ we have $f^n(x)=x$ and 
$|(f^n)'(x) -1|>\gm$, then it implies that the interval around $x$ 
of size $\|f\|_{C^1}^{-n}\gm$ is free from periodic points
of the same period (see Proposition 1.1 \cite{KH}). Quantity $\gm$
is called {\it hyperbolicity} and $x$ is called $(n,\gm)$-hyperbolic.
This quantity was introduced by Gromov \cite{Go} and Yomdin \cite{Y}. 
If one can estimate hyperbolicity for all points of period $n$ from
below, then one can estimate the number of periodic points of period
$n$. Upper bound (\ref{growthbound}) is obtained by proving
lower bound on the rate of decay of hyperbolicity with period
for prevalent diffeomorphisms. This is the reason the proof needs
to control derivative along trajectories.
  
The Closing Lemma above also gives an idea of how much we must change
the parameter $u$ to make a point $x_0$ that is $(n,\gm)$-periodic not
be $(n,\gm)$-periodic for a given $\gm > 0$, which as we described
above is one way to make a map that is ``bad'' for the initial
condition $x_0$ become ``good''.  To make use of our other alternative
we must determine how much we need to perturb a map $f$ to make a
given $x_0$ be $(n,\gm)$-hyperbolic for some $\gm > 0$.

\smallskip
\noindent
{\bf Perturbation of hyperbolicity.}\ {\it 
Let $\{x_k=f^k(x_0)\}_{k=0}^{n-1}$ be a trajectory 
of length $n$ of a $C^1$ map $f:I \to I$.  Then for the map 
\beq \label{movehyperbolicity}
f_v(x)=f(x)+v(x-x_{n-1})\prod_{k=0}^{n-2}(x-x_k)^2
\eneq
such that $v\in \R$ and 
\beq
\left|\vphantom{{f'}^2}|(f^n_v)'(x_0)| - 1\right| =
\left|\vphantom{{{\prod_0^n}^2}^2}
\left|\prod_{k=0}^{n-1}f'(x_k)+
v\prod_{k=0}^{n-2}(x_{n-1}-x_k)^2\prod_{k=0}^{n-2}f'(x_k)\right|
- 1\right| > \gm
\eneq
we have that $x_0$ is an $(n,\gm)$-hyperbolic point of $f_v$.}
\noindent

One more time we can see the product of distances
$\prod_{k=0}^{n-2}|x_{n-1}-x_k|$ along the trajectory
is important quantitative characteristic of how much freedom
we have to perturb.

The perturbations (\ref{closing}) and (\ref{movehyperbolicity}) are
reminiscent of Newton interpolation polynomials.  Let us put these
formulas into a general setting using singularity theory.

\subsection{Distance to the diagonal in the multijet space}\label{blowup}

Consider the $2n$-parameter family of perturbation of a map
$f:I \hookrightarrow I$ by polynomials of degree $2n-1$
\beq \label{2ndegree}
f_{\eps}(x)=f(x)+\sum_{k=0}^{2n-1}\eps_k x^k.
\eneq 
Define a map 
\beq
\beal\label{multijet}
{{\Cal J}}^{1}_nf:\underbrace{I \times \dots \times I}_{n \ \textup{times}}
\times \R^{2n} & \to 
\underbrace{I \times \dots \times I}_{n \  \textup{times}}
\times \underbrace{(I \times \R) \times \dots 
\times (I\times \R)}_{n \ \textup{times}}\\
{{\Cal J}}^{1,n}f(x_0, \dots , x_{n-1}, \eps)= & \\
\left(x_0, \dots , x_{n-1},f_\eps(x_0),f'_\eps(x_0),\right.&
\left.\dots,f_\eps(x_{n-1}),f'_\eps(x_{n-1})\right).
\enal
\eneq
This map is called {\it the $n$-tuple 1-jet map}. The $1$-jet of 
a function means that we take into account not only the image of 
a point, but also its derivative.
The $1$-jet of a function is usually denoted by 
$j^1f_\eps(x)=(x,f_\eps(x),f_\eps'(x))$. The space of 1-jets of 
functions on the interval $I$ is denoted by ${{\Cal J}}^1(I,\R)$. 
The product of $n$ copies of ${{\Cal J}}^1(I,\R)$ is
{\it multijet space} and is denoted by
\beq
{{\Cal J}}^{1,n}(I,\R)=\underbrace{{{\Cal J}}^1(I,\R) \times \dots 
\times {{\Cal J}}^1(I,\R)}_{n \ \textup{times}}.
\eneq
We need to include into our consideration derivatives, because 
we are interested in hyperbolicity (property of derivatives) 
of periodic points. The set of points 
\beq
\Delta_n(I)=\left\{\{x_0, \dots , x_{n-1}\}\times \R^{2n} \subset 
\underbrace{I \times \dots \times I}_{n \ \textup{times}}
\times \R^{2n}:\exists\ \textup{s.t.}\ i \neq j\ x_i=x_j\right\}
\quad \eneq
is called {\it the diagonal} in the space of multijets.
In singularity theory the space of multijets is defined outside
of the diagonal $\Delta_n(I)$ and is usually denoted by
${{\Cal J}}^{1,n}(I,\R)={{\Cal J}}^{1,n}(I,\R) \setminus \Delta_n(I)$ 
(see \cite{GG}).

It is easy to see that {\it a recurrent trajectory 
$\{x_k\}_{k\in \Z_+}$ is located in a neighborhood of the diagonal  
$\Delta_n(I)$ in the space of multijets for a sufficiently large $n$}. 
If $\{x_k\}_{k=0}^{n-1}$ is a part of a recurrent trajectory of 
length $n$, then the product of distances along the trajectory  
\beq \label{productformula}
\prod_{k=0}^{n-2} \left| x_{n-1}-x_k \right|
\eneq  
measures how close $\{x_k\}_{k=0}^{n-1}$ to the diagonal 
$\Delta_n(I)$, or how independently one can perturb points of a 
trajectory. One can also say that (\ref{productformula})
is a quantitative characteristic of how recurrent a trajectory 
of length $n$ is. Introduction of this {\it product of distances
along a trajectory} is a new central point of the method.

\subsection{
Newton interpolation and blow-up
along the diagonal in multijet space}

Now look at Grigoriev-Yakovenko's construction \cite{GY} in 
the $1$-dimensional case with more details. This construction puts 
the ``Closing Lemma'' and ``Perturbation of Hyperbolicity'' 
statements above into a general framework. 

Again consider the $2n$-parameter family (\ref{2ndegree}) of 
perturbations of a $C^1$ map $f:I \to I$ by polynomials of degree $2n-1$. 
Our goal now is to describe how such perturbations affect the
$n$-tuple $1$-jet of $f$, and since the operator $j^{1,n}$ is linear
in $f$, for the time being we consider only the perturbations
$\phi_\eps$ and their $n$-tuple $1$-jets.  For each $n$-tuple
$\{x_k\}_{k=0}^{n-1}$ there is a natural transformation
${\Cal J}^{1,n} : I^n \times \R^{2n} \to {\Cal J}^{1,n}(I,\R)$ from
$\eps$-coordinates to jet-coordinates, given by
\beq
\label{epstojet}
{\Cal J}^{1,n}(x_0, \dots, x_{n-1}, \eps) =
j^{1,n} \phi_\eps(x_0, \dots, x_{n-1}).
\eneq

Instead of working directly with the transformation ${\Cal J}^{1,n}$, we
introduce intermediate $u$-coordinates based on Newton interpolation
polynomials.  The relation between $\eps$-coordinates and
$u$-coordinates is given implicitly by
\beq\label{identity}
\phi_\eps(x) = \sum_{k=0}^{2n-1}\eps_k x^k=
\sum_{k=0}^{2n-1}u_k \prod_{j=0}^{k-1} (x-x_{j(\textup{mod}\ n)}).
\eneq
Based on this identity, we can define functions
$\D^{1}_n : I^n \times \R^{2n} \to I^n \times \R^{2n}$ 
and
$\pi^{1}_n:I^n \times \R^{2n}\to {\Cal J}^{1,n}(I,\R)$
so that ${\Cal J}^{1,n} = \pi^{1}_n \circ \D^{1}_n$, or in other words 
the diagram in Figure~4.2 commutes. This definition coincides with
the one we gave before. We will show later that
$\D^{1}_n$ is invertible, while $\pi^{1}_n$ is invertible away from
the diagonal $\Delta_n(I)$ and defines a blow-up along it in the space
of multijets ${\Cal J}^{1,n}(I,\R)$. Consider diagram 4.2 for
$m=n=N=1$.

Recall that the intermediate space, denoted by 
${\cal {DD}}^{1}_n(I,\R)$, is called {\em the space of divided 
differences\/} and consists of $n$-tuples of points 
$\{x_k\}_{k=0}^{n-1}$ and $2n$ real coefficients 
$\{u_k\}_{k=0}^{2n-1}$.  Here are explicit coordinate-by-coordinate
formulas defining $\pi^{1}_{n} : {\cal {DD}}^{1}_n(I,\R) \to
{\Cal J}^{1,n}(I,\R)$.
\beq
\beal\label{Newtonexpression}
\phi_\eps(x_0) = &\, u_0, \\
\phi_\eps(x_1) = &\, u_0+u_1(x_1-x_0),\\
\phi_\eps(x_2) = &\, u_0+u_1(x_2-x_0)+u_2(x_2-x_0)(x_2-x_1),\\
\vdots\,&\\
\phi_\eps(x_{n-1}) = &\, u_0+u_1(x_{n-1}-x_0)+\dots 
+u_{n-1}(x_{n-1}-x_0)\dots (x_{n-1}-x_{n-2}),\\
\phi_\eps'(x_0) = &\, \frac{\pa}{\pa x}\left(\sum_{k=0}^{2n-1}
u_k \prod_{j=0}^{k-1} (x-x_{j(\textup{mod}\ n)})\right)\Big|_{x=x_0}, \\
\vdots\,&\\
\phi_\eps'(x_{n-1}) = &\, \frac{\pa}{\pa x} \left(\sum_{k=0}^{2n-1}
u_k \prod_{j=0}^{k-1} (x-x_{j(\textup{mod}\ n)})\right)\Big|_{x=x_{n-1}}, 
\enal
\eneq

These formulas are very useful for dynamics.  For a given base map $f$
and initial point $x_0$, the image $f_\eps(x_0) = f(x_0) +
\phi_\eps(x_0)$ of $x_0$ depends only on $u_0$.  Furthermore the
image can be set to any desired point by choosing $u_0$ appropriately
--- we say then that it depends nontrivially on $u_0$.  If
$x_0$, $x_1$, and $u_0$ are fixed, the image $f_\eps(x_1)$ of $x_1$
depends only on $u_1$, and as long as $x_0 \neq x_1$ it depends
nontrivially on $u_1$.  More generally for $0 \leq k \leq n-1$, if 
pairwise distinct points $\{x_j\}_{j=0}^k$ and coefficients
$\{u_j\}_{j=0}^{k-1}$ are fixed, then the image $f_\eps(x_k)$ of
$x_k$ depends only and nontrivially on $u_k$.

Suppose now that an $n$-tuple of points $\{x_j\}_{j=0}^{n}$ not on the
diagonal $\Delta_n(I)$ and Newton coefficients $\{u_j\}_{j=0}^{n-1}$
are fixed.  Then derivative $f'_\eps(x_0)$ at $x_0$ depends only and
nontrivially on $u_n$.  Likewise for $0 \leq k \leq n-1$, if distinct
points $\{x_j\}_{j=0}^{n}$ and Newton coefficients
$\{u_j\}_{j=0}^{n+k-1}$ are fixed, then the derivative $f'_\eps(x_k)$ 
at $x_k$ depends only and nontrivially on $u_{n+k}$. 

As Figure 4.3 illustrates, these considerations show that for any map
$f$ and any desired trajectory of distinct points with any given
derivatives along it, one can choose Newton coefficients
$\{u_k\}_{k=0}^{2n-1}$ and explicitly construct a map $f_\eps = f +
\phi_\eps$ with such a trajectory.  Thus we have shown that $\pi^{1}_n$ 
is invertible away from the diagonal $\Delta_n(I)$ and defines a blow-up
along it in the space of multijets ${\Cal J}^{1,n}(I,\R)$.

\begin{figure}[htbp]\label{NIP}
  \begin{center}
   \begin{psfrags}
     \psfrag{x_0}{$x_0$}
     \psfrag{x_1}{$x_1$}
     \psfrag{x_k}{$x_k$}
     \psfrag{x_{k+1}}{$x_{k+1}$}
     \psfrag{f_u(x_0)}{$f_u(x_0)$}
     \psfrag{f_u(x_k)}{$f_u(x_k)$}
     \psfrag{f'_u(x_0)}{$f'_u(x_0)$}
     \psfrag{f'_u(x_k)}{$f'_u(x_k)$}
     \psfrag{u_0}{$u_0$}
     \psfrag{u_k}{$u_k$}
     \psfrag{u_n}{$u_n$}
     \psfrag{u_{n+k}}{$u_{n+k}$}
     \psfrag{dots}{$ \cdots $}      
    \includegraphics[width= 4in,angle=0]{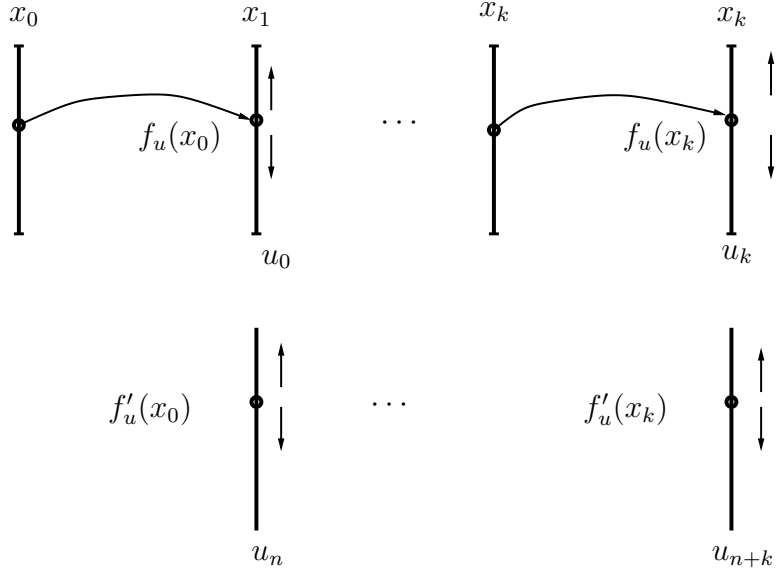}
    \end{psfrags}
   \caption{Newton coefficients and their action}
  \end{center}
\end{figure}

The function $\D^{1,n} : I^n \times \R^{2n} \to 
{\cal {DD}}^{1,n}(I,\R)$ was explicitly defined using so-called 
divided differences above.  

Recall that $\D^{1}_n$ was defined implicitly by (\ref{identity}).
We have described how to use divided differences to construct a degree
$2n-1$ interpolating polynomial of the form on the right-hand side of
(\ref{identity}) for an arbitrary $C^\infty$ function $g$.  Our
interest then is in the case $g = \phi_\eps$, which as a degree
$2n-1$ polynomial itself will have no remainder term and coincide
exactly with the interpolating polynomial.  Thus $\D^{1,n}$ is given
coordinate-by-coordinate by
\beq 
\beal\label{Newtonmap}
u_m = &\, \Delta^m \left( \sum_{k=0}^{2n-1}\eps_k x^k\right)
(x_0,\dots,x_{m\ (mod\ n)}) \\
= &\, \eps_m+\sum_{k=m+1}^{2n-1} \eps_k p_{k,m}(x_0,\dots , 
x_{m\ (mod\ n)})
\enal
\eneq
for $m = 0, \dots, 2n-1$.  We call the transformation given by
(\ref{Newtonmap}) the {\em Newton map\/}.  Notice that for fixed
$\{x_k\}_{k=0}^{2n-1}$, the Newton map is linear and given by 
an upper triangular matrix with units on the diagonal.  Hence it 
is Lebesgue volume-preserving and invertible, whether or not
$\{x_k\}_{k=0}^{2n-1}$ lies on the diagonal $\Delta_n(I)$.

We call the basis of monomials  
\beq \label{Newtonbasis}
\prod_{j=0}^k (x-x_{j(\textup{mod}\ n)}) \ \ \ 
\textup{for}\ \ \ k=0,\dots, 2n-1
\eneq
in the space of polynomials of degree $2n-1$ the {\em Newton basis\/}
defined by the $n$-tuple $\{x_k\}_{k=0}^{n-1}$.  The Newton map and
the Newton basis, and their analogues in dimension $N$, are useful
tools for perturbing trajectories and proving (\ref{growthbound}).

{\it{Acknowledgments:}}\ In this lecture I have used fragments 
of the announcement \cite{KH}. Good presentation of it  is
absolutely due to my coauthor Brian Hunt. Needless to say
that numerous communications with him and also John Mather
were very important for me.


\begin{thebibliography}{DGOP}
\def\bi#1{\bibitem[#1]{#1}}

\bi{A1} V. Arnold at el.,  Some unsolved problems in the theory
of differential equations and mathematical physics, Russ. Math.
Surveys, {\bf 44}, (1989), no. 4;  

\bi{A2} V. Arnold, Problems for Arnold's seminar, 1989;

\bi{AA} D. Anosov, V. Arnold,  Dynamical systems. I, 
Encyclopedia Math. Sci., 1, Springer, Berlin, 1988;

\bi{AGV} V. Arnold,  S. Gusein-Zade, A. Varchenko 
Singularities of differentiable maps. Vol. I. Monographs in
Mathematics, 82, Birkhauser Boston, 1985;

\bi{AM} M.\ Artin,\ B.\ Mazur,\ 
Periodic orbits, Annals of Mathematics, {\bf 81}, 1965, 82--99;

\bi{Be} I. Bendixon, Sur les corbes definies par des equations
diffirentialles, Acta Math. {\bf 24}, (1901), 1--88;

\bi{BZ} I. Berezin, N. Zhidkov, Computing Methods, Vol. 1,
Pergamon, Oxford, 1965;

\bi{BCR} J. Bochnak, M. Coste, M.-F. Roy, Real Algebraic Geometry,
Springer-Verlag, Berlin, 1998. 

\bi{Bo} A. Bolibrukh, present volume; 

\bi{Bu} A. Buium, present volume;

\bi{DeR} Z. Denkowska, R. Roussarie, A Method of Desingularization
for Analytic two-dimensional vector field families, Bol. Soc. 
Bras. Mat. {\bf 22}, 1, (1991), 93--126;

\bi{DeW} Z. Denkowska, K. Wachta, A Construction of a subanalytic
stratification under the condition (w), Bull. Polish Acad. Sci.
Math. {\bf 35}, (1987), no. 7--8, 401--405;

\bi{Dr1} L. van den Dries, present volume;

\bi{Dr2} L. van den Dries,  o-minimal structures and real analytic
geometry. Current developments in mathematics, 1998
(Cambridge, MA), 105--152, Int. Press, Somerville, MA, 1999;

\bi{Du} H. Dulac, Sur les cycles limites, Bull. Soc. Math. France
{\bf 51}, (1923), 45--188;

\bi{D} F. Dumortier, Singularities of vector fields on the plane,
J. Diff.Equations, {\bf 23}, (1977), 53--106;

\bi{DMR} F. Dumortier, M. El Morsalani, C. Rousseau, Hilbert's 16th 
problem for quadratic systems and cyclicity of elementary graphics,
Nonlinearity {\bf 9} (1996), no. 5, 1209--1261;

\bi{DRR} F. Dumortier, R. Roussarie, C. Rousseau, Hilbert's 16th
problem for quadratic vector fields. J. Differential Equations,
{\bf 110}, (1994), no. 1, 86--133;

\bi{E} J. Ecalle, Introduction aux fonctions analysables et preuve
consrustive de la conjecture de Dulac, Herman, Paris, (1992);

\bi{FP} J.-P. Francoise, C. Pugh, Keeping track of limit cycles. J. 
Diff. Eqns {\bf 65}, (1986), no. 2, 139--157;

\bi{GK} A. Gabrielov, A. Khovanskii, Multiplicity of a Notherian 
Intersection, Geometry of Differential Equations,
Amer. Math. Soc. Transl., 1998, 119-131;

\bi{Ga} L. Gavrilov, The infinitesimal 16th Hilbert problem in 
the quadratic case. Invent. Math. {\bf 143} (2001), no. 3, 449--497;

\bi{GWPL} C. Gibson, K. Wirthmuller, A. du Plessis,
E. Loojenga, Topological Stability of Smooth Mappings, LNM, {\bf 552}, 
Springer, 1976;

\bi{GG} M. Golubitsky, V. Guillemin, Stable Mappings and 
Their Singularities, Graduate Texts in Mathematics {\bf 14},
Springer-Verlag, 1973;

\bi{GST1} S. Gonchenko, L. Shil'nikov, D. Tuvaev,
On models with non-rough Poincar\'e homoclinic curves,
{\it Physica D} {\bf 62} (1993), 1--14;

\bi{GST2} S. Gonchenko, L. Shil'nikov, D. Tuvaev,
Homoclinic tangencies of an arbitrary order in Newhouse's
domains (in Russian), preprint;

\bi{GM} M. Goresky,  R. MacPherson, Stratified 
Morse Theory, Springer, 1987;

\bi{GY} A. Grigoriev, S. Yakovenko, Topology of generic multijet
preimages and blow-up via Newton interpolation, J. Diff. Eqns
150, (1998), no. 2, 349--362;

\bi{Go} M. Gromov, On entropy of holomorphic maps,
preprint;

\bi{Gr} T. Grozovskii, Bifurcations of polycycles an ``apple'' and
a ``half-apple'' in generic two-parameter families,
Diff. Equations (in Russian), {\bf{32}}, (1996), no. 4, 458--469;

\bi{GH} J. Guckenheimer, P. Holmes, Nonlinear Oscillations,
Dynamical Systems, and Bifurcations of Vector Fields, Springer-Verlag,
New York, 1983; 

\bi{Hi} H. Hironaka, Number Theory, Algebraic
Geometry and Commutative Algebra, Volume in Honor of
Y. Akizuki, Kinokunia, Tokyo, 1973;

\bi{HSY} B. Hunt, T. Sauer, J. Yorke, Prevalence: 
a translation-invariant `almost every' on infinite-dimensional 
spaces" Bull. Amer. Math. Soc. (N.S.) {\bf 27}, (1992), no. 2, 
217--238 \& {\bf 28}, (1993), no. 2, 306--307. 

\bi{I1} Yu. Ilyashenko, The multiplicity of limit cycles that arise
in the perturbation of a Hamiltonian equation of the class
$dw/dz=P_2/Q_1$, in the real and complex domain, Amer. Math.Soc.
Transl. (2), {\bf 118}, (1982), 191--202; 

\bi{I2} Yu. Ilyashenko,  Dulac's memoir "On limit cycles" and related 
questions of the local theory of differential equations, Russ. Math.
Surveys, {\bf 40}, (1985), no. 6(246), 41--78; 

\bi{I3} Yu. Ilyashenko, Finiteness theorem for limit cycles,
Amer. Math. Soc., Providence, 1991;

\bi{I4} Yu. Ilyashenko, Normal forms for local families and nonlocal 
bifurcations, Asterisque, {\bf{222}}, (1994), 233--258;


\bi{IK} Yu. Ilyashenko, V. Kaloshin, Bifurcation of planar and spatial 
polycycles: Arnold's program and its development. The Arnoldfest,
241--271, Fields Inst. Commun., 24, Amer. Math. Soc., Providence, RI,
1999;

\bi{IL} Yu. Ilyashenko, W. Li, Nonlocal bifurcations. Math. Surv. 
and Mon., 66. American Mathematical Society, Providence, RI, 1999; 

\bi{IY1} Yu.Ilyashenko, S. Yakovenko, Finitely smooth normal
forms of local families of diffeomorphisms and vector fields,
Russian Math. Surveys {\bf{46}}, (1991), no. 1, 1--43;

\bi{IY2} Yu.Ilyashenko, S. Yakovenko, Finite Cyclicity of Elementary
Polycycles in Generic Families, Amer.Math.Soc.Transl, {\bf 165},
(1995), 1--20 \& 21--95;

\bi{Ja}  N. Jakobson, Basic Algebra, vol. 1, 1974; 

\bi{Ka1} V. Kaloshin, The Hilbert-Arnold Problem and estimates on
cyclicity of elementary polycycles and multiplicity of generic germs,
preprint;

\bi{Ka2} V. Kaloshin, Bifurcations of spatial polycycles,
(in Russian) (in preparation);

\bi{Ka3} V. Kaloshin, The Hilbert-Arnold Problem and 
estimates on cyclicity of planar and spacial polycycles,
Funct. Anal. and Appl. {\bf 35}, (2001), 78--81;

\bi{Ka4} V. Kaloshin, A Geometric Proof of Existence of
Whitney's stratifications, preprint;

\bi{Ka5} V. Kaloshin, An Extension of Artin-Mazur Theorem,
Ann of Math. {\bf 50}, (1999), no. 2, 729--741; 

\bi{Ka6} V. Kaloshin, Generic diffeomorphisms with 
superexponential growth of number of periodic points, 
Comm. in Math. Physics {\bf 211}, (2000) 1, 253--271;

\bi{Ka7} V. Kaloshin, Some prevalent properties of smooth 
dynamical systems, Proc. of Steklov Math. Inst. {\bf 213}, 
(1997), 123--151; 

\bi{Ka8} V. Kaloshin, Prevalence in spaces of finitely smooth 
mappings, Funct. Anal. Appl. {\bf 31}, (1997), no. 2, 95--99;

\bi{Ka9} V. Kaloshin, Streched exponential estimate on 
the rate of growth of the number of periodic points for prevalent 
diffeomorphisms, Thesis, Princeton, 2001;

\bi{KH} V. Kaloshin, B. Hunt, Streched exponential estimate on 
the rate of growth of the number of periodic points for prevalent 
diffeomorphisms, Electronic Research Announcements of AMS, 
part I, {\bf 7}, 17--27, 2001 \& part II, {\bf 7}, 28--36, 2001;

\bi{Kh1} A. Khovanskii, Real analytic manifolds with the property
of finiteness and complex abelian integrals, Func. Anal and Appl.
{\bf 18}, (1984), no.2, 40--50;

\bi{Kh2} A. Khovanskii, Fewnomials, 
Amer.Math.Soc. Transl., Providence, RI, 1991;

\bi{Kl} O. Kleban, Order of the topologically sufficient jet of a
smooth vector field on the real plane at a singular point of finite 
multiplicity. Amer. Math. Soc. Transl, Ser.2, {\bf 165}, (1995),
131--153; 

\bi{KS} A. Kotova, V. Stanzo, Few-Parameter Generic Families on the
Sphere, Amer. Math. Soc. Translations, Providence, RI, Ser. 2,
{\bf{213}}, (1996), 155--202;

\bi{Ku} T.-C. Kuo, The ratio test for analytic Whitney
stratifications, Lecture Notes, No. 192, 141-149;

\bi{Lo1} S. Lojasiewicz, Ensemble Semi-Analytiques,
IHES Lecture Notes, 1965;

\bi{Lo2} S. Lojasiewicz, Sur le g\'eometrie semi- et 
sous-analytic, (French) Ann. Inst. Fourier (Grenoble), 
{\bf 43}, (1993), no. 5, 1575--1595. 

\bi{LSW} S. Lojasiewicz, J. Stasica, K. Wachta, 
Subanalytic stratifications. Verdier's condition, Bull. 
Polish Acad. Sci. Math. {\bf 34}, (1986), no. 9--10, 531--539;

\bi{Mr} P. Mardesi\'c, An explicit bound for the multiplicity of
zeroes of generic Abelian integrals, Nonlinearity, {\bf 4}, (1991),
no.3, 845--852;

\bi{Ma} J. Mather, Notes on topological stability, Harvard University, 
1970;

\bi{Mi1} J. Milnor, Singularities of Complex 
Hypersurfaces, Ann. of Math. Studies, no. 61, 1968;

\bi{Mi2} J. Milnor, Topology from Differentiable viewpoint,
Princeton University Press, 1997 ;

\bi{MR} R. Moussu, C. Roche, Khovanskii's theory and the Dulac
problem, Invent. Math. {\bf 105} (1991), no. 2, 431--441;

\bi{Mu} D. Mumford, Algebraic Geometry I, Spri\-n\-ger, 
New York, 1976;

\bi{NY1} D. Novikov, S. Yakovenko, Simple Exponential
Estimate for the number of zeroes of complete Abelian Integrals,
Ann. Inst. Fourier (Grenoble), {\bf 45}, (1995), no. 4,
897--927;

\bi{NY2} D. Novikov, S. Yakovenko, present volume;

\bi{O} J. Oxtoby, Measure and category. A survey of the analogies 
between topological and measure spaces. Graduate Texts in 
Mathematics, 2. Springer-Verlag, New York-Berlin, 1980. 

\bi{Pa} W. Pawlucki, The Puiseux Theorem for subanalytic 
mappings, Polish Acad. Sci. Math. {\bf 32}, (1984), no. 9--10, 
555--560;

\bi{P} G. Petrov, Nonoscillation of elliptic integrals, Func. Anal.
and Appl., {\bf 24}, (1990), no.3, 205--210;

\bi{PW} A. du Plessis, T. Wall, The Geometry
of Topological Stability, Oxford, 1995;

\bi{R1} R. Roussarie, Cyclicite finie et le 16 problem d'Hilbert,
Dynamical systems, (Volparaiso, (1986)), (R.Bauon, R.Lavarca, and 
J.Palis, eds.), LNM, {\bf 1331}, Springer-Verlag, Berlin and New York, 
(1988), 161--188;

\bi{R2} R. Roussarie, Bifurcations of Planar Vector Fields and 
Hilbert's Sixteenth Problem, Progress in Mathematics, {\bf 164},
Birkhauser, 1998;

\bi{Se} A. Seidenberg, Reduction of singularities of the 
differential equations $Ady= Bdx$, Amer. J. Math. {\bf 90},
(1968), 248--269;

\bi{Sh} L. Shilnikov, A case of the existence of a denumerable set of 
periodic motions. (Russian) Dokl. Akad. Nauk SSSR 160, (1965), 558--561; 

\bi{S} D. Shlomiuck(ed), Bifurcations and periodic orbits of vector
fields, NATO AS1, Series C (Math. and Phys. Sciences), vol. {\bf 408},
Kluwel, Dordrecht, Boston, London, 1993;

\bi{Sm} S. Smale, Differentiable Dynamical Systems,
Bull. Amer. Math.\ Soc. {\bf 73} (1967), 747--817;

\bi{Ta} F. Takens, Unfoldings of certain singularities of vector
fields: generalized Hopf bifurcations, J. Differential Equations 
{\bf{14}} (1973), 476--493;

\bi{Th1} R. Thom, Ensembles et morphismes stratifies. 
Bull. Amer. Math. Soc. {\bf{75}}, (1969), 240--284;

\bi{Th2} R. Thom, Propri\'et\'e Diff\'erentielle 
Locales des Ensembles Analytiques, Seminaire Bourbaki,
1964/65, exp. 281;

\bi{Tr} S. Trifonov, Desingularization in Families of Analytic
Differential Equations, Amer. Math. Soc. Transl, Ser.2, {\bf 165}, 
(1995), 97--129;

\bi{V} A. Varchenko, Estimation of the number of zeros of 
an abelian integral depending on a parameter, and limit cycles, 
Func. Anal and Appl. {\bf 18}, (1984), no.2, 14--25;

\bi{Wa} T. Wall,  Regular Stratifications,
Lecture Notes in Mathematics, No. 468,  332-344;

\bi{Wh} H. Whitney,  Tangents to an Analytic Variety,
Ann. of Math. {\bf 81}, (1965),  496--549.

\bi{Y} Y. Yomdin, A quantitative version of the Kupka-Smale
Theorem, Ergod. Th. Dynam. Sys. {bf 5}, (1985), 449--472.  
\end{thebibliography}
\end{document}